\documentclass[10pt]{amsart}
\usepackage[margin=2.5cm]{geometry}
\usepackage{tikz-cd}
\usepackage{graphicx}					
\usepackage{amssymb}
\usepackage{amsmath}
\usepackage{subfig, graphicx}
\usepackage{mathrsfs}
\usepackage{amssymb, amsthm, indentfirst}
\usepackage{relsize}
\usepackage{hyperref}
\usepackage{stmaryrd}
\usepackage{lineno}

\numberwithin{equation}{section}
\usepackage{color}
\theoremstyle{plain}
\newtheorem{thm}{Theorem}[section]

\newtheorem{conj}[thm]{Conjecture}
\newtheorem{lemma}[thm]{Lemma}
\newtheorem{prop}[thm]{Proposition}
\newtheorem{corollary}[thm]{Corollary}
\theoremstyle{definition}

\newtheorem{rmk}[thm]{Remark}
\def\Gal{\operatorname{Gal}}

\newtheorem*{hypothesis*}{Hypothesis}

\author{SHIH-YU CHEN}
\address{Institute of Mathematics~\\Academia Sinica~\\ 6F, Astronomy-Mathematics Building, No.\,1, Sec.\,4, Roosevelt Road, Taipei 10617, Taiwan, ROC}
\address{Department of Mathematics, Kyoto University, Kitashirakawa Oiwake-cho, Sakyo-ku, Kyoto 606-8502, Japan}
\email{sychen.math@gmail.com}

\def\GL{{\rm{GL}}}
\def\GSp{{\rm GSp}}

\def\o{\frak{o}}

\def\A{{\mathbb A}}
\def\C{{\mathbb C}}
\def\E{{\mathbb E}}
\def\F{{\mathbb F}}
\def\K{{\mathbb E}}
\def\L{{\mathbb L}}
\def\R{{\mathbb R}}
\def\Q{{\mathbb Q}}
\def\Z{{\mathbb Z}}

\def\<{\langle}
\def\>{\rangle}

\def\bp{\begin{pmatrix}}
\def\ep{\end{pmatrix}}

\def\<{\langle}
\def\>{\rangle}

\def\GL{\operatorname{GL}}

\def\GSp{\operatorname{GSp}}

\def\1{\mathbf{1}}

\def\itPi{\mathit{\Pi}}
\def\itPsi{\mathit{\Psi}}
\def\itSigma{\mathit{\Sigma}}

\title{On Deligne's conjecture for symmetric fourth $L$-functions of Hilbert modular forms}

\begin{document}
\begin{abstract}
We prove an automorphic analogue of Deligne's conjecture for symmetric fourth $L$-functions of Hilbert modular forms. We extend the result of Morimoto \cite{Morimoto2021} based on generalization and refinement of the results of Grobner and Lin \cite{GL2020} to cohomological irreducible essentially conjugate self-dual cuspidal automorphic representations of $\GL_2$ and $\GL_3$ over CM-fields.
\end{abstract}

\maketitle
\tableofcontents

\section{Introduction}

Let $\F$ be a totally real number field with $[\F:\Q]=d$.
Let $\itPi$ be a cohomological irreducible cuspidal automorphic representation of $\GL_2(\A_\F)$ with central character $\omega_\itPi$. We have $|\omega_\itPi| = |\mbox{ }|_{\A_\F}^{\sf w}$ for some ${\sf w} \in \Z$.
Let $f_\itPi$ and $f_{\itPi^\vee}$ be the normalized newforms of $\itPi$ and $\itPi^\vee$, respectively.
The Petersson norm of $f_\itPi$ is defined by
\[
\Vert f_\itPi \Vert = \int_{\A_\F^\times\GL_2(\F)\backslash \GL_2(\A_\F)} f_\itPi(g)f_{\itPi^\vee}(g\cdot{\rm diag}(-1,1)_\infty)\,dg^{\rm Tam}.
\]
Here $dg^{\rm Tam}$ is the Tamagawa measure on $\A_\F^\times \backslash \GL_2(\A_\F)$.
Note that by definition we have $\Vert f_\itPi \Vert =\Vert f_{\itPi^\vee} \Vert$.
For each archimedean place $v$ of $\F$, we have
\[
\itPi_v = D_{\kappa_v} \otimes |\mbox{ }|_\R^{{\sf w}/2}
\]
for some $\kappa_v \in \Z_{\geq 2}$ such that $\kappa_v \equiv {\sf w} \,({\rm mod}\,2)$. Here $D_\kappa$ is the discrete series representation of $\GL_2(\R)$ with weight $\kappa \in \Z_{\geq 2}$.
For a finite order Hecke character $\chi$ of $\A_\F^\times$, let 
\[
L(s,\itPi,{\rm Sym}^4 \otimes \chi)
\]
be the twisted symmetric fourth $L$-function of $\itPi$ by $\chi$. We denote by $L^{(\infty)}(s,\itPi,{\rm Sym}^4 \otimes \chi)$ the $L$-function obtained by excluding the archimedean $L$-factors. A critical point for $L(s,\itPi,{\rm Sym}^4 \otimes \chi)$ is an integer $m$ which is not a pole of the archimedean local factors $L(s,\itPi_v,{\rm Sym}^4 \otimes \chi_v)$ and $L(1-s,\itPi_v^\vee,{\rm Sym}^4 \otimes \chi_v^{-1})$ for all archimedean places $v$ of $\F$.
More precisely, $L(s,\itPi,{\rm Sym}^4 \otimes \chi)$ has no critical points if the signature of $\chi$ is not parallel. Suppose $\chi$ has parallel signature ${\rm sgn}(\chi)\in \{\pm1\}$, then the set of critical points for $L(s,\itPi,{\rm Sym}^4 \otimes \chi)$ is the disjoint union of ${\rm Crit}^+({\rm Sym}^4\itPi;\,{\rm sgn}(\chi))$ and ${\rm Crit}^-({\rm Sym}^4;\,{\rm sgn}(\chi))$, where
\begin{align*}
{\rm Crit}^+({\rm Sym}^4\itPi;\,\pm) & = \left\{1-2{\sf w} \leq  m \leq \min_{v \mid \infty}\{\kappa_v\}-1-2{\sf w}\,\left\vert\, (-1)^m=\pm\right.\right\},\\
{\rm Crit}^-({\rm Sym}^4\itPi;\,\pm) & = \left\{1-m-4{\sf w} \, \vert \, m \in {\rm Crit}^+({\rm Sym}^4\itPi;\,\pm)\right\}.
\end{align*}
For $\sigma \in {\rm Aut}(\C)$, let ${}^\sigma\!\itPi$ be the unique cohomological irreducible cuspidal automorphic representation of $\GL_2(\A_\F)$ such that ${}^\sigma\!\itPi_f$ is the $\sigma$-conjugate of $\itPi_f = \bigotimes_{v \nmid \infty}\itPi_v$.
We have the following automorphic analogue of Deligne's conjecture  \cite{Deligne1979} on the algebraicity of the critical values of $L(s,\itPi,{\rm Sym}^4 \otimes \chi)$. 
\begin{conj}\label{C:Deligne}
Let $\chi$ be a finite order Hecke character of $\A_\F^\times$ with parallel signature.
\begin{itemize}
\item[(1)] Let $m-2{\sf w} \in {\rm Crit}^+({\rm Sym}^4\itPi;\,{\rm sgn}(\chi))$. We have
\begin{align*}
&\sigma \left(\frac{L^{(\infty)}(m-2{\sf w},\itPi,{\rm Sym}^4 \otimes \chi)}{|D_\F|^{1/2}\cdot(2\pi\sqrt{-1})^{3dm+3\sum_{v \mid \infty}\kappa_v}\cdot (\sqrt{-1})^{d{\sf w}}\cdot G(\chi\cdot\omega_\itPi^2)^3\cdot \Vert f_\itPi \Vert^3} \right)\\
& = \frac{L^{(\infty)}(m-2{\sf w},{}^\sigma\!\itPi,{\rm Sym}^4 \otimes {}^\sigma\!\chi)}{|D_\F|^{1/2}\cdot(2\pi\sqrt{-1})^{3dm+3\sum_{v \mid \infty}\kappa_v}\cdot (\sqrt{-1})^{d{\sf w}}\cdot G({}^\sigma\!\chi\cdot{}^\sigma\!\omega_\itPi^2)^3\cdot \Vert f_{{}^\sigma\!\itPi} \Vert^3}
\end{align*}
for all $\sigma \in {\rm Aut}(\C)$.
\item[(2)] Let $m-2{\sf w} \in {\rm Crit}^-({\rm Sym}^4\itPi;\,{\rm sgn}(\chi))$. We have
\begin{align*}
&\sigma \left(\frac{L^{(\infty)}(m-2{\sf w},\itPi,{\rm Sym}^4 \otimes \chi)}{(2\pi\sqrt{-1})^{2dm+3\sum_{v \mid \infty}\kappa_v}\cdot (\sqrt{-1})^{d{\sf w}}\cdot G(\chi\cdot\omega_\itPi^2)^2\cdot \Vert f_\itPi \Vert^3} \right)\\
& = \frac{L^{(\infty)}(m-2{\sf w},{}^\sigma\!\itPi,{\rm Sym}^4 \otimes {}^\sigma\!\chi)}{(2\pi\sqrt{-1})^{2dm+3\sum_{v \mid \infty}\kappa_v}\cdot (\sqrt{-1})^{d{\sf w}}\cdot G({}^\sigma\!\chi\cdot{}^\sigma\!\omega_\itPi^2)^2\cdot \Vert f_{{}^\sigma\!\itPi} \Vert^3}
\end{align*}
for all $\sigma \in {\rm Aut}(\C)$.
\end{itemize}
Here $D_\F$ is the discriminant of $\F$ and $G(\chi\cdot\omega_\itPi^2)$ is the Gauss sum of $\chi\cdot\omega_\itPi^2$.
\end{conj}

\begin{rmk}
The conjecture is true when $\itPi$ is of CM-type (cf.\,\cite[\S\,4]{RS2007c}).
\end{rmk}

\begin{rmk}
For critical points in ${\rm Crit}^+({\rm Sym}^4\itPi;\,{\rm sgn}(\chi))$, the algebraicity of the critical values was proved by Pitale--Saha--Schmidt \cite[\S\,7.4]{PSS2020} in terms of Petersson norm of the Ramakrishnan--Shahidi lift of $\itPi$ (cf.\,\cite{RS2007b}) under certain assumptions. These assumptions were lifted by the recent work \cite{Liu2019b} of Liu. Similar result was also proved by Ibukiyama and Katsurada \cite{IK2014} for elliptic modular forms of full level.
\end{rmk}

Following is the first main result of this paper.
\begin{thm}\label{T:Deligne}
If $\kappa_v \geq 3$ for all archimedean places $v$ of $\F$, then Conjecture \ref{C:Deligne} holds.
\end{thm}

Conjecture \ref{C:Deligne}-(1) was proved by Morimoto \cite{Morimoto2021} under additional assumptions that $\omega_\itPi$ is trivial and $\sigma \in  {\rm Aut}(\C/\mathbb{K})$ for some totally imaginary number field $\mathbb{K}$ containing $\F$. Morimoto proved the algebraicity based mainly on the result \cite[Theorem C]{GL2020} of Grobner and Lin in the special case $\GL_3 \times \GL_2$, and together with other results on the algebraicity of critical $L$-values including \cite{Shimura1978}, \cite{Im1991}, \cite{GH1993}, \cite{Liu2019b}, and \cite{PSS2020}.
We prove Theorem \ref{T:Deligne} by following similar idea. The novelty here is that we generalize and refine the results \cite[Theorems A, B, and C]{GL2020} of Grobner and Lin for $\GL_3 \times \GL_2$ (cf.\,Theorems \ref{T:RS 0}, \ref{T:Asai 0}, and Corollary \ref{T:GGP 2} respectively). For instance, we extend \cite[Theorem C]{GL2020} to cohomological irreducible essentially conjugate self-dual cuspidal automorphic representations of $\GL_2$ and $\GL_3$ over CM-fields and prove the Galois-equivariance property in full generality.
We shall introduce our second main result Theorems \ref{T:RS 0} and \ref{T:Asai 0} in the next paragraph.

Let $\K$ be a totally imaginary quadratic extension of $\F$.
Let $\omega_{\K/\F}$ be the quadratic Hecke character of $\A_\F^\times$ associated to $\K/\F$ by class field theory.
Let $\itPi$ and $\itSigma$ be cohomological irreducible cuspidal automorphic representations of $\GL_2(\A_\K)$ and $\GL_3(\A_\K)$ with central characters $\omega_\itPi$ and $\omega_\itSigma$, respectively.
We have $|\omega_\itPi| = |\mbox{ }|_{\A_\E}^{{\sf w}(\itPi)}$ and $|\omega_\itSigma| = |\mbox{ }|_{\A_\E}^{3{\sf w}(\itSigma)/2}$ for some ${\sf w}(\itPi),{\sf w}(\itSigma) \in \Z$.
For each archimedean place $v$ of $\E$, we have
\[
\itPi_v = {\rm Ind}_{B_2(\C)}^{\GL_2(\C)}(\chi_{\kappa_{1,v}}\boxtimes \chi_{\kappa_{2,v}})\otimes |\mbox{ }|_\C^{{\sf w}(\itPi)/2},\quad \itSigma_v = {\rm Ind}_{B_3(\C)}^{\GL_3(\C)}(\chi_{\ell_{1,v}} \boxtimes \chi_{\ell_{2,v}} \boxtimes \chi_{\ell_{3,v}})\otimes |\mbox{ }|_\C^{{\sf w}(\itSigma)/2}
\]
for some $\underline{\kappa}_v = (\kappa_{1,v},\kappa_{2,v}) \in \Z^2$ and $\underline{\ell}_v = (\ell_{1,v},\ell_{2,v},\ell_{3,v}) \in \Z^3$ such that
\begin{align*}
\kappa_{1,v} > \kappa_{2,v},\quad &\kappa_{1,v} \equiv \kappa_{2,v}\equiv 1+ {\sf w}(\itPi)\,({\rm mod}\,2),\\
\ell_{1,v} > \ell_{2,v} > \ell_{3,v},\quad &\ell_{1,v} \equiv \ell_{2,v} \equiv \ell_{3,v} \equiv {\sf w}(\itSigma) \,({\rm mod}\,2).
\end{align*}
Here $\chi_\kappa$ is the character of $\C^\times$ defined by $\chi_\kappa(z) = (z/\overline{z})^{\kappa/2}$ for $\kappa \in \Z$.
Let 
\[
L(s,\itSigma \times \itPi)
\] be the Rankin--Selberg $L$-function of $\itSigma \times \itPi$ and denote by $L^{(\infty)}(s,\itSigma \times \itPi)$ the $L$-function obtained by excluding the archimedean $L$-factors. 
A critical point for $L(s,\itSigma \times \itPi)$ is a half-integer $m+\tfrac{1}{2}$ which is not a pole of the archimedean local factors $L(s,\itSigma_v \times \itPi_v)$ and $L(1-s,\itSigma^\vee_v \times \itPi^\vee_v)$ for all archimedean places $v$ of $\E$.
For $\sigma \in {\rm Aut}(\C)$, let ${}^\sigma\!\itPi$ and ${}^\sigma\!\itSigma$ be the unique cohomological irreducible cuspidal automorphic representations of $\GL_2(\A_\E)$ and $\GL_3(\A_\E)$ such that ${}^\sigma\!\itPi_f$ and ${}^\sigma\!\itSigma_f$ are the $\sigma$-conjugate of $\itPi_f = \bigotimes_{v \nmid \infty}\itPi_v$ and $\itSigma_f = \bigotimes_{v \nmid \infty}\itSigma_v$, respectively.
With respect to fixed choice of generators in the relative Lie algebra cohomology groups for $\GL_3(\C)$ and $\GL_2(\C)$ (cf.\,(\ref{E:generator GL_2}) and (\ref{E:generator GL_3})), we can define the bottom degree Whittaker periods 
\[
p^b(\itSigma),\quad p^b(\itPi)
\] 
of $\itSigma$ and $\itPi$ respectively. These periods are defined by comparing the rational structures on $\itSigma_f$ and $\itPi_f$ given by the Whittaker models and by the isotypic components in the bottom degree cuspidal cohomology of locally symmetric spaces with appropriate coefficients.
Following is our result on the algebraicity of the critical values of $L(s,\itSigma \times \itPi)$ in terms of the bottom degree Whittaker periods.

\begin{thm}[Theorem \ref{T:RS}]\label{T:RS 0}
Assume $\ell_{1,v} > -\kappa_{2,v} > \ell_{2,v} > -\kappa_{1,v} > \ell_{3,v}$ for all $v \in S_\infty$.
Let $m+\tfrac{1}{2}$ be critical for $L(s,\itSigma \times \itPi)$.
We have
\begingroup
\begin{align*}
&\sigma\left(\frac{L^{(\infty)}(m+\tfrac{1}{2},\itSigma \times \itPi)}{(2\pi\sqrt{-1})^{6dm}\cdot G(\omega_{\E/\F})\cdot G(\omega_\itPi)\cdot\pi^{3d(1+{\sf w}(\itSigma)+{\sf w}(\itPi))+\sum_{v \in S_\infty}(\ell_{1,v}-\ell_{3,v}+(\kappa_{1,v}-\kappa_{2,v})/2)}\cdot  p^b(\itSigma)\cdot p^b(\itPi)} \right)\\
&= \frac{L^{(\infty)}(m+\tfrac{1}{2},{}^\sigma\!\itSigma \times {}^\sigma\!\itPi)}{(2\pi\sqrt{-1})^{6dm}\cdot G(\omega_{\E/\F})\cdot G({}^\sigma\!\omega_\itPi)\cdot\pi^{3d(1+{\sf w}(\itSigma)+{\sf w}(\itPi))+\sum_{v \in S_\infty}(\ell_{1,v}-\ell_{3,v}+(\kappa_{1,v}-\kappa_{2,v})/2)}\cdot p^b({}^\sigma\!\itSigma)\cdot p^b({}^\sigma\!\itPi)}
\end{align*}
\endgroup
for all $\sigma \in {\rm Aut}(\C)$. 
\end{thm}


We assume further that $\itPi$ and $\itSigma$ are essentially conjugate self-dual, that is, 
\[
\itPi^\vee = \itPi^c \otimes \chi\circ {\rm N}_{\K/\F},\quad \itSigma^\vee = \itSigma^c \otimes \eta\circ {\rm N}_{\K/\F}
\]
for some algebraic Hecke characters $\chi$ and $\eta$ of $\A_\F^\times$ with parallel signatures ${\rm sgn}(\chi)$ and ${\rm sgn}(\eta)$, respectively.
After replacing $\chi$ by $\omega_{\E/\F}\cdot\chi$ and $\eta$ by $\omega_{\E/\F}\cdot\eta$ if necessary, we may assume that
\[
{\rm sgn}(\chi) = (-1)^{{\sf w}(\itPi)},\quad {\rm sgn}(\eta) = (-1)^{{\sf w}(\itSigma)}.
\]
Let 
\[
L(s,\itPi,{\rm As}^\pm \otimes \chi),\quad L(s,\itSigma,{\rm As}^\pm \otimes \eta)
\] be the twisted Asai $L$-functions of $\itPi$ and $\itSigma$ by $\chi$ and $\eta$, respectively. 
Note that these $L$-functions are holomorphic and non-vanishing at $s=1$ by \cite[Corollary 2.5.9]{Mok2015}.
Following is our result on the algebraicity of twisted Asai $L$-functions at $s=1$ in terms of the bottom degree Whittaker periods.

\begin{thm}[Theorem \ref{T:Asai 2}]\label{T:Asai 0}
Assume $\itPi$ and $\itSigma$ are essentially conjugate self-dual.
\begin{itemize}
\item[(1)] We have
\begin{align*}
&\sigma\left(\frac{L^{(\infty)}(1,\itPi,{\rm As}^+ \otimes \chi)}{|D_\F|^{1/2}\cdot (\sqrt{-1})^{d{\sf w}(\itPi)}\cdot  G(\chi)\cdot\pi^{3d+\sum_{v \in S_\infty}(\kappa_{1,v}-\kappa_{2,v})/2}\cdot  p^b(\itPi)}\right) \\
& = \frac{L^{(\infty)}(1,{}^\sigma\!\itPi,{\rm As}^+ \otimes {}^\sigma\!\chi)}{|D_\F|^{1/2}\cdot (\sqrt{-1})^{d{\sf w}(\itPi)}\cdot G({}^\sigma\!\chi)\cdot\pi^{3d+\sum_{v \in S_\infty}(\kappa_{1,v}-\kappa_{2,v})/2}\cdot  p^b({}^\sigma\!\itPi)}
\end{align*}
for all $\sigma \in {\rm Aut}(\C)$.
\item[(2)] We have
\begin{align*}
&\sigma\left(\frac{L^{(\infty)}(1, \itSigma,{\rm As}^- \otimes \eta)}{(\sqrt{-1})^{d{\sf w}(\itSigma)}\cdot G(\omega_{\E/\F})\cdot G(\eta)^3\cdot \pi^{6d+\sum_{v \in S_\infty}(\ell_{1,v}-\ell_{3,v})}\cdot p^b(\itSigma)} \right)\\
& = \frac{L^{(\infty)}(1, \itSigma,{\rm As}^- \otimes \eta)}{(\sqrt{-1})^{d{\sf w}(\itSigma)}\cdot G(\omega_{\E/\F})\cdot G({}^\sigma\!\eta)^3\cdot\pi^{6d+\sum_{v \in S_\infty}(\ell_{1,v}-\ell_{3,v})}\cdot p^b({}^\sigma\!\itSigma)}
\end{align*}
for all $\sigma \in {\rm Aut}(\C)$.
\end{itemize}
\end{thm}

In \cite{Raghuram2016} and \cite{GHL2016}, we have algebraicity results for critical values of Rankin--Selberg $L$-functions for $\GL_n(\A_\E) \times \GL_{n-1}(\A_\E)$ and Asai $L$-functions for $\GL_n(\A_\E)$ in terms of Whittaker periods. In these literatures, the choice of generators in the relative Lie algebra cohomology are not specified. Therefore, there are undetermined archimedean local factors appearing in the denominator of the ratios.
In \cite[Theorems A and B]{GL2020}, the problem was resolved by Grobner and Lin under some assumptions.
More precisely, the Galois-equivariant property was proved over the Galois closure of $\E$, and some non-vanishing hypotheses on central $L$-values were imposed for the central critical point.
They also assume the archimedean local components are conjugate self-dual.
We lifted these restrictions in the case $\GL_3 \times \GL_2$. Our proof of Theorems \ref{T:RS 0} and \ref{T:Asai 0} are different from theirs for their method relies on the theory of Eisenstein cohomology 
and a detailed comparison between critical $L$-values associated to automorphically induced representations and CM-periods. Our approach is based on the classical theory of cuspidal cohomology and explicit computations of the archimedean local factors appearing in the algebraicity results \cite{Raghuram2016}, \cite{GHL2016}, and \cite{BR2017}.
It would be interesting to compare the generators of relative Lie algebra cohomology defined in this paper and in \cite{GL2020}.

This paper is organized as follows. In \S\,\ref{S:proof}, we proof Theorem \ref{T:Deligne} assuming the validity of Theorems \ref{T:RS 0} and \ref{T:Asai 0}. In \S\,\ref{S:Whittaker periods}, we define the Whittaker periods for cohomological representations of $\GL_2$ and $\GL_3$ over CM-fields. The classical results on cuspidal cohomology are recalled in \S\,\ref{SS:cuspidal cohomology}. The Whittaker periods depend on the normalization of archimedean Whittaker functions in the minimal types which we specified in \S\,\ref{SS:Whittaker periods}. 
The main results of \S\,\ref{S:RS}, \S\,\ref{S:adjoint}, and \S\,\ref{S:Asai} are Theorems \ref{T:RS}, \ref{T:adjoint}, and \ref{T:Asai}, which are refinements of the results \cite{Raghuram2016}, \cite{BR2017}, and \cite{GHL2016}, respectively. In these theorems, we express the algebraicity of various special $L$-values in terms of Whittaker periods. 
In appendix \ref{S:appendix}, we proved period relations for the conjectural periods for $\GL_n$ under duality. The result is used to pass from the right-half critical range to the left-half critical range in the proof of Theorem \ref{T:Deligne}.
As other consequences of Theorem \ref{T:RS}, we prove the algebraicity of ratios of critical $L$-values in Corollary \ref{C:ratio}, and establish period relations for Whittaker periods under base change in Theorem \ref{T:period relation}.

\section{Notation and convention}
\subsection{Notation}
Let $\GL_n$ be the general linear group of rank $n$ over $\Q$. Let $N_n$ and $T_n$ be the maximal unipotent subgroup and the maximal torus of $\GL_n$ consisting of upper triangular unipotent matrices and diagonal matrices, respectively.
For $a_1,\cdots,a_n \in \GL_1$, we write
\[
{\rm diag}(a_1,\cdots,a_n) = \bp a_1 &\cdots &0\\ \vdots&\ddots&\vdots \\ 0&\cdots& a_n \ep \in T_n.
\]
Denote by $X^+(T_n)$ the set of dominant integral weights of $T_n$ with respect to the standard Borel subgroup $B_n$ consisting of upper triangular matrices. We regard $X^+(T_n)$ as the set of $n$-tuples of integers $\mu = (\mu_1,\cdots,\mu_n)$ such that $\mu_1 \geq \cdots \geq \mu_n$. 
For $\mu \in X^+(T_n)$, its dual weight $\mu^\vee \in X^+(T_n)$ is defined by $\mu^\vee_i = -\mu_{n-i+1}$ for $1 \leq i \leq n$.

Let $\F$ be a number field. 
Let $D_\F$ and ${\rm Reg}_\F$ be the discriminant and regulator of $\F$, respectively.
Let $\zeta_\F(s)$ be the completed Dedekind zeta function of $\F$.
For each place $v$ of $\F$, let $\F_v$ be the completion of $\F$ at $v$ and $\frak{o}_{\F_v}$ the ring of integers of $\F_v$ when $v$ is finite.
Let $\A_\F$ be the ring of adeles of $\F$ and $\A_{\F,f}$ its finite part. 
Let $\psi_\Q = \bigotimes_v \psi_{\Q_v} : \Q\backslash\A_\Q \rightarrow \C$ be the additive character defined so that
\begin{align*}
\psi_{\Q_p}(x) & = e^{-2\pi \sqrt{-1}x} \mbox{ for }x \in \Z[p^{-1}],\\
\psi_{\R}(x) & = e^{2\pi \sqrt{-1}x} \mbox{ for }x \in \R.
\end{align*}
Let $\psi_{n,\F}=\bigotimes_v \psi_{n,\F_v} : N_n(\F)\backslash N_n(\A_\F) \rightarrow \C$ be the non-degenerated additive character defined by
\[
\psi_{n,\F}(u) = \psi_\Q\circ {\rm tr}_{\F/\Q}(u_{12}+u_{23}+\cdots u_{n-1,n})
\]
for $u=(u_{ij}) \in N_n(\A_\F)$.
We write $\psi_\F = \psi_{1,\F}$.
Let 
\[
\Gamma_\R(s) = \pi^{-s/2}\Gamma\left(\tfrac{s}{2}\right),\quad \Gamma_\C(s) = 2(2\pi)^{-s}\Gamma(s),
\]
where $\Gamma(s)$ is the gamma function.
We recall Barne's first and second lemmas, which will be used in our computation of archimedean integrals.
\begin{lemma}
Let $a,b,c,d,e \in \C$.
\begin{itemize}
\item[(1)] We have
\begin{align*}
\int_L \frac{ds}{2\pi\sqrt{-1}}\,\Gamma_\C(s+a)\Gamma_\C(s+b)\Gamma_\C(-s+c)\Gamma_\C(-s+d) & = 2 \cdot \frac{\Gamma_\C(a+c)\Gamma_\C(a+d)\Gamma_\C(b+c)\Gamma_\C(b+d)}{\Gamma_\C(a+b+c+d)}.
\end{align*}
Here $L$ is any vertical path from south to north which keeps the poles of $\Gamma_\C(s+a)\Gamma_\C(s+b)$ and $\Gamma_\C(-s+c)\Gamma_\C(-s+d)$ on its left and right, respectively.
\item[(2)] We have
\begin{align*}
&\int_L \frac{ds}{2\pi\sqrt{-1}}\,\frac{\Gamma_\C(s+a)\Gamma_\C(s+b)\Gamma_\C(s+c)\Gamma_\C(-s+d)\Gamma_\C(-s+e)}{\Gamma_\C(s+a+b+c+d+e)} \\
& =2\cdot \frac{\Gamma_\C(a+d)\Gamma_\C(a+e)\Gamma_\C(b+d)\Gamma_\C(b+e)\Gamma_\C(c+d)\Gamma_\C(c+e)}
{\Gamma_\C(b+c+d+e)\Gamma_\C(a+c+d+e)\Gamma_\C(a+b+d+e)}.
\end{align*}
Here $L$ is any vertical path from south to north which keeps the poles of $\Gamma_\C(s+a)\Gamma_\C(s+b)\Gamma_\C(s+c)$ and $\Gamma_\C(-s+d)\Gamma_\C(-s+e)$ on its left and right, respectively.
\end{itemize}
\end{lemma}

\subsection{Lie algebras}

Let $K_n$ be the compact modulo center subgroup of $\GL_n(\C)$ defined by
\[
K_{n} = \C^\times \cdot {\rm U}(n).
\]
Here we regard $\C^\times$ as the center of $\GL_n(\C)$. 
We denote by  $\frak{su}(n)$, $\frak{u}(n)$, $\frak{k}_n$, and $\frak{g}_n$ the Lie algebras of the real Lie groups ${\rm SU}(n)$, ${\rm U}(n)$, $K_n$, and $\GL_n(\C)$, respectively. Their complexifications are denoted by $\frak{su}(n)_\C$, $\frak{u}(n)_\C$, $\frak{k}_{n,\C}$, and $\frak{g}_{n,\C}$, respectively. 
In $\frak{g}_n$, let $e_{ij}$ be the $n$ by $n$ matrix with $1$ in the $(i,j)$-entry and zeros otherwise.
By abuse of notation, the image of $X \in \frak{g}_n$ in $\frak{g}_n / \frak{k}_n$ is denoted by the same symbol $X$.
For $n=2$, let $\{Y_{(0,0)},Y_{\pm(1,-1)}\} \subset \frak{g}_{2,\C}/\frak{k}_{2,\C}$ be a basis defined by
\begin{align*}
Y_{(0,0)} = e_{11}-e_{22},\quad Y_{\pm(1,-1)} = \frac{1}{2}\left[(\sqrt{-1}\,e_{12}-\sqrt{-1}\,e_{21})\otimes\sqrt{-1}\pm (-e_{12}-e_{21})\right],
\end{align*}
and $\{Y_{(0,0)}^*,Y_{\pm(1,-1)}^*\} \subset (\frak{g}_{2,\C}/\frak{k}_{2,\C})^*$ be the corresponding dual basis.
Here the subscripts $(0,0)$ and $\pm(1,-1)$ refer to the weights under the adjoint action of ${\rm U}(1) \times {\rm U}(1)$ on $\frak{g}_{2,\C}/\frak{k}_{2,\C}$.
We denote by
\[
\left(\frak{g}_{2,\C}/\frak{k}_{2,\C}\right)^*_\Q 
\]
the $\Q$-vector space spanned by $\{Y_{(0,0)}^*,Y_{\pm(1,-1)}^*\}$.
Let $E_\pm \in \frak{su}(2)_\C$ defined by
\[
E_\pm = \frac{1}{2}\left[ -(\sqrt{-1}\,e_{12}+\sqrt{-1}\,e_{21})\otimes\sqrt{-1} \pm (e_{12} - e_{21}) \right].
\]
For $n=3$, let 
\[
\{Z_{12},Z_{23},X_{\pm(1,-1,0)},X_{\pm(0,1,-1)},X_{\pm(1,0,-1)}\} \subset \frak{g}_{3,\C}/\frak{k}_{3,\C}
\]
be a basis defined by
\begin{align*}
Z_{12}& = e_{11}-e_{22},\quad Z_{23} = e_{22}-e_{33},\\
X_{\pm(1,-1,0)} &=\frac{1}{2}\left[(\sqrt{-1}\,e_{12}-\sqrt{-1}\,e_{21})\otimes\sqrt{-1}\pm (-e_{12}-e_{21})\right],\\
X_{\pm(0,1,-1)} &=\frac{1}{2}\left[(\sqrt{-1}\,e_{23}-\sqrt{-1}\,e_{32})\otimes\sqrt{-1}\pm (-e_{23}-e_{32})\right],\\
X_{\pm(1,0,-1)} &= \frac{1}{2}\left[(\sqrt{-1}\,e_{13}-\sqrt{-1}\,e_{31})\otimes\sqrt{-1}\pm (-e_{13}-e_{31})\right],
\end{align*}
and $\{Z_{12}^*,Z_{23}^*,X_{\pm(1,-1,0)}^*,X_{\pm(0,1,-1)}^*,X_{\pm(1,0,-1)}^*\} \subset (\frak{g}_{3,\C}/\frak{k}_{3,\C})^*$ be the corresponding dual basis. 
Here the subscripts $\pm(1,-1,0)$, $\pm(0,1,-1)$, and $\pm(1,0,-1)$ refer to the weights under the adjoint action of ${\rm U}(1) \times {\rm U}(1) \times {\rm U}(1)$ on $\frak{g}_{3,\C}/\frak{k}_{3,\C}$.
We denote by
\[
\left(\frak{g}_{3,\C}/\frak{k}_{3,\C}\right)^*_\Q 
\]
the $\Q$-vector space spanned by $\{Z_{12}^*,Z_{23}^*,X_{\pm(1,-1,0)}^*,X_{\pm(0,1,-1)}^*,X_{\pm(1,0,-1)}^*\}$.
For $1 \leq i,j \leq 3$, let $E_{ij} \in \frak{su}(3)_\C$ defined by
\[
E_{ij} = \frac{1}{2}\left[  -(\sqrt{-1}\,e_{ij}+\sqrt{-1}\,e_{ji}) \otimes \sqrt{-1} + (e_{ij}-e_{ji}) \right].
\]

\subsection{Measures}\label{SS:measure}
Let $\F$ be a number field.
Let $v$ be a place of $\F$. The Haar measure $dx_v$ on $\F_v$ is normalized so that ${\rm vol}(\o_{\F_v},dx_v)=1$ if $v$ is finite, and $dx_v$ is (resp.\,twice) the Lebesgue measure on $\F_v$ if $\F_v = \R$ (resp.\,$\F_v=\C$). The Haar measure $d^\times x_v$ on $\F_v^\times$ is normalized so that ${\rm vol}(\o_{\F_v}^\times,d^\times x_v)=1$ if $v$ is finite, and 
\[
d^\times x_v = \Gamma_{\F_v}(1)\cdot \frac{dx_v}{|x_v|_{\F_v}}
\]
if $v$ is archimedean.
Here $|\mbox{ }|_\R = |\mbox{ }|$ is the usual absolute value and $|z|_\C = |z\overline{z}|_\R$.
The Haar measure $dg = \prod_{v}dg_v$ on $\GL_n(\A_\F)$ is defined as follows: When $v$ is finite, we assume ${\rm vol}(\GL_n(\o_{\F_v}),dg_v)=1$. When $v$ is archimedean, we have
\begin{align*}
\int_{\GL_n(\F_v)}f(g_v)\,dg_v = \int_{N_n(\F_v)}du\int_{\R_{+}^n}d^\times a_1\cdots d^\times a_n\int_{K_v}dk\,f(u\cdot{\rm diag}(a_1\cdots a_n,a_2\cdots a_n,\cdots,a_n)\cdot k)
\end{align*}
for $f \in L^1(\GL_n(\F_v))$. Here $du$ is defined by the product measure on $\F_v^{n(n-1)/2}$ and ${\rm vol}(K_v,dk)=1$ with
\[
K_v = \begin{cases}
{\rm O}(n) & \mbox{ if $\F_v=\R$},\\
{\rm U}(n) & \mbox{ if $\F_v=\C$}.
\end{cases}
\]

\subsection{Gauss sums}
Let $\F$ be a number field.
Let $\chi$ be an algebraic Hecke character of $\A_\F^\times$.
The signature of $\chi$ at a real place $v$ is the value $\chi_v(-1) \in \{\pm1\}$. The signature ${\rm sgn}(\chi)$ of $\chi$ is the sequence of signs $(\chi_v(-1))_{v \mid \infty,\,v\,{\rm real}}$ indexed by real places.
We say $\chi$ has parallel signature if it has the same signature at all real places.
The Gauss sum $G(\chi)$ of $\chi$ is defined by
\[
G(\chi) = |D_\F|^{-1/2}\prod_{v \nmid \infty}\varepsilon(0,\chi_v,\psi_{\F_v}),
\]
where $\varepsilon(s,\chi_v,\psi_{\F_v})$ is the $\varepsilon$-factor of $\chi_v$ with respect to $\psi_{\F_v}$ defined in \cite{Tate1979}.
For $\sigma \in {\rm Aut}(\C)$, let ${}^\sigma\!\chi$ be the unique algebraic Hecke character  of $\A_\F^\times$ such that ${}^\sigma\!\chi(x) = \sigma(\chi(x))$ for $x \in \A_{\F,f}^\times$.
It is easy to verify that 
\begin{align}\label{E:Galois Gauss sum}
\begin{split}
\sigma(G(\chi)) &= {}^\sigma\!\chi(u_\sigma)G({}^\sigma\!\chi),\\
\sigma\left(\frac{G(\chi\chi')}{G(\chi)G(\chi')}\right) &= \frac{G({}^\sigma\!\chi{}^\sigma\!\chi')}{G({}^\sigma\!\chi)G({}^\sigma\!\chi')}
\end{split}
\end{align}
for algebraic Hecke characters $\chi,\chi'$ of $\A_\F^\times$,
where $u_\sigma \in \prod_p \Z_p^\times \subset \A_{\F,f}^\times$ is the unique element such that $\sigma(\psi_{\Q}(x)) = \psi_{\Q}(u_\sigma x)$ for $x \in \A_{\Q,f}$.

\subsection{CM-fields}\label{SS:CM-fields}

Throughout this paper, we fix a totally real number field $\F$ with $[\F:\Q]=d$. Let $S_\infty$ be the set of archimedean places of $\F$. 
For $v \in S_\infty$, we identify $\F_v$ with $\R$ via the real embedding of $\F$ associated to $v$.
Let $\E$ be a totally imaginary quadratic extension of $\F$.
We identify $S_\infty$ with the set of archimedean places of $\E$ in a natural way.
Let $\omega_{\E/\F}$ be the quadratic Hecke character of $\A_\F^\times$ associated to $\E/\F$ by class field theory.
For $v \in S_\infty$, let $\{\iota_v,\overline{\iota}_v\}$ be the pair of complex embeddings of $\E$ corresponding to $v$, and we identify $\E_v$ with $\C$ via $\iota_v$.
We fix a CM-type $\{\iota_v\,\vert\,v \in S_\infty\}$ of $\E$.
For $v \in S_\infty$ and $\sigma \in {\rm Aut}(\C)$, let $\sigma\circ v \in S_\infty$ be the place associated to the embedding $\sigma \circ \iota_v$.
For $n \in \Z_{\geq 1}$, let
\[
\frak{g}_{n,\infty} = \bigoplus_{v \in S_\infty} \frak{g}_n,\quad \frak{k}_{n,\infty} = \bigoplus_{v \in S_\infty} \frak{k}_n
\]
be the Lie algebras of the real Lie groups
\[
\GL_n(\E_\infty) = \prod_{v \in S_\infty}\GL_n(\E_v),\quad K_{n,\infty} = \prod_{v \in S_\infty} K_n,
\]
respectively.
We recall in the following lemma a Galois-equivariant property of the quadratic Gauss sum $G(\omega_{\E/\F})$.

\begin{lemma}\label{L:Gauss sum}
We have
\[
\frac{\sigma (G(\omega_{\E/\F}))}{G(\omega_{\E/\F})} = \frac{\sigma (|D_\E|^{1/2}\cdot (\sqrt{-1})^d)}{|D_\E|^{1/2}\cdot (\sqrt{-1})^d} = (-1)^{{}^\sharp\{v \in S_\infty\,\vert\,\sigma \circ \iota_v = \overline{\iota}_{\sigma\circ v}\}}
\]
for all $\sigma \in {\rm Aut}(\C)$.
\end{lemma}

\begin{proof}
The first equality is well-known. Indeed, it follows from the class number formulas for $\F$ and $\E$, and (cf.\,\cite[Proposition 3.1]{Shimura1978})
\[
\frac{L^{(\infty)}(1,\omega_{\E/\F})}{|D_\F|^{1/2}\cdot (2\pi\sqrt{-1})^d\cdot G(\omega_{\E/\F})} \in \Q^\times.
\]
To prove the second equality, fix $\alpha\in\E^\times$ such that ${\rm tr}_{\E/\F}(\alpha)=0$. Let $D_{\E/\F}$ be the relative discriminant of $\E/\F$. Note that
\[
D_\E = N_{\F/\Q}D_{\E/\F}\cdot D_\F^2.
\]
Thus
\[
|D_\E| = |N_{\F/\Q}D_{\E/\F}|\cdot |D_\F|^2 \in |N_{\E/\Q}(\alpha)|\cdot \Q^{\times,2} = \prod_{v \in S_\infty}|\iota_v(\alpha)|^2 \cdot \Q^{\times,2}.
\]
On the other hand, we have
\begin{align*}
\left( \frac{(\sqrt{-1})^d\cdot \prod_{v \in S_\infty}|\iota_v(\alpha)|}{\prod_{v \in S_\infty}\iota_v(\alpha)}\right)^2 = \frac{(-1)^d\cdot \prod_{v \in S_\infty}\iota_v(\alpha) \cdot \prod_{v \in S_\infty}\overline{\iota_v(\alpha)}}{\prod_{v \in S_\infty}\iota_v(\alpha)^2} = 1.
\end{align*}
Therefore,
\[
|D_\E|^{1/2}\cdot (\sqrt{-1})^d \in \left(\prod_{v \in S_\infty}\iota_v(\alpha) \right)\cdot \Q^\times.
\]
We conclude that
\begin{align*}
\frac{\sigma (|D_\E|^{1/2}\cdot (\sqrt{-1})^d)}{|D_\E|^{1/2}\cdot (\sqrt{-1})^d} = \frac{\prod_{v \in S_\infty}\sigma\circ\iota_v(\alpha)}{\prod_{v \in S_\infty}\iota_v(\alpha)} = (-1)^{{}^\sharp\{v \in S_\infty\,\vert\,\sigma \circ \iota_v = \overline{\iota}_{\sigma\circ v}\}}
\end{align*}
for all $\sigma \in {\rm Aut}(\C)$.
This completes the proof.
\end{proof}

\section{Proof of Theorem \ref{T:Deligne}}\label{S:proof}

We prove our main result Theorem \ref{T:Deligne} in this section. First we recall the result due to Liu \cite{Liu2019b} and Pitale--Saha--Schmidt \cite{PSS2020} on the algebraicity of the critical values of the twisted standard $L$-functions for $\GSp_4(\A_\F)$.

\begin{thm}[Liu, Pitale--Saha--Schmidt]\label{T:Liu}
Let $\itPsi$ be a cohomological irreducible cuspidal automorphic representation of $\GSp_4(\A_\F)$ such that $\itPsi_v$ is a holomorphic discrete series representation for each $v \in S_\infty$. 
There exist non-zero complex numbers $p({}^\sigma\!\itPsi)$ defined for each $\sigma \in {\rm Aut}(\C)$ satisfying the following property:
Let $\chi$ be a finite order Hecke character of $\A_\F^\times$ with parallel signature and $m \in \Z_{\geq 1}$ be a critical point of the twisted standard $L$-function $L(s,\itPsi , \,{\rm std}\otimes \chi)$ such that $m \neq 1$ if $\F=\Q$ and $\chi^2=1$. We have
\begin{align*}
L^{(\infty)}(m,\itPsi,\, {\rm std}\otimes\chi)\sim_{\Q(\itPsi)}(2\pi\sqrt{-1})^{3dm}\cdot G(\chi)^3\cdot p(\itPsi),
\end{align*}
and the ratio is equivariant under ${\rm Aut}(\C)$.
\end{thm}

\begin{rmk}
The result of Pitale--Saha--Schmidt was stated for $\F=\Q$ with certain parity conditions on $m$ and $\chi$.
The general case was settled in \cite{Liu2019b} based on explicit computation of archimedean doubling local zeta integrals.
The other results in the literatures on the algebraicity of twisted standard $L$-functions were mainly for scalar-valued Hilbert--Siegel modular forms. For our purpose here, we need to consider vector-valued Hilbert--Siegel modular forms.
\end{rmk}

Now we begin the proof of Theorem \ref{T:Deligne}. 
We may assume $\itPi$ is non-CM.
Let ${\rm Sym}^3\itPi$ be the functorial lift of $\itPi$ with respect to the symmetric cube representation of $\GL_2$. The functoriality was established by Kim and Shahidi \cite{KS2002}.
Note that ${\rm Sym}^3\itPi$ is a cohomological irreducible cuspidal automorphic representation of $\GL_4(\A_\F)$. We have the factorization of twisted exterior square $L$-function of ${\rm Sym}^3\itPi$ by $\omega_\itPi^{-3}$:
\[
L(s,{\rm Sym}^3\itPi,\wedge^2 \otimes \omega_\itPi^{-3}) = L(s,\itPi,{\rm Sym}^4\otimes \omega_\itPi^{-2})\cdot \zeta_\F(s).
\]
Since $L(s,\itPi,{\rm Sym}^4\otimes \omega_\itPi^{-2})$ is holomorphic and non-vanishing at $s=1$ (to be explained below), we see that $L(s,{\rm Sym}^3\itPi,\wedge^2 \otimes \omega_\itPi^{-3})$ has a pole at $s=1$. 
Therefore, by \cite[Theorem 12.1]{GT2011}, there exists an irreducible globally generic cuspidal automorphic representation $\itPsi_{\rm gen} = \bigotimes_v \itPsi_{{\rm gen},v}$ of $\GSp_4(\A_\F)$ such that ${\rm Sym}^3\itPi$ is the strong lift of $\itPsi_{\rm gen}$. 
In particular, $\itPsi_{{\rm gen},v}$ is a generic discrete series representation of $\GSp_4(\F_v)$ for all $v \in S_\infty$.
Let 
\begin{align}\label{E:holomorphic descend}
\itPsi_{\rm hol} = \bigotimes_{v} \itPsi_{{\rm hol},v}
\end{align}
be an irreducible admissible representation of $\GSp_4(\A_\F)$ defined so that $\itPsi_{{\rm hol},v} = \itPsi_{{\rm gen},v}$ if $v$ is finite and $\itPsi_{{\rm hol},v}$ is the holomorphic discrete series representation in the $L$-packet of $\GSp_4(\F_v)$ containing $\itPsi_{{\rm gen},v}$ if $v \in S_\infty$. Since ${\rm Sym}^3\itPi$ is cuspidal, if follows from the global multiplicity formula of Gee and Ta\"{i}bi \cite[Theorem 7.4.1]{GT2019} that $\itPsi_{\rm hol}$ appears in the automorphic discrete spectrum of $\GSp_4(\A_\F)$. It then follows from the temperedness of $\itPsi_{{\rm hol},v}$ for all $v \in S_\infty$ and the result of Wallach \cite[Theorem 4.3]{Wallach1984} that $\itPsi_{\rm hol}$ is cuspidal. Moreover, $\itPsi_{\rm hol}$ is cohomological (cf.\,\cite[Proposition 5.3]{RS2018}).
Since ${\rm Sym}^3\itPi$ is the strong lift of $\itPsi_{\rm hol}$, it is easy to verify that
\[
L(s, \itPsi_{\rm hol}, \,{\rm std}\otimes\chi) = L(s,\itPi,{\rm Sym}^4 \otimes\omega_\itPi^{-2}\chi)
\]
for any Hecke character $\chi$ of $\A_\F^\times$.
Also note that $L(s,\itPi,{\rm Sym}^4\otimes\chi)$ is holomorphic and non-vanishing for ${\rm Re}(s) \geq 1-2{\sf w}$ whenever $\chi$ is unitary.
Indeed, since $\itPi$ is non-CM and cohomological, the functorial lift ${\rm Sym}^4\itPi$ of $\itPi$ with respect to the symmetric fourth representation of $\GL_2$ exists and is cuspidal by the results of Kim \cite[Theorem B]{Kim2003} and Kim--Shahidi \cite[Theorem 3.3.7]{KS2002b}.
The non-vanishing assertion thus follows from the results of Jacquet and Shalika \cite{JS1976} and \cite[Theorem 5.3]{JS1981}.
In particular, we see that $L(m,\itPi,{\rm Sym}^4\otimes\chi) \neq 0$ for all finite order Hecke characters $\chi$ of $\A_\F^\times$ with parallel signature and all $m \in {\rm Crit}^+({\rm Sym}^4\itPi;\,{\rm sgn}(\chi))$.

\subsection{Case $\F \neq \Q$}\label{SS:3.1}
First we consider the case when $\F \neq \Q$.
Since all the critical values are non-zero, by Theorem \ref{T:Liu}, we see that Conjecture \ref{C:Deligne}-(1) holds if and only if it holds for some $\chi$ and some $m$.
Now we verify the conjecture for 
\[
\chi = (|\mbox{ }|_{\A_\F}^{-{\sf w}}\omega_\itPi)^{-2}\chi',\quad m=1-2{\sf w}
\]
for some quadratic Hecke character $\chi'$ with parallel signature $-1$.
Note that this case is excluded in Theorem \ref{T:Liu} when $\F=\Q$ and this is the only place where the assumption $\F\neq\Q$ is used.
Let $\E$ be a totally imaginary quadratic extension of $\F$ such that the base change lifts ${\rm BC}_{\E}(\itPi)$ and ${\rm BC}_{\E}({\rm Sym}^2\itPi)$ of $\itPi$ and ${\rm Sym}^2\itPi$ to $\GL_2(\A_\E)$ and $\GL_3(\A_\E)$, respectively, are cuspidal.
Then ${\rm BC}_{\E}(\itPi)$ and ${\rm BC}_{\E}({\rm Sym}^2\itPi)$ are cohomological with
\begin{align*}
{\rm BC}_{\E}(\itPi)_v &= {\rm Ind}_{B_2(\C)}^{\GL_2(\C)}(\chi_{\kappa_v-1}\boxtimes\chi_{1-\kappa_v})\otimes|\mbox{ }|_\C^{{\sf w}/2},\\
{\rm BC}_{\E}({\rm Sym}^2\itPi)_v &= {\rm Ind}_{B_3(\C)}^{\GL_3(\C)}(\chi_{2\kappa_v-2}\boxtimes{\bf 1}\boxtimes\chi_{2-2\kappa_v})\otimes|\mbox{ }|_\C^{\sf w}
\end{align*}
for each $v \in S_\infty$. Moreover, ${\rm BC}_{\E}(\itPi)$ and ${\rm BC}_{\E}({\rm Sym}^2\itPi)$ are essentially conjugate self-dual with
\[
{\rm BC}_{\E}(\itPi)^\vee = {\rm BC}_{\E}(\itPi)^c \otimes \omega_\itPi^{-1}\circ{\rm N}_{\E/\F},\quad {\rm BC}_{\E}({\rm Sym}^2\itPi)^\vee = {\rm BC}_{\E}({\rm Sym}^2\itPi)^c \otimes \omega_\itPi^{-2}\circ{\rm N}_{\E/\F}.
\]
We have the following factorizations of Rankin--Selberg $L$-function and twisted Asai $L$-functions:
\begin{align*}
L(s,{\rm BC}_{\E}({\rm Sym}^2\itPi)\times{\rm BC}_{\E}(\itPi)) &= L(s,{\rm Sym}^2\itPi \times \itPi)\cdot L(s,{\rm Sym}^2\itPi \times \itPi\otimes\omega_{\E/\F}),\\
L(s,{\rm BC}_{\E}(\itPi),{\rm As}^+\otimes \omega_\itPi^{-1}) &= L(s,\itPi,{\rm Sym}^2\otimes \omega_\itPi^{-1})\cdot L(s,\omega_{\E/\F}),\\
L(s,{\rm BC}_{\E}({\rm Sym}^2\itPi),{\rm As}^-\otimes \omega_\itPi^{-2}) &= L(s,\itPi,{\rm Sym}^4\otimes \omega_\itPi^{-2}\omega_{\E/\F})\cdot L(s,\itPi,{\rm Sym}^2\otimes \omega_\itPi^{-1})\cdot L(s,\omega_{\E/\F}).
\end{align*}
Also it is clear that
$
2\kappa_v-2 > \kappa_v-1 > 0 > 1-\kappa_v > 2-2\kappa_v
$
for all $v \in S_\infty$. 
Therefore, the assumption in Theorem \ref{T:RS 0} is satisfied and we conclude from Theorems \ref{T:RS 0} and \ref{T:Asai 0} that
\begin{align}\label{E:Main proof 1}
\begin{split}
&L^{(\infty)}(m+\tfrac{1}{2},{\rm Sym}^2\itPi \times \itPi)\cdot L^{(\infty)}(m+\tfrac{1}{2},{\rm Sym}^2\itPi \times \itPi\otimes\omega_{\E/\F})\\
&\sim_{\Q(\itPi)}\\
&|D_\F|^{1/2}\cdot (2\pi\sqrt{-1})^{3d(2m-2+3{\sf w})}\cdot G(\omega_\itPi)^9\cdot L^{(\infty)}(1,\itPi,{\rm Sym}^2\otimes \omega_\itPi^{-1})^2\cdot L^{(\infty)}(1,\omega_{\E/\F})^2\\
&\times L^{(\infty)}(1,\itPi,{\rm Sym}^4\otimes \omega_\itPi^{-2}\omega_{\E/\F})
\end{split}
\end{align}
for all critical points $m+\tfrac{1}{2}$ of $L(s,{\rm Sym}^2\itPi \times \itPi)$, and the ratios are equivariant under ${\rm Aut}(\C)$.
Here we have used the fact that $\omega_{{\rm BC}_{\E}(\itPi)} = \omega_\itPi\circ{\rm N}_{\E/\F}$ and 
\[
\sigma \left( \frac{G(\omega_\itPi\circ {\rm N}_{\E/\F})}{G(\omega_\itPi)^2}\right) = \frac{G({}^\sigma\!\omega_\itPi\circ {\rm N}_{\E/\F})}{G({}^\sigma\!\omega_\itPi)^2}
\]
for all $\sigma \in {\rm Aut}(\C)$, which is a direct consequence of (\ref{E:Galois Gauss sum}).
In the following theorem, we recall the known results in the literatures on the algebraicity of the critical $L$-values appearing in the ratio (\ref{E:Main proof 1}) except for $L(1,\itPi,{\rm Sym}^4\otimes \omega_\itPi^{-2}\omega_{\E/\F})$.

\begin{thm}\label{T:algebraicity results}
\noindent
\begin{itemize}
\item[(1)]{\rm (Garrett--Harris}\,\cite[Theorem 4.6]{GH1993}{\rm,\,C.-}\,\cite[Theorem 1.2]{Chen2021}{\rm)}
Let $\chi$ be a finite order Hecke character of $\A_\F^\times$ and $m+\tfrac{1}{2} \in \Z+\tfrac{1}{2}$ be a critical point of the triple product $L$-function $L(s,\itPi \times \itPi \times \itPi\otimes\chi)$. 
We have
\begin{align*}
&L^{(\infty)}(m+\tfrac{1}{2},\itPi \times \itPi \times \itPi\otimes\chi)\\
&\sim_{\Q(\itPi)\Q(\chi)}\\
&(2\pi\sqrt{-1})^{4dm+\sum_{v\in S_\infty}(3\kappa_v+6{\sf w}+2)}\cdot(\sqrt{-1})^{d{\sf w}}\cdot G(\omega_\itPi)^6 \cdot G(\chi)^4\cdot \Vert f_\itPi \Vert^3,
\end{align*}
and the ratio is equivariant under ${\rm Aut}(\C)$.
\item[(2)]{\rm (Shimura}\,\cite[Theorem 4.3]{Shimura1978}{\rm)} There exists a sequence of non-zero complex numbers $(p({}^\sigma\!\itPi,\underline{\varepsilon}))_{\sigma \in {\rm Aut}(\C)}$ defined for each $\underline{\varepsilon} \in \{\pm1\}^{S_\infty}$ satisfying the following properties:
\begin{itemize}
\item[(i)] Let $\chi$ be a finite order Hecke character of $\A_\F^\times$ and $m+\tfrac{1}{2} \in \Z+\tfrac{1}{2}$ be a critical point of the tiwsted standard $L$-function $L(s,\itPi \otimes \chi)$. We have
\begin{align*}
L^{(\infty)}(m+\tfrac{1}{2},\itPi \times \chi)&\sim_{\Q(\itPi)\Q(\chi)}\\
&(2\pi\sqrt{-1})^{dm+\sum_{v \in S_\infty}(\kappa_v+{\sf w})/2}\cdot G(\chi) \cdot p(\itPi,(-1)^m\cdot{\rm sgn}(\chi)),
\end{align*}
and the ratio is equivariant under ${\rm Aut}(\C)$.
\item[(ii)] Let $\underline{\varepsilon}\in \{\pm1\}^{S_\infty}$.
We have
\begin{align*}
p(\itPi,\underline{\varepsilon})\cdot p(\itPi,-\underline{\varepsilon})\sim_{\Q(\itPi)}(2\pi\sqrt{-1})^d\cdot(\sqrt{-1})^{d{\sf w}}\cdot G(\omega_\itPi)\cdot\Vert f_\itPi \Vert,
\end{align*}
and the ratio is equivariant under ${\rm Aut}(\C)$.
\end{itemize}
\item[(3)]{\rm (Sturm}\,\cite{Sturm1989}{\rm,\,Im}\,\cite{Im1991}{\rm)} Let $m \in \Z_{\geq 1}$ be a critical point of the twisted symmetric square $L$-function $L(s,\itPi,{\rm Sym}^2\otimes\omega_\itPi^{-1})$. We have
\begin{align*}
L^{(\infty)}(m,\itPi,{\rm Sym}^2 \otimes \omega_\itPi^{-1})\sim_{\Q(\itPi)}(2\pi\sqrt{-1})^{2dm+\sum_{v \in S_\infty}\kappa_v}\cdot(\sqrt{-1})^{d{\sf w}}\cdot \Vert f_\itPi \Vert,
\end{align*}
and the ratio is equivariant under ${\rm Aut}(\C)$.
\item[(4)]{\rm (Shimura}\,\cite[Proposition 3.1]{Shimura1978}{\rm)} Let $\chi$ be a finite order Hecke character of $\A_\F^\times$ and $m \in \Z_{\geq 1}$ be a critical point of $L(s,\chi)$. We have
\begin{align*}
L^{(\infty)}(m,\chi)\sim_{\Q(\chi)}|D_\F|^{1/2}\cdot (2\pi\sqrt{-1})^{dm}\cdot G(\chi),
\end{align*}
and the ratio is equivariant under ${\rm Aut}(\C)$.
\end{itemize} 
\end{thm}

Let $m+\tfrac{1}{2}$ be a critical point of $L(s,{\rm Sym}^2\itPi \times \itPi )$ with $m+\tfrac{1}{2} \neq -\tfrac{3}{2}{\sf w}+\tfrac{1}{2}$, which is the possible central critical point, then 
\[
L(m+\tfrac{1}{2},{\rm Sym}^2\itPi \times \itPi\otimes\chi) \neq 0,\quad L(m+\tfrac{1}{2},\itPi\otimes\omega_\itPi\cdot\chi)\neq0
\]
for any unitary Hecke character $\chi$ of $\A_\F^\times$ by \cite[Theorem 5.3]{JS1981} and \cite[Theorem 5.2]{Shahidi1981}.
The existence of such critical point is guaranteed by our assumption that $\kappa_v \geq 3$ for all $v \in S_\infty$.
Note that we have the factorization
\[
L(s,\itPi \times \itPi \times \itPi \otimes\chi) = L(s,{\rm Sym}^2\itPi\times\itPi\otimes\chi)\cdot L(s,\itPi \otimes\omega_\itPi\cdot\chi).
\]
It is clear that $m+\tfrac{1}{2}$ is also a critical point of $L(s,\itPi\otimes\omega_\itPi\cdot\chi)$.
We thus deduce from Theorem \ref{T:algebraicity results}-(1) and (2) that
\begin{align*}
&L^{(\infty)}(m+\tfrac{1}{2},{\rm Sym}^2\itPi \times \itPi\otimes\chi)\\
&\sim_{\Q(\itPi)\Q(\chi)}\\
&(2\pi\sqrt{-1})^{3dm+\sum_{v \in S_\infty}(5\kappa_v+9{\sf w}+2)/2}\cdot G(\omega_\itPi)^4\cdot G(\chi)^3\cdot \Vert f_\itPi \Vert^2 \cdot p(\itPi,(-1)^{m+1}\cdot {\rm sgn}(\chi)),
\end{align*}
for all finite order Hecke characters $\chi$ of $\A_\F^\times$, and the ratios are equivariant under ${\rm Aut}(\C)$.
In particular, combining with Theorem \ref{T:algebraicity results}-(2)-(ii), we have
\begin{align}\label{E:Main proof 2}
\begin{split}
&\sigma \left(\frac{L^{(\infty)}(m+\tfrac{1}{2},{\rm Sym}^2\itPi \times \itPi)\cdot L^{(\infty)}(m+\tfrac{1}{2},{\rm Sym}^2\itPi \times \itPi\otimes\omega_{\E/\F})}{(2\pi\sqrt{-1})^{6dm+\sum_{v \in S_\infty}(5\kappa_v+9{\sf w}+3)}\cdot (\sqrt{-1})^{d{\sf w}}\cdot G(\omega_\itPi)^9\cdot G(\omega_{\E/\F})\cdot \Vert f_\itPi \Vert^5} \right)\\
& = \frac{L^{(\infty)}(m+\tfrac{1}{2},{\rm Sym}^2{}^\sigma\!\itPi \times {}^\sigma\!\itPi)\cdot L^{(\infty)}(m+\tfrac{1}{2},{\rm Sym}^2{}^\sigma\!\itPi \times {}^\sigma\!\itPi\otimes\omega_{\E/\F})}{(2\pi\sqrt{-1})^{6dm+\sum_{v \in S_\infty}(5\kappa_v+9{\sf w}+3)}\cdot (\sqrt{-1})^{d{\sf w}}\cdot G({}^\sigma\!\omega_\itPi)^9\cdot G(\omega_{\E/\F})\cdot \Vert f_{{}^\sigma\!\itPi} \Vert^5}
\end{split}
\end{align}
for all $\sigma \in {\rm Aut}(\C)$.
By Theorem \ref{T:algebraicity results}-(3) and (4), we have
\begin{align}\label{E:Main proof 3}
\begin{split}
&\sigma \left( \frac{L^{(\infty)}(1,\itPi,{\rm Sym}^2\otimes \omega_\itPi^{-1})^2\cdot L^{(\infty)}(1,\omega_{\E/\F})^2}{(2\pi\sqrt{-1})^{6d+2\sum_{v \in S_\infty}\kappa_v}\cdot\Vert f_\itPi \Vert^2}\right)\\
& = \frac{L^{(\infty)}(1,{}^\sigma\!\itPi,{\rm Sym}^2\otimes {}^\sigma\!\omega_\itPi^{-1})^2\cdot L^{(\infty)}(1,\omega_{\E/\F})^2}{(2\pi\sqrt{-1})^{6d+2\sum_{v \in S_\infty}\kappa_v}\cdot\Vert f_{{}^\sigma\!\itPi} \Vert^2}
\end{split}
\end{align}
for all $\sigma \in {\rm Aut}(\C)$.
It then follows from (\ref{E:Main proof 1})-(\ref{E:Main proof 3}) that
\begin{align*}
&\sigma \left(\frac{L^{(\infty)}(1,\itPi,{\rm Sym}^4 \otimes \omega_\itPi^{-2}\omega_{\E/\F})}{|D_\F|^{1/2}\cdot(2\pi\sqrt{-1})^{3d+3\sum_{v \in S_\infty}\kappa_v}\cdot (\sqrt{-1})^{d{\sf w}}\cdot G(\omega_{\E/\F})\cdot \Vert f_\itPi \Vert^3} \right)\\
& = \frac{L^{(\infty)}(1,{}^\sigma\!\itPi,{\rm Sym}^4 \otimes\omega_\itPi^{-2}\omega_{\E/\F})}{|D_\F|^{1/2}\cdot(2\pi\sqrt{-1})^{3d+3\sum_{v \in S_\infty}\kappa_v}\cdot (\sqrt{-1})^{d{\sf w}}\cdot G(\omega_{\E/\F})\cdot \Vert f_{{}^\sigma\!\itPi} \Vert^3}
\end{align*}
for all $\sigma \in {\rm Aut}(\C)$.
This completes the proof of Conjecture \ref{C:Deligne}-(1).
Now we show that Conjecture \ref{C:Deligne}-(2) holds. Let $m-2{\sf w}\in {\rm Crit}^-({\rm Sym}^4\itPi;\,{\rm sgn}(\chi))$. It is easy to verify that (cf.\,Lemma \ref{L:critical gamma factor})
\[
\prod_{v \in S_\infty} \gamma(m-2{\sf w},\itPi_v,{\rm Sym}^4\otimes \chi_v, \psi_{\F_v}) \in (2\pi\sqrt{-1})^{5dm-3d}\cdot \Q^\times.
\]
Note that the central character of the cuspidal automorphic representation ${\rm Sym}^4\itPi \otimes \chi$ of $\GL_5(\A_\F)$ is equal to $\chi^5\cdot\omega_\itPi^{10}$. 
By Lemma \ref{L:Galois-equiv. root number} below, we have
\begin{align*}
\sigma \left( \frac{\prod_{v \nmid \infty}\varepsilon(m-2{\sf w},\itPi_v,{\rm Sym}^4\otimes \chi_v, \psi_{\F_v})}{|D_\F|^{1/2}\cdot G(\chi\cdot\omega_\itPi^2)^5}\right) = \frac{\prod_{v \nmid \infty}\varepsilon(m-2{\sf w},{}^\sigma\!\itPi_v,{\rm Sym}^4\otimes {}^\sigma\!\chi_v, \psi_{\F_v})}{|D_\F|^{1/2}\cdot G({}^\sigma\!\chi\cdot{}^\sigma\!\omega_\itPi^2)^5}
\end{align*}
for all $\sigma \in {\rm Aut}(\C)$.
Therefore, Conjecture \ref{C:Deligne}-(2) follows from Conjecture \ref{C:Deligne}-(1) for the critical value $L^{(\infty)}(1-m+2{\sf w},\itPi^\vee,{\rm Sym}^4\otimes\chi^{-1})$ and the global functional equation for $\GL_5(\A_\F)$ that
\begin{align*}
L^{(\infty)}(m-2{\sf w}, \itPi,{\rm Sym}^4\otimes \chi)& = \prod_{v \nmid \infty}\varepsilon(m-2{\sf w},\itPi_v,{\rm Sym}^4\otimes \chi_v, \psi_{\F_v})\cdot \prod_{v \in S_\infty} \gamma(m-2{\sf w},\itPi_v,{\rm Sym}^4\otimes \chi_v, \psi_{\F_v})\\
&\times L^{(\infty)}(1-m+2{\sf w}, \itPi^\vee,{\rm Sym}^4\otimes \chi^{-1}).
\end{align*}
This completes the proof when $\F \neq \Q$.

\subsection{Case $\F=\Q$}
Now we prove the case when $\F=\Q$ by base change trick. Let $\kappa\in\Z_{\geq 3}$ be the weight of $\itPi_\infty$. Let $\L$ be a totally real cyclic extension over $\Q$ such that $r = [\L:\Q] \geq 2$ is a prime number.
Let $\omega_{\mathbb{L}/\Q}$ be a non-trivial Hecke character of $\A_\Q^\times$ that vanishing on ${\rm N}_{\mathbb{L}/\Q}(\A_{\mathbb{L}}^\times)$. Note that by class field theory we have $\omega_{\mathbb{L}/\Q}^r=1$ and
\[
\zeta_{\mathbb{L}}(s) = \prod_{i=1}^r L(s,\omega_{{\mathbb{L}}/\Q}^i).
\]
In particular, it follows from Theorem \ref{T:algebraicity results}-(4) that 
\begin{align}\label{E:Main proof 4}
\sigma\left(\frac{G(\omega_{{\mathbb{L}}/\Q})^{r(r-1)/2}}{|D_{\mathbb{L}}|^{1/2}}\right) = \frac{G({}^\sigma\!\omega_{{\mathbb{L}}/\Q})^{r(r-1)/2}}{|D_{\mathbb{L}}|^{1/2}}
\end{align}
for all $\sigma \in {\rm Aut}(\C)$.
Let ${\rm BC}_{\mathbb{L}}(\itPi)$ be the base change lift of $\itPi$ to $\GL_2(\A_{\mathbb{L}})$. Then ${\rm BC}_{\mathbb{L}}(\itPi)$ is cohomological and cuspidal with 
\[
{\rm BC}_{\mathbb{L}}(\itPi)_v = \itPi_\infty
\]
for all archimedean places $v$ of $\mathbb{L}$. Note that $\omega_{{\rm BC}_{\mathbb{L}}(\itPi)} = \omega_\itPi\circ{\rm N}_{\mathbb{L}/\Q}$.
Let $f_{{\rm BC}_{\mathbb{L}}(\itPi)}$ be the normalized newform of ${{\rm BC}_{\mathbb{L}}(\itPi)}$.
By Theorem \ref{T:algebraicity results}-(2), one can deduce that
\begin{align}\label{E:Main proof 5}
\sigma \left( \frac{\Vert f_{{\rm BC}_{\mathbb{L}}(\itPi)}\Vert}{\Vert f_\itPi\Vert^r}\right) = \frac{\Vert f_{{\rm BC}_{\mathbb{L}}({}^\sigma\!\itPi)}\Vert}{\Vert f_{{}^\sigma\!\itPi}\Vert^r}
\end{align}
for all $\sigma \in {\rm Aut}(\C)$.
Let $\chi$ be a finite order Hecke character of $\A_\Q^\times$. It follows easily from the definition of Gauss sum that
\begin{align}\label{E:Main proof 6}
\sigma \left( \frac{G(\chi\circ {\rm N}_{\mathbb{L}/\Q})}{G(\chi)^r} \right) =  \frac{G({}^\sigma\!\chi\circ {\rm N}_{\mathbb{L}/\Q})}{G({}^\sigma\!\chi)^r}
\end{align}
for all $\sigma \in {\rm Aut}(\C)$. We have the following factorization of twisted symmetric fourth $L$-function:
\begin{align}\label{E:Main proof 7}
\begin{split}
L(s,{\rm BC}_{\mathbb{L}}(\itPi),{\rm Sym}^4\otimes\chi\circ{\rm N}_{\mathbb{L}/\Q})  = \prod_{i=1}^r L(s,\itPi,{\rm Sym}^4\otimes\chi\omega_{\mathbb{L}/\Q}^i).
\end{split}
\end{align}
Consider the critical value 
\[
L^{(\infty)}(2-2{\sf w},{\rm BC}_{\mathbb{L}}(\itPi),{\rm Sym}^4\otimes {\bf 1}).
\]
As we have proved in \S\,\ref{SS:3.1}, Conjecture \ref{C:Deligne} holds for this critical value. 
Therefore, by Theorem \ref{T:Liu} and (\ref{E:Main proof 4})-(\ref{E:Main proof 7}), we obtain the period relation
\[
\sigma\left( \frac{p(\itPsi_{\rm hol})^r}{(2\pi\sqrt{-1})^{3r\kappa}\cdot (\sqrt{-1})^{r{\sf w}}\cdot\Vert f_\itPi\Vert^{3r}}\right) = \frac{p({}^\sigma\!\itPsi_{\rm hol})^r}{(2\pi\sqrt{-1})^{3r\kappa}\cdot  (\sqrt{-1})^{r{\sf w}}\cdot\Vert f_{{}^\sigma\!\itPi}\Vert^{3r}}
\]
for all $\sigma \in {\rm Aut}(\C)$. Here $\itPsi_{\rm hol}$ is the irreducible cuspidal automorphic representation of $\GSp_4(\A_\Q)$ in (\ref{E:holomorphic descend}).
Taking $r=2,3$, we then deduce the period relation
\[
\sigma\left( \frac{p(\itPsi_{\rm hol})}{(2\pi\sqrt{-1})^{3\kappa}\cdot (\sqrt{-1})^{{\sf w}}\cdot\Vert f_\itPi\Vert^{3}}\right) = \frac{p({}^\sigma\!\itPsi_{\rm hol})}{(2\pi\sqrt{-1})^{3\kappa}\cdot (\sqrt{-1})^{{\sf w}}\cdot\Vert f_{{}^\sigma\!\itPi}\Vert^{3}}
\]
for all $\sigma \in {\rm Aut}(\C)$.
By this period relation and Theorem \ref{T:Liu} again, we conclude that Conjecture \ref{C:Deligne} holds for all $m-2{\sf w} \in {\rm Crit}^+({\rm Sym}^4\itPi;\,{\rm sgn}(\chi))$ with $m \neq 1$ if $\chi^2=1$. Now we consider the remaining cases in the right-half critical region. Assume $\chi^2=1$ and ${\rm sgn}(\chi)=-1$. Take $r=3$. Since Conjecture \ref{C:Deligne} holds for the critical values
\[
L^{(\infty)}(1-2{\sf w},{\rm BC}_{\L}(\itPi),{\rm Sym}^4\otimes\chi\circ{\rm N}_{\L/\Q}),\quad L^{(\infty)}(1-2{\sf w},\itPi,{\rm Sym}^4\otimes\chi\omega_{\mathbb{L}/\Q}),\quad L^{(\infty)}(1-2{\sf w},\itPi,{\rm Sym}^4\otimes\chi\omega_{\mathbb{L}/\Q}^2),
\]
we deduce from (\ref{E:Main proof 4})-(\ref{E:Main proof 7}) again that Conjecture \ref{C:Deligne} also holds for $L^{(\infty)}(1-2{\sf w},\itPi,{\rm Sym}^4\otimes\chi)$. This completes the proof of Conjecture \ref{C:Deligne}-(1). Similarly as in \S\,\ref{SS:3.1}, we then conclude Conjecture \ref{C:Deligne}-(2) from the global functional equation for $\GL_5(\A_\Q)$. This completes the proof.
\section{Whittaker periods for $\GL_2$ and $\GL_3$}\label{S:Whittaker periods}

Let $\E$ be a totally imaginary quadratic extension of $\F$.
We follow the notation in \S\,\ref{SS:CM-fields}.

\subsection{Irreducible algebraic representations of $\GL_n$}\label{SS:rep GLn}

For $\mu\in X^+(T_n)$, let $(\rho_\mu,M_\mu)$ be the irreducible algebraic representation of $\GL_n$ with highest weight $\mu$. Let $M_{\mu,\C} = M_\mu \otimes_\Q \C$ be the base change to a representation of $\GL_n(\C)$. We denote by $V^{(n)}_\mu = M_{\mu,\C}\vert_{{\rm U}(n)}$ the irreducible representation of ${\rm U}(n)$ obtained by restriction of $M_{\mu,\C}$. Note that $M_\mu^\vee = M_{\mu^\vee}$.

Let 
\[
\mu = \prod_{v \in S_\infty} \mu_v \in \prod_{v \in S_\infty}(X^+(T_n) \times X^+(T_n))
\]
be a sequence of pairs of dominant integral weights of $T_n$ indexed by $S_\infty$. Write $\mu_v= (\mu_{\iota_v},\mu_{\overline{\iota}_v}) \in X^+(T_n) \times X^+(T_n)$ for $v \in S_\infty$.
We say $\mu$ is pure if there exists an integer ${\sf w}$ such that $\mu_{\overline{\iota}_v} = \mu_{\iota_v}^\vee + ({\sf w},\cdots,{\sf w})$ for all $v \in S_\infty$.
For $v \in S_\infty$, let $(\rho_{\mu_v},M_{\mu_v,\C})$ be the complex representation of $\GL_n(\C)$ defined by $M_{\mu_v,\C} = M_{\mu_{\iota_v},\C}\otimes M_{\mu_{\overline{\iota}_v},\C}$ and 
\[
\rho_{\mu_v}(g)\cdot ({\bf v}_{\iota_v} \otimes {\bf v}_{\overline{\iota}_v}) = \rho_{\mu_{\iota_v}}(g){\bf v}_{\iota _v}\otimes \rho_{\mu_{\overline{\iota}_v}}(\overline{g}){\bf v}_{\overline{\iota}_v}.
\]
Let $\mu_v^c = (\mu_{\overline{\iota}_v},\mu_{\iota_v})$ and
\[
{\bf c}^{(n)}_{\mu_v} : M_{\mu_v,\C} \longrightarrow M_{\mu_v^c,\C}
\]
be the isomorphism defined by
\[
{\bf c}^{(n)}_{\mu_v} ({\bf v}_{\iota_v}\otimes {\bf v}_{\overline{\iota}_v}) = {\bf v}_{\overline{\iota}_v} \otimes {\bf v}_{\iota_v}.
\]
Then it is clear that ${\bf c}_{\mu_v}^{(n)}$ is $\GL_n(\C)$-conjugate-equivariant, that is, 
\[
{\bf c}^{(n)}_{\mu_v}\circ\rho_{\mu_v}(g) = \rho_{\mu_v^c}(\overline{g})\circ{\bf c}^{(n)}_{\mu_v}
\]
for all $g \in \GL_n(\C)$.
Define the representation $(\rho_\mu,M_{\mu,\C})$ of $\GL_n(\E_\infty)$ by
\[
(\rho_\mu,M_{\mu,\C}) = (\otimes_{v \in S_\infty}\rho_{\mu_v},\otimes_{v \in S_\infty}M_{\mu_v,\C}).
\]  
For $\sigma \in {\rm Aut}(\C)$, let ${}^\sigma\!\mu = \prod_{v \in S_\infty}{}^\sigma\!\mu_v\in \prod_{v \in S_\infty}(X^+(T_n) \times X^+(T_n))$ defined by ${}^\sigma\!\mu_v = (\mu_{\sigma^{-1}\circ \iota_v},\mu_{\sigma^{-1}\circ \overline{\iota}_v})$ for $v \in S_\infty$. We have the $\sigma$-linear isomorphism $M_{\mu,\C} \rightarrow M_{{}^\sigma\!\mu,\C}$ defined by
\begin{align}\label{E:sigma linear}
\bigotimes_{v \in S_\infty}z_v\cdot({\bf v}_{\iota_v} \otimes {\bf v}_{\overline{\iota}_v}) = \bigotimes_{v \in S_\infty} \sigma(z_v) \cdot ({\bf v}_{\sigma^{-1}\circ\iota_v} \otimes {\bf v}_{\sigma^{-1}\circ\overline{\iota}_v})
\end{align}
for $z_v \in \C$, ${\bf v}_{\iota_v} \in M_{\mu_{\iota_v}}$ and ${\bf v}_{\overline{\iota}_v} \in M_{\mu_{\overline{\iota}_v}}$.

\subsubsection{Irreducible algebraic representations of $\GL_2$}\label{SS:rep GL2}
Let $n=2$ and $\mu = (\mu_1,\mu_2) \in X^+(T_2)$. We fix a model of $(\rho_\mu,M_\mu)$ as follows: 
Let $M_\mu$ be the $\Q$-vector space consisting of homogeneous polynomials over $\Q$ of degree $\mu_1-\mu_2$ in variables $x$ and $y$. For $g \in \GL_2(\Q)$, we define $\rho_\mu(g) \in {\rm End}_\Q(M_\mu)$ by
\[
\rho_\mu(g)\cdot P(x,y) = \det(g)^{\mu_2}\cdot P((x,y)g)
\]
for $P \in M_\mu$. Note that a highest weight vector of $M_\mu$ is given by $x^{\mu_1-\mu_2}$. We have the following relations in $V_\mu^{(2)}$:
\begin{align}\label{E:Lie algebra action 1}
\begin{split}
E_+\cdot x^{\mu_1-\mu_2-i}y^i &= i\cdot x^{\mu_1-\mu_2-i+1}y^{i-1},\\
E_-\cdot x^{\mu_1-\mu_2-i}y^i &= (\mu_1-\mu_2-i) \cdot x^{\mu_1-\mu_2-i-1}y^{i+1}
\end{split}
\end{align}
for $0 \leq i \leq \mu_1-\mu_2$.
Using the above relation, it is easy to see that a generator of $(M_\mu \otimes M_{\mu^\vee})^{\GL_2}$ is given by
\begin{align}\label{E:generator 1}
\sum_{i=0}^{\mu_1-\mu_2}(-1)^i{\mu_1-\mu_2\choose i}\cdot x^{\mu_1-\mu_2-i}y^i \otimes x^iy^{\mu_1-\mu_2-i}.
\end{align}
Let 
\[
c^{(2)}_\mu : V_\mu^{(2)} \longrightarrow V_{\mu^\vee}^{(2)}
\]
be the ${\rm SU}(2)$-conjugate-equivariant isomorphism normalized so that
\begin{align*}
c_\mu^{(2)}(x^{\mu_1-\mu_2}) = y^{\mu_1-\mu_2}.
\end{align*}
Since $c_\mu^{(2)}\circ E_\pm = -E_\mp \circ c_\mu^{(2)}$, it follows from (\ref{E:Lie algebra action 1}) that 
\begin{align}\label{E:conjugate equiv. 1}
c_\mu^{(2)}(x^{\mu_1-\mu_2-i}y^i) =(-1)^i\cdot x^iy^{\mu_1-\mu_2-i}
\end{align}
for $0 \leq i \leq \mu_1-\mu_2$.
For example, the adjoint representation of ${\rm U}(2)$ on $\frak{g}_{2,\C}/\frak{k}_{2,\C}$ is isomorphic to $V_{(1,-1)}^{(2)}$. More precisely, we have the ${\rm U}(2)$-equivariant isomorphism $\frak{g}_{2,\C}/\frak{k}_{2,\C} \rightarrow V_{(1,-1)}^{(2)}$ given by
\begin{align}\label{E:conjugate equiv. 2}
 Y_{(1,-1)}\longmapsto x^2,\quad  Y_{(0,0)} \longmapsto 2xy ,\quad Y_{(-1,1)} \longmapsto y^2.
\end{align}
Also $c_{(1,-1)}^{(2)} : \frak{g}_{2,\C}/\frak{k}_{2,\C} \rightarrow \frak{g}_{2,\C}/\frak{k}_{2,\C}$ is given by
\begin{align}\label{E:conjugate equiv. Lie 1}
Y_{(0,0)}\longmapsto -Y_{(0,0)},\quad Y_{\pm(1,-1)} \longmapsto Y_{\mp(1,-1)}.
\end{align}

\subsubsection{Irreducible algebraic representations of $\GL_3$}\label{SS:rep GL3}
Let $n=3$ and $\mu = (\mu_1,\mu_2,\mu_3) \in X^+(T_3)$. We fix a model of $(\rho_\mu,M_\mu)$ as follows:
Let $M_\mu$ be the $\Q$-vector space generated by homogeneous polynomials of degree $\mu_1+\mu_2-2\mu_3$ in variable $x = (x_{ij})_{1 \leq i \leq 2,\, 1 \leq j \leq 3}$ of the form
\[
{\bf v}_{(\mu;\, n_1,n_2,n_3,n_{23},n_{13},n_{12})} = \prod_{i=1}^3 x_{2i}^{n_i} \prod_{1 \leq j < j' \leq 3}\det \bp x_{1j} & x_{1j'} \\ x_{2j} & x_{2j'}\ep^{n_{j,j'}}
\]
with $n_1+n_2+n_3 = \mu_1-\mu_2$ and $n_{12}+n_{23}+n_{13} = \mu_2-\mu_3$.
For $g \in \GL_3(\Q)$, we define $\rho_\mu(g) \in {\rm End}_\Q(M_\mu)$ by
\[
\rho_\mu(g)\cdot P(x) = \det(g)^{\mu_3}\cdot P(xg)
\]
for $P \in M_\mu$.
Note that a highest weight vector of $M_\mu$ is given by ${\bf v}_{(\mu;\,\mu_1-\mu_2,0,0,0,0,\mu_2-\mu_3)}$.
We have the following relation in $V_\mu^{(3)}$:
For $1 \leq i \neq j \leq 3$, we have
\begin{align}\label{E:Lie algebra action 2}
E_{ij}\cdot {\bf v}_{(\mu;\,\underline{n})} = n_j\cdot {\bf v}_{(\mu;\,\underline{n}+f_i-f_j)} - (-1)^k n_{kj}\cdot {\bf v}_{(\mu;\,\underline{n}+f_{ki}-f_{kj})}
\end{align}
for $\underline{n} = (n_1,n_2,n_3,n_{23},n_{13},n_{12}) \in \Z_{\geq 0}^6$ with $n_1+n_2+n_3 = \mu_1-\mu_2$ and $n_{12}+n_{23}+n_{13} = \mu_2-\mu_3$. 
Here ${\bf v}_{(\mu;\,\underline{n}')}=0$ if $\underline{n}' \notin \Z_{\geq 0}^6$, $\{i,j,k\} = \{1,2,3\}$, and 
\begin{gather*}
f_1 = (1,0,0,0,0,0),\quad f_2 = (0,1,0,0,0,0),\quad f_3 = (0,0,1,0,0,0),\\ 
f_{23}=f_{32} = (0,0,0,1,0,0),\quad f_{13}=f_{31} = (0,0,0,0,1,0),\quad f_{12}=f_{21} = (0,0,0,0,0,1).
\end{gather*}
Let
\[
c_\mu^{(3)} : V_\mu^{(3)} \longrightarrow V_{\mu^\vee}^{(3)}
\]
be the ${\rm SU}(3)$-conjugate-equivariant isomorphism normalized so that
\begin{align*}
c_\mu^{(3)}({\bf v}_{(\mu;\,0,0,\mu_1-\mu_2,0,0,\mu_2-\mu_3)}) = {\bf v}_{(\mu^\vee;\,0,0,\mu_2-\mu_3,0,0,\mu_1-\mu_2)}.
\end{align*}
For example, the adjoint representation of ${\rm U}(3)$ on $\frak{g}_{3,\C} / \frak{k}_{3,\C}$ is isomorphic to $V_{(1,0,-1)}^{(3)}$.
More precisely, we have the ${\rm U}(3)$-equivariant isomorphism $\frak{g}_{3,\C}/\frak{k}_{3,\C} \rightarrow V_{(1,0,-1)}^{(3)}$ given by
\begin{align*}
\begin{split}
Z_{12}&\longmapsto 2 \cdot{\bf v}_{((1,0,-1);\,1,0,0,1,0,0)} + {\bf v}_{((1,0,-1);\,0,0,1,0,0,1)},\\
Z_{23} &\longmapsto -{\bf v}_{((1,0,-1);\,1,0,0,1,0,0)} -2\cdot {\bf v}_{((1,0,-1);\,0,0,1,0,0,1)},\\
X_{(1,-1,0)} &\longmapsto {\bf v}_{((1,0,-1);\,1,0,0,0,1,0)},\quad X_{(-1,1,0)}\longmapsto {\bf v}_{((1,0,-1);\,0,1,0,1,0,0)},\\
X_{(0,1,-1)} &\longmapsto -{\bf v}_{((1,0,-1);\,0,1,0,0,0,1)},\quad X_{(0,-1,1)}\longmapsto -{\bf v}_{((1,0,-1);\,0,0,1,0,1,0)},\\
X_{(1,0,-1)} &\longmapsto -{\bf v}_{((1,0,-1);\,1,0,0,0,0,1)},\quad X_{(-1,0,1)}\longmapsto {\bf v}_{((1,0,-1);\,0,0,1,1,0,0)}.
\end{split}
\end{align*}
We refer to \cite[Table 1]{HO2009} for the table of the adjoint representation.
Also $c_{(1,0,-1)}^{(3)} : \frak{g}_{3,\C}/\frak{k}_{3,\C} \rightarrow \frak{g}_{3,\C}/\frak{k}_{3,\C}$ is given by 
\begin{gather}
Z_{12} \longmapsto Z_{12},\quad Z_{23} \longmapsto Z_{23},\label{E:conjugate equiv. Lie 2}\\
X_{\pm (1,-1,0)} \longmapsto - X_{\mp (1,-1,0)},\quad X_{\pm (0,1,-1)} \longmapsto - X_{\mp (0,-1,1)},\quad X_{\pm (1,0,-1)} \longmapsto - X_{\mp (-1,0,1)}.\label{E:conjugate equiv. Lie 3}
\end{gather}
We also need to consider the restriction of $M_\mu$ to $\GL_{2}$, where we identify $\GL_2$ as a subgroup of $\GL_3$ via the embedding $g \mapsto \bp g & 0 \\ 0 & 1\ep$.
For $\lambda = (\lambda_1,\lambda_2) \in X^+(T_2)$, we write $\mu \succ \lambda$ if
\[
\mu_1 \geq \lambda_1 \geq \mu_2 \geq \lambda_2 \geq \mu_3.
\]
Then we have (cf.\,\cite[Theorem 8.1.2]{{GW2009}})
\[
M_\mu \vert_{\GL_2} = \bigoplus_{\mu \succ \lambda} M_\lambda.
\]
The subsequent archimedean computations in \S\,\ref{S:RS}--\S\,\ref{S:Asai} involve mainly the following $\lambda$:
\[
\lambda = (\mu_1,\mu_2),\quad \lambda = (\mu_2,\mu_2),\quad \lambda = (\mu_2,\mu_3).
\]
The highest (resp.\,lowest) weight vectors of $M_\mu$ appear in $M_{(\mu_1,\mu_2)}$ (resp.\,$M_{(\mu_2,\mu_3)}$), and $M_{(\mu_2,\mu_2)}$ is the unique $\GL_2$-type of $M_\mu$ with dimension one.
For $\mu \succ \lambda$, we have the $\GL_2$-equivariant homomorphism
\[
M_\lambda \longrightarrow M_\mu, \quad x^{\lambda_1-\lambda_2} \longmapsto {\bf v}_{(\mu;\,\lambda_1-\mu_2,0,\mu_1-\lambda_1,0,\mu_2-\lambda_2,\lambda_2-\mu_3)}.
\] 
For $\mu \succ \lambda$ and $0 \leq i \leq \lambda_1-\lambda_2$, let ${\bf v}_{(\mu;\,\lambda,i)}$ denotes the image of $x^{\lambda_1-\lambda_2-i}y^i$ under the above homomorphism. 
A highest and a lowest weight vectors of $M_\mu$ are given by
\[
{\bf v}_{(\mu;\,(\mu_1,\mu_2),0)},\quad {\bf v}_{(\mu;\,(\mu_2,\mu_3),\mu_2-\mu_3)}
\]
respectively.
A generator of $(M_\mu \otimes M_{\mu^\vee})^{\GL_3}$ is of the form
\begin{align}\label{E:generator 2}
\sum_{\mu \succ \lambda}C_{(\mu;\,\lambda)}\sum_{i=0}^{\lambda_1-\lambda_2}(-1)^i{\lambda_1-\lambda_2 \choose i}\cdot {\bf v}_{(\mu;\,\lambda,i)}\otimes {\bf v}_{(\mu^\vee;\,\lambda^\vee,\lambda_1-\lambda_2-i)}.
\end{align}
Here the sequence of non-zero rational numbers $\{C_{(\mu;\,\lambda)}\}_{\mu \succ \lambda}$ is uniquely determined up to scalars and we normalize it so that
\[
C_{(\mu;\,(\mu_2,\mu_2))}=1.
\]
Then it is clear that
\begin{align}\label{E:scalar relation}
C_{(\mu;\,\lambda)} = C_{(\mu^\vee;\,\lambda^\vee)}
\end{align}
for $\mu \succ \lambda$.
We list in the following lemma some relations among the vectors ${\bf v}_{(\mu;\,\lambda,i)}$.
\begin{lemma}
Let $\mu \succ \lambda$ and $0 \leq i \leq \lambda_1-\lambda_2$. We have
\begin{align}
\begin{split}
{\bf v}_{(\mu;\,\lambda,0)} &=\tfrac{(\mu_1-\lambda_1)!(\lambda_2-\mu_3)!}{(\mu_1-\mu_2)!(\mu_2-\mu_3)!}\cdot E_{32}^{\mu_2-\lambda_2}\circ E_{13}^{\lambda_1-\mu_2}\cdot {\bf v}_{(\mu;\,(\mu_2,\mu_2),0)},\\
{\bf v}_{(\mu;\,\lambda,\lambda_1-\lambda_2)} &= (-1)^{\mu_2-\lambda_2}\cdot \tfrac{(\mu_1-\lambda_1)!(\lambda_2-\mu_3)!}{(\mu_1-\mu_2)!(\mu_2-\mu_3)!}\cdot E_{23}^{\lambda_1-\mu_2}\circ E_{31}^{\mu_2-\lambda_2}\cdot {\bf v}_{(\mu;\,(\mu_2,\mu_2),0)},
\end{split}\label{E:relation 1}\\
\begin{split}
c_\mu^{(3)}({\bf v}_{(\mu;\,\lambda,i)}) &= (-1)^{\mu_2-\lambda_2+i}\cdot{\bf v}_{(\mu^\vee;\,\lambda^\vee,\lambda_1-\lambda_2-i)},
\end{split}\label{E:conjugate equiv. 3}\\
\begin{split}
{\bf v}_{(\mu;\,\lambda,i)}
& = {\lambda_1-\lambda_2 \choose i}^{-1}\sum_{m= \max\{0,\lambda_2-\mu_2+i\}}^{\min\{i,\lambda_1-\mu_2\}}{\lambda_1-\mu_2 \choose m}{\mu_2-\lambda_2 \choose i-m}\\
&\quad\times {\bf v}_{(\mu;\, \lambda_1-\mu_2-m, m,\mu_1-\lambda_1,i-m, \mu_2-\lambda_2-i+m, \lambda_2-\mu_3)}.
\end{split}\label{E:U(2) to U(3)}
\end{align}
\end{lemma}

\begin{proof}
Relation (\ref{E:relation 1}) follows directly from (\ref{E:Lie algebra action 2}). Since $c_\mu^{(3)}\circ E_{ij} = -E_{ji}\circ c_\mu^{(3)}$, (\ref{E:conjugate equiv. 3}) is a consequence of (\ref{E:conjugate equiv. 1}) and (\ref{E:relation 1}).
Note that 
${\bf v}_{(\mu;\,\lambda,i)} = \frac{(\lambda_1-\lambda_2-i)!}{(\lambda_1-\lambda_2)!}\cdot E_{21}^i\cdot {\bf v}_{(\mu;\,\lambda,0)}$. By an induction argument using (\ref{E:Lie algebra action 2}), we then have (\ref{E:U(2) to U(3)}).

\end{proof}

\subsection{Cuspidal cohomology for $\GL_n$}\label{SS:cuspidal cohomology}
Let $\mu = \prod_{v \in S_\infty} \mu_v \in \prod_{v \in S_\infty} (X^+(T_n) \times X^+(T_n))$ with $\mu_v = (\mu_{\iota_v},\mu_{\overline{\iota}_v})$ for $v \in S_\infty$. 
Let $K_f$ be a neat open compact subgroup of $\GL_n(\A_{\K,f})$ and $\mathcal{S}_{n,K_f}$ be the locally symmetric space defined by
\[
\mathcal{S}_{n,K_f} = \GL_n(\K)\backslash \GL_n(\A_\K) / K_{n,\infty}K_f.
\]
The complex representation $M_{\mu,\C}^\vee$ defines a locally constant sheaf
$\mathcal{M}_{\mu,\C}^\vee$ of $\C$-vector spaces on $\mathcal{S}_{n,K_f}$ (cf.\,\cite[\S\,2.2.7]{HR2020}). For $q \in \Z_{\geq 0}$, we denote by
\[
H^q_B(\mathcal{S}_{n,K_f},{M}_{\mu,\C}^\vee),\quad H^q(\mathcal{S}_{n,K_f},\mathcal{M}_{\mu,\C}^\vee)
\]
the $q$-th singular cohomology of $\mathcal{S}_{n,K_f}$ with coefficients in ${M}_{\mu,\C}^\vee$ and the $q$-th sheaf cohomology of $\mathcal{M}^\vee_{\mu,\C}$, respectively.
The two cohomology groups are canonically isomorphic under the de Rham isomorphism.
For $\sigma \in {\rm Aut}(\C)$, the $\sigma$-linear isomorphism $M_{\mu,\C} \rightarrow M_{{}^\sigma\!\mu,\C}$ in (\ref{E:sigma linear}) naturally induces a $\sigma$-linear isomorphism $H^q_B(\mathcal{S}_{n,K_f},{M}_{\mu,\C}^\vee) \rightarrow H^q_B(\mathcal{S}_{n,K_f},{M}_{{}^\sigma\!\mu,\C}^\vee)$. This in term induces a $\sigma$-linear isomorphism
\[
T_{\sigma,K_f} : H^q(\mathcal{S}_{n,K_f},\mathcal{M}_{\mu,\C}^\vee) \longrightarrow H^q(\mathcal{S}_{n,K_f},\mathcal{M}_{{}^\sigma\!\mu,\C}^\vee)
\]
via the de Rham isomorphism.
Passing to the limit and define
\[
H^q(\mathcal{S}_n,\mathcal{M}_{\mu,\C}^\vee) = \varinjlim_{K_f}H^q(\mathcal{S}_{n,K_f},\mathcal{M}_{\mu,\C}^\vee).
\]
Then $H^q(\mathcal{S}_n,\mathcal{M}_{\mu,\C}^\vee)$ is canonically isomorphic to the relative Lie algebra cohomology
\[
H^q(\frak{g}_{n,\infty},K_{n,\infty};C^{\infty}(\GL_n(\K)\backslash \GL_n(\A_\K)) \otimes M_{\mu,\C}^\vee).
\]
We have the $\GL_n(\A_{\K,f})$-module structure on $H^q(\mathcal{S}_n,\mathcal{M}_{\mu,\C}^\vee)$ induced by the right translation of $\GL_n(\A_{\K,f})$ on $C^{\infty}(\GL_n(\K)\backslash \GL_n(\A_\K))$.
For $\sigma \in {\rm Aut}(\C)$, we have the $\sigma$-linear isomorphism 
\[
T_\sigma = \varinjlim_{K_f}T_{\sigma,K_f} : H^q(\mathcal{S}_n,\mathcal{M}_{\mu,\C}^\vee) \longrightarrow H^q(\mathcal{S}_n,\mathcal{M}_{\mu,\C}^\vee)
\]
which commutes with the $\GL_n(\A_{\K,f})$-action (cf.\,\cite[\S\,2.3.5]{Raghuram2016}).
The cuspidal cohomology is define to be the relative Lie algebra cohomology
\[
H_{\rm cusp}^q(\mathcal{S}_n,\mathcal{M}_{\mu,\C}^\vee) = H^q(\frak{g}_{n,\infty},K_{n,\infty};\mathcal{A}_{\rm cusp}(\GL_n(\K)\backslash \GL_n(\A_\K)) \otimes M_{\mu,\C}^\vee),
\]
where $\mathcal{A}_{\rm cusp}(\GL_n(\K)\backslash \GL_n(\A_\K))$ is the space of cusp forms on $\GL_n(\A_\K)$.
As explained in \cite[p.\,126]{Clozel1990}, the natural inclusion 
\[
\mathcal{A}_{\rm cusp}(\GL_n(\K)\backslash \GL_n(\A_\K)) \subset C^\infty(\GL_n(\K)\backslash \GL_n(\A_\K))
\]
induces a $\GL_n(\A_{\K,f})$-equivariant injective homomorphism
\[
H_{\rm cusp}^q(\mathcal{S}_n,\mathcal{M}_{\mu,\C}^\vee) \longrightarrow H^q(\mathcal{S}_n,\mathcal{M}_{\mu,\C}^\vee).
\]
We listed in the following theorem some general results on the cuspidal cohomology for $\GL_n(\A_\K)$. The proofs are contained in \cite[\S\,3.5 and Lemme 4.9]{Clozel1990}.

\begin{thm}\label{T:Clozel}
Let $\mu = \prod_{v \in S_\infty} \mu_v \in \prod_{v \in S_\infty} (X^+(T_n) \times X^+(T_n))$ and $q \in \Z_{\geq 0}$.
\begin{itemize}
\item[(1)] For $\sigma\in{\rm Aut}(\C)$, we have $T_\sigma\left(H_{\rm cusp}^q(\mathcal{S}_n,\mathcal{M}_{\mu,\C}^\vee)\right) = H_{\rm cusp}^q(\mathcal{S}_n,\mathcal{M}_{{}^\sigma\!\mu,\C}^\vee)$.
\item[(2)] $H_{\rm cusp}^q(\mathcal{S}_n,\mathcal{M}_{\mu,\C}^\vee)=0$ unless $\mu$ is pure and $d\cdot\tfrac{n(n-1)}{2}\leq q \leq d\cdot\tfrac{(n+2)(n-1)}{2}$.
\item[(3)] Suppose $\mu$ is pure with $\mu_{\iota_v} = (\mu_{1,v},\cdots,\mu_{n,v})$ and $\mu_{\overline{\iota}_v} = \mu_{\iota_v}^\vee + ({\sf w},\cdots,{\sf w})$ for $v \in S_\infty$. Let $\underline{\kappa}_v = (\kappa_{1,v},\cdots,\kappa_{n,v}) \in X^+(T_n)$ defined by 
\[
\kappa_{i,v} = 2\mu_{i,v}+n+1-2i-{\sf w}
\]
for $1 \leq i \leq n$ and $v \in S_\infty$.
Then
\[
H_{\rm cusp}^q(\mathcal{S}_n,\mathcal{M}_{\mu,\C}^\vee) = \bigoplus_\itPi H^q(\frak{g}_{n,\infty},K_{n,\infty};\itPi \otimes M_{\mu,\C}^\vee),
\]
where $\itPi$ runs through irreducible cuspidal automorphic representations of $\GL_n(\A_\K)$ such that
\[
\itPi_v = {\rm Ind}^{\GL_n(\C)}_{B_n(\C)}(\chi_{\kappa_{1,v}} \boxtimes \cdots \boxtimes \chi_{\kappa_{n,v}}) \otimes |\mbox{ }|_\C^{{\sf w}/2}
\]
for $v \in S_\infty$.
Here $\chi_\kappa$ is the character of $\C^\times$ defined by $\chi_\kappa(z) = (z/\overline{z})^{\kappa/2}$ for $\kappa \in \Z$.
\item[(4)] Let $\itPi$ be an irreducible cuspidal automorphic representation which contributes to the cuspidal cohomology of $\GL_n(\A_\K)$ with coefficients in $M^\vee_{\mu,\C}$. 
Then
\[
H^q(\frak{g}_{n,\infty},K_{n,\infty};\itPi \otimes M_{\mu,\C}^\vee) = \left(\itPi \otimes \wedge^q (\frak{g}_{n,\infty,\C} / \frak{k}_{n,\infty,\C})^* \otimes M_{\mu,\C}^\vee\right)^{K_{n,\infty}} \simeq \itPi_f
\]
for $q=d\cdot\tfrac{n(n-1)}{2},\, d\cdot\tfrac{(n+2)(n-1)}{2}$.
\end{itemize}
\end{thm}

We call $\tfrac{n(n-1)}{2}\cdot d$ and $\tfrac{(n+2)(n-1)}{2}\cdot d$ the bottom and top cuspidal degree, respectively.
Let $\itPi = \bigotimes_v \itPi_v$ be an irreducible cuspidal automorphic representation which contributes to the cuspidal cohomology of $\GL_n(\A_\K)$ with coefficients in $M^\vee_{\mu,\C}$. Let $\itPi_f  = \bigotimes_{v \nmid \infty} \itPi_v$ be the finite part of $\itPi$. For $\sigma \in {\rm Aut}(\C)$ and $v \in S_\infty$, let ${}^\sigma\!\itPi_v$ be the representation of $\GL_n(\K_v) = \GL_n(\C)$ defined by 
\[
{}^\sigma\!\itPi_v = \begin{cases}
\itPi_{\sigma^{-1}\circ v} & \mbox{ if $\sigma^{-1}\circ \iota_v = \iota_{\sigma^{-1}\circ v}$},\\
\itPi_{\sigma^{-1}\circ v}^c & \mbox{ if $\sigma^{-1}\circ \iota_v = \overline{\iota}_{\sigma^{-1}\circ v}$}.
\end{cases}
\]
Here $\itPi_v^c$ is the complex conjugation of $\itPi_v$.
Let ${}^\sigma\!\itPi$ be the irreducible admissible representation of $\GL_n(\A_\K)$ defined by 
\[
{}^\sigma\!\itPi = \bigotimes_{v \in S_\infty}{}^\sigma\!\itPi_v \otimes {}^\sigma\!\itPi_f.
\]
By Theorem \ref{T:Clozel} and the strong multiplicity one theorem for $\GL_n$, ${}^\sigma\!\itPi$ is automorphic cuspidal and contributes to the cuspidal cohomology of $\GL_n(\A_\K)$ with coefficients in $M^\vee_{{}^\sigma\!\mu,\C}$.
Moreover, ${}^\sigma\!\itPi_f$ appears in the cuspidal cohomology $H_{\rm cusp}^q(\mathcal{S}_3,\mathcal{M}_{{}^\sigma\!\mu,\C}^\vee)$ with multiplicity one for $q=\tfrac{n(n-1)}{2}\cdot d, \,\tfrac{(n+2)(n-1)}{2}\cdot d$, and the ${}^\sigma\!\itPi_f$-isotypic component is equal to $H^q(\frak{g}_{n,\infty},K_{n,\infty};{}^\sigma\!\itPi \otimes M_{{}^\sigma\!\mu,\C}^\vee)$. We have the $\sigma$-linear $\GL_n(\A_{\K,f})$-equivariant isomorphism 
\begin{align}\label{E:sigma linear map 1}
T_\sigma : H^q(\frak{g}_{n,\infty},K_{n,\infty}; \itPi \otimes M_{\mu,\C}^\vee) \longrightarrow H^q(\frak{g}_{n,\infty},K_{n,\infty};{}^\sigma\!\itPi \otimes M_{{}^\sigma\!\mu,\C}^\vee).
\end{align}
The rationality field $\Q(\itPi)$ of $\itPi$ is define to be the fixed field of
$\left\{\sigma \in {\rm Aut}(\C) \, \vert \, {}^\sigma\!\itPi = \itPi \right\}$ and is a number field.
We then have the $\Q(\itPi)$-rational structure on the cuspidal cohomology $H^q(\frak{g}_{n,\infty},K_{n,\infty}; \itPi \otimes M_{\mu,\C}^\vee)$ by taking the ${\rm Aut}(\C/\Q(\itPi))$-invariants:
\begin{align*}
&H^q(\frak{g}_{n,\infty},K_{n,\infty}; \itPi \otimes M_{\mu,\C}^\vee)^{{\rm Aut}(\C/\Q(\itPi))} \\
&= \left.\left\{\varphi \in H^q(\frak{g}_{n,\infty},K_{n,\infty}; \itPi \otimes M_{\mu,\C}^\vee) \,\right\vert\, T_\sigma \varphi=\varphi  \mbox{ for all }\sigma \in {\rm Aut}(\C/\Q(\itPi))\right\}.
\end{align*}

\subsection{Rational structures via the Whittaker models}\label{SS:Whittaker periods}

\subsubsection{Whittaker models}
Let $\itPi = \bigotimes_v \itPi_v$ be a cohomological irreducible cuspidal automorphic representation of $\GL_n(\A_\K)$. 
For each place $v$ of $\K$, let $\mathcal{W}(\itPi_v,\psi_{n,\K_v})$ be the space of Whittaker functions of $\itPi_v$ with respect to $\psi_{n,\K_v}$. 
Recall that $\mathcal{W}(\itPi_v,\psi_{n,\K_v})$ is contained in the space of locally constant (resp.\,smooth and moderate growth) functions $W : \GL_n(\K_v) \rightarrow \C$ which satisfy
\[
W(ug) = \psi_{n,\K_v}(u)W(g)
\]
for all $u \in N_n(\K_v)$ and $g \in \GL_n(\K_v)$ when $v$ is finite (resp.\,$v$ is archimedean).
We denote by $\rho$ the right translation action of $\GL_n(\K_v)$ (resp.\,$(\frak{g}_{n,v},K_{n,v})$) on $\mathcal{W}(\itPi_v,\psi_{n,\K_v})$ if $v$ is finite (resp.\,$v$ is archimedean).
Let 
\[
\mathcal{W}(\itPi_\infty,\psi_{n,\E_\infty}) = \bigotimes_{v \in S_\infty} \mathcal{W}(\itPi_v,\psi_{n,\C})
\]
and $\mathcal{W}(\itPi_f,\psi_{n,\K}^{(\infty)})$ be the space of Whittaker function of $\itPi_f = \bigotimes_{v \nmid \infty} \itPi_v$ with respect to $\psi_{n,\K}^{(\infty)} = \bigotimes_{v \nmid \infty} \psi_{n,\K_v}$. 
We have the $\GL_n(\A_{\K,f})$-equivariant isomorphism
\[
\bigotimes_{v\nmid \infty}\mathcal{W}(\itPi_v,\psi_{n,\K_v}) \longrightarrow \mathcal{W}(\itPi_f,\psi_{n,\K}^{(\infty)}),\quad \bigotimes_{v \nmid \infty} W_v \longmapsto \prod_{v \nmid \infty} W_v.
\]
Here we take the restricted tensor product on the left-hand side with respect to the normalized unramified Whittaker functions $W_v^\circ \in \mathcal{W}(\itPi_v,\psi_{n,\K_v})$, that is, the right $\GL_n(\frak{o}_{\K_v})$-invariant Whittaker function normalized so that $W_v^\circ(1)=1$ for all but finitely many finite places $v$ such that $\itPi_v$ is unramified.
For $\sigma \in {\rm Aut}(\C)$, let 
\begin{align}\label{E:sigma linear map 2}
t_\sigma : \mathcal{W}(\itPi_f,\psi_{n,\K}^{(\infty)}) \longrightarrow \mathcal{W}({}^\sigma\!\itPi_f,\psi_{n,\K}^{(\infty)})
\end{align}
be the $\sigma$-linear $\GL_n(\A_{\K,f})$-equivariant isomorphism defined by
\[
t_\sigma W(g) = \sigma ( W ({\rm diag}(u_\sigma^{-n+1},u_\sigma^{-n+2},\cdots,1)g))
\]
for $g \in \GL_n(\A_{\K,f})$. Here $u_\sigma \in  \prod_p \Z_p^\times \subset \A_{\E,f}^\times$ is the unique element depending on $\sigma$ such that $\sigma(\psi_\Q(x)) = \psi_\Q(u_\sigma x)$ for all $x \in \A_{\Q,f}$.
We then have the $\Q(\itPi)$-rational structure on $\mathcal{W}(\itPi_f,\psi_{n,\K}^{(\infty)})$ by taking the ${\rm Aut}(\C/\Q(\itPi))$-invariants 
\[
\mathcal{W}(\itPi_f,\psi_{n,\K}^{(\infty)})^{{\rm Aut}(\C/\Q(\itPi))} = \left.\left\{W \in \mathcal{W}(\itPi_f,\psi_{n,\K}^{(\infty)}) \,\right\vert\, t_\sigma W=W  \mbox{ for all }\sigma \in {\rm Aut}(\C/\Q(\itPi))\right\}.
\]
For $\varphi \in \itPi$, let $W_{\varphi,\psi_{n,\E}}$ be the Whittaker function of $\varphi$ with respect to $\psi_{n,\E}$ defined by
\[
W_{\varphi,\psi_{n,\E}}(g) = \int_{N_n(\E)\backslash N_n(\A_\E)}\varphi(ug)\overline{\psi_{n,\E}(u)}\,du^{\rm Tam}.
\]
Here $du^{\rm Tam}$ is the Tamagawa measure on $N_n(\A_\E)$. 
This defines a $\left((\frak{g}_{n,\infty},K_{n,\infty})\times\GL_n(\A_{\E,f})\right)$-equivariant isomorphism
\begin{align}\label{E:Whittaker isomorphism}
\itPi \longrightarrow \mathcal{W}(\itPi_\infty,\psi_{n,\E_\infty})\otimes \mathcal{W}(\itPi_f,\psi_{n,\E}^{(\infty)}),\quad \varphi \longmapsto W_{\varphi,\psi_{n,\F}}.
\end{align}

\subsubsection{Whittaker periods for $\GL_n$}\label{SS:Whittaker GL_n}

Let $\itPi$ be a cohomological irreducible cuspidal automorphic representation of $\GL_n(\A_\E)$. 
Let $\mu = \prod_{v \in S_\infty} \mu_v \in \prod_{v \in S_\infty} (X^+(T_n) \times X^+(T_n))$ such that $\itPi$ contributes to the cuspidal cohomology of $\GL_n(\A_\E)$ with coefficients in $M^\vee_{\mu,\C}$.
Let $*=b$ or $*=t$ and put 
\[
*_n=\begin{cases}
\tfrac{n(n-1)}{2} & \mbox{ if $*=b$},\\
\tfrac{(n+2)(n-1)}{2} & \mbox{ if $*=t$}.
\end{cases}
\]
The scripts $b$ and $t$ refer to the bottom and top degree cuspidal cohomology.
For each $v \in S_\infty$, we fix generators
\begin{align}\label{E:generators}
\begin{split}
[\itPi_v]_*  &\in H^{*_n}(\frak{g}_n,K_n;\mathcal{W}(\itPi_v,\psi_{n,\C}) \otimes M_{\mu_v,\C}^\vee) = \left(\mathcal{W}(\itPi_v,\psi_{n,\C})\otimes \wedge^{*_n}(\frak{g}_{n,\C} / \frak{k}_{n,\C})^* \otimes M_{\mu_v,\C}^\vee \right)^{K_n},\\
[\itPi_v^c]_*  &\in H^{*_n}(\frak{g}_n,K_n;\mathcal{W}(\itPi_v^c,\psi_{n,\C}) \otimes M_{\mu_v^c,\C}^\vee) = \left(\mathcal{W}(\itPi_v^c,\psi_{n,\C})\otimes \wedge^{*_n}(\frak{g}_{n,\C} / \frak{k}_{n,\C})^* \otimes M_{\mu_v^c,\C}^\vee \right)^{K_n}.
\end{split}
\end{align}
We assume $[\itPi_v]_* = [\itPi_v^c]_*$ if $\itPi_v = \itPi_v^c$.
We will specify the choices for $n=2$ and $n=3$ in \S\,\ref{SS:Whittaker period GL_2} and \S\,\ref{SS:Whittaker period GL_3}.
For $\sigma \in {\rm Aut}(\C)$, we then have the generator
\begin{align*}
[{}^\sigma\!\itPi_\infty]_* &= \left(\bigotimes_{v \in S_\infty,\,\sigma^{-1}\circ \iota_v = \iota_{\sigma^{-1}\circ v}}[\itPi_{\sigma^{-1}\circ v}]_*\right)\otimes \left(\bigotimes_{v \in S_\infty,\,\sigma^{-1}\circ \iota_v = \overline{\iota}_{\sigma^{-1}\circ v}}[\itPi_{\sigma^{-1}\circ v}^c]_*\right)\\
&\in H^{*_nd}\left(\frak{g}_{n,\infty},K_{n,\infty};\mathcal{W}({}^\sigma\!\itPi_\infty,\psi_{n,\E_\infty}) \otimes M_{{}^\sigma\!\mu,\C}^\vee\right).
\end{align*}
Let
\begin{align}\label{E:Galois equiv}
\mathcal{F}_{\sigma,*}: \mathcal{W}({}^\sigma\!\itPi_f,\psi_{n,\E}^{(\infty)}) \longrightarrow H^{*_nd}\left(\frak{g}_{n,\infty},K_{n,\infty};{}^\sigma\!\itPi \otimes M_{{}^\sigma\!\mu,\C}^\vee\right)
\end{align}
be the $\GL_n(\A_{\E,f})$-equivariant isomorphism defined as follows:
For $W \in \mathcal{W}({}^\sigma\!\itPi_f,\psi_{n,\E}^{(\infty)})$, we have
\[
 [{}^\sigma\!\itPi_\infty]_* \otimes W \in H^{*_nd}\left(\frak{g}_{n,\infty},K_{n,\infty};\mathcal{W}({}^\sigma\!\itPi_\infty,\psi_{n,\E_\infty})\otimes\mathcal{W}({}^\sigma\!\itPi_f,\psi_{n,\E}^{(\infty)}) \otimes M_{{}^\sigma\!\mu,\C}^\vee\right).
\]
Then $\mathcal{F}_{\sigma,*}(W)$ is the image of $ [{}^\sigma\!\itPi_\infty]_* \otimes W$ under the $\GL_n(\A_{\E,f})$-equivariant isomorphism induced by the inverse of the map (\ref{E:Whittaker isomorphism}) for ${}^\sigma\!\itPi$. We simply write $\mathcal{F}_{\sigma,*}=\mathcal{F}_*$ if $\sigma$ is the identity map.
Comparing the $\Q(\itPi)$-rational structures given by the cuspidal cohomology and by the Whittaker model, we have the following lemma/definition of the Whittaker periods (cf.\,\cite[Proposition 3.3]{RS2008}).

\begin{lemma}\label{L:Whittaker periods}
Let $*=b$ or $*=t$.
There exist $p^*({}^\sigma\!\itPi) \in \C^\times$ for each $\sigma \in {\rm Aut}(\C)$, unique up to $\Q({}^\sigma\!\itPi)^\times$, such that
\begin{align*}
\frac{\mathcal{F}_{\sigma,*}\left(\mathcal{W}({}^\sigma\!\itPi_f,\psi_{n,\E}^{(\infty)})^{{\rm Aut}(\C/\Q({}^\sigma\!\itPi))}\right)}{p^*({}^\sigma\!\itPi)} = H^{*_nd}(\frak{g}_{n,\infty},K_{n,\infty};{}^\sigma\!\itPi \otimes M_{{}^\sigma\!\mu,\C}^\vee)^{{\rm Aut}(\C/\Q({}^\sigma\!\itPi))}.
\end{align*}
Moreover, we can normalize the periods so that the diagram
\[
\begin{tikzcd}
\mathcal{W}(\itPi_f,\psi_{n,\F}^{(\infty)}) \arrow[rr, "\frac{1}{p^*(\itPi)}\cdot \mathcal{F}_*"] \arrow[d, "t_\sigma"] & &H^{*_nd}\left(\frak{g}_{n,\infty},K_{n,\infty};\itPi \otimes M_{\mu,\C}^\vee\right)\arrow[d, "T_\sigma"]\\
\mathcal{W}({}^\sigma\!\itPi_f,\psi_{n,\F}^{(\infty)})\arrow[rr, "\frac{1}{p^*({}^\sigma\!\itPi)}\cdot \mathcal{F}_{\sigma,*}"] & & H^{*_nd}\left(\frak{g}_{n,\infty},K_{n,\infty};{}^\sigma\!\itPi \otimes M_{{}^\sigma\!\mu,\C}^\vee\right)
\end{tikzcd}
\]
commutes.
Here $T_\sigma$ and $t_\sigma$ are the $\GL_n(\A_\E)$-equivariant $\sigma$-linear isomorphism in (\ref{E:sigma linear map 1}) and (\ref{E:sigma linear map 2}), respectively.
\end{lemma}
We call $p^b(\itPi)$ and $p^t(\itPi)$ the bottom degree Whittaker period and the top degree Whittaker period of $\itPi$, with respect to the choice of generators in (\ref{E:generators}). 

\subsubsection{Whittaker periods for $\GL_2$}\label{SS:Whittaker period GL_2}
Let $\itPi$ be a cohomological irreducible cuspidal automorphic representations of $\GL_2(\A_\K)$ with central characters $\omega_\itPi$.
We have $|\omega_\itPi| = |\mbox{ }|_{\A_\E}^{{\sf w}}$ for some ${\sf w} \in \Z$.
For $v \in S_\infty$, we have
\[
\itPi_v = {\rm Ind}_{B_2(\C)}^{\GL_2(\C)}(\chi_{\kappa_{1,v}}\boxtimes \chi_{\kappa_{2,v}})\otimes |\mbox{ }|_\C^{{\sf w}/2}
\]
for some $\underline{\kappa}_v = (\kappa_{1,v},\kappa_{2,v}) \in \Z^2$ such that
\begin{align*}
&\kappa_{1,v} > \kappa_{2,v},\quad \kappa_{1,v} \equiv \kappa_{2,v}\equiv 1+ {\sf w}\,({\rm mod}\,2).
\end{align*}
Note that $V_{\underline{\kappa}_v}^{(2)}$ is the minimal ${\rm U}(2)$-type of $\itPi_v$.
Then $\itPi$ contributes to the cuspidal cohomology of $\GL_2(\A_\K)$ with coefficients in $M_{\lambda,\C}^\vee$, where $\lambda = \prod_{v \in S_\infty}\lambda_v $ with 
\[
\lambda_{\iota_v} = \left(\tfrac{\kappa_{1,v}-1+{\sf w}}{2}, \tfrac{\kappa_{2,v}+1+{\sf w}}{2}\right),\quad
\lambda_{\overline{\iota}_v} = \lambda_{\iota_v}^\vee + ({\sf w},{\sf w}).
\]


We have $*_2=1$ (resp.\,$*_2=2$) if $*=b$ (resp.\,$*=t$).
Let $v \in S_\infty$. 
Note that $V_{\underline{\kappa}_v^\vee}^{(2)}$ appears in 
$
\wedge^{*_2}(\frak{g}_{2,\C} / \frak{k}_{2,\C})^* \otimes M_{\lambda_v,\C}^\vee
$
with multiplicity one.
Let
\[
\xi_{*,v}^{(2)} : V_{\underline{\kappa}_v^\vee}^{(2)}\longrightarrow \wedge^{*_2}(\frak{g}_{2,\C} / \frak{k}_{2,\C})^* \otimes M_{\lambda_v,\C}^\vee
\]
be the ${\rm U}(2)$-equivariant homomorphism such that
\begin{align}\label{E:U(2) equiv.}
\xi_{*,v}^{(2)}(x^{\kappa_{1,v}-\kappa_{2,v}}) =
\begin{cases}
Y_{(-1,1)}^*\otimes (x^{\lambda_{1,v}-\lambda_{2,v}}\otimes y^{\lambda_{1,v}-\lambda_{2,v}}) & \mbox{ if $*=b$},\\
Y_{(0,0)}^*\wedge Y_{(-1,1)}^*\otimes (x^{\lambda_{1,v}-\lambda_{2,v}}\otimes y^{\lambda_{1,v}-\lambda_{2,v}}) & \mbox{ if $*=t$}.
\end{cases}
\end{align}
Let
\begin{align}\label{E:U(2) equiv. 2}
\xi_{*,v}^{c,(2)} : V_{\underline{\kappa}_v}^{(2)}\longrightarrow  \wedge^{*_2}(\frak{g}_{2,\C} / \frak{k}_{2,\C})^* \otimes M_{\lambda_v^c,\C}^\vee
\end{align}
be the ${\rm U}(2)$-equivariant homomorphism defined so that the diagram
\[
\begin{tikzcd}
V_{\underline{\kappa}_v}^{(2)} \arrow[r, "\xi_{*,v}^{c,(2)}"] \arrow[d, "c_{\underline{\kappa}_v}^{(2)}"] & \wedge^{*_2}(\frak{g}_{2,\C} / \frak{k}_{2,\C})^* \otimes M_{\lambda_v^c,\C}^\vee\\
V_{\underline{\kappa}_v^\vee}^{(2)} \arrow[r, "\xi_{*,v}^{(2)}"]  & \wedge^{*_2}(\frak{g}_{2,\C} / \frak{k}_{2,\C})^* \otimes M_{\lambda_v,\C}^\vee\arrow[u, "\wedge^{*_2}(c_{(1,-1)}^{(2)})^* \otimes \,{\bf c}_{\lambda_v^\vee}^{(2)}"']
\end{tikzcd}
\]
commutes.
Here $c_{\underline{\kappa}_v}^{(2)}$, $c_{(1,-1)}^{(2)}$, and ${\bf c}_{\lambda_v^\vee}^{(2)}$ are the conjugate-equivariant isomorphisms defined in \S\,\ref{SS:rep GLn}.
We recall the following well-known result on the explicit formula for the Whittaker functions of $\itPi_v$ in the minimal ${\rm U}(2)$-type (cf.\,\cite[Proposition 6.1]{Cheng2020}).
\begin{lemma}
There exists a unique ${\rm U}(2)$-equivariant homomorphism 
\begin{align*}
V_{\underline{\kappa}_v}^{(2)} \longrightarrow \mathcal{W}(\itPi_v,\psi_{2,\C}),\quad x^{\kappa_{1,v}-\kappa_{2,v}-i}y^i \longmapsto W_{(\underline{\kappa}_v,{\sf w};\,i)}
\end{align*}
such that
\begin{align*}
W_{(\underline{\kappa}_v,{\sf w};\,i)}({\rm diag}(a_1a_2,a_2)) = (\sqrt{-1})^{\kappa_{1,v}-i}(a_1a_2^2)^{{\sf w}}\int_{L}\frac{ds}{2\pi\sqrt{-1}}\,a_1^{-s+1}\Gamma_\C(\tfrac{s+i}{2})\Gamma_\C(\tfrac{s+\kappa_{1,v}-\kappa_{2,v}-i}{2})
\end{align*}
for all $a_1,a_2>0$.
\end{lemma}
Following (\ref{E:generator 1}), let the classes $[\itPi_v]_*$ and $[\itPi_v^c]_*$ in (\ref{E:generators}) be defined by
\begin{align}\label{E:generator GL_2}
\begin{split}
[\itPi_v]_* & = \sum_{i=0}^{\kappa_{1,v}-\kappa_{2,v}}(-1)^i{\kappa_{1,v}-\kappa_{2,v} \choose i}\cdot W_{(\underline{\kappa}_v,{\sf w};\,i)}\otimes \xi_{*,v}^{(2)}(x^iy^{\kappa_{1,v}-\kappa_{2,v}-i}),\\
[\itPi_v^c]_* & = \sum_{i=0}^{\kappa_{1,v}-\kappa_{2,v}}(-1)^i{\kappa_{1,v}-\kappa_{2,v} \choose i}\cdot W_{(\underline{\kappa}_v^\vee,{\sf w};\,i)}\otimes (-1)^{{\sf w}}\cdot\xi_{*,v}^{c,(2)}(x^iy^{\kappa_{1,v}-\kappa_{2,v}-i}).
\end{split}
\end{align}
In the following lemma, we verify that the classes are well-defined.
\begin{lemma}
Let $v \in S_\infty$. Suppose $\itPi_v = \itPi_v^c$, then
$[\itPi_v]_* = [\itPi^c_v]_*$.
\end{lemma}

\begin{proof}
We drop the subscript $v$ for brevity.
Define the ${\rm U}(2)$-equivariant homomorphism
\[
\overline{\xi}_{*}^{(2)} : V_{\underline{\kappa}}^{(2)}\longrightarrow \wedge^{*_2}(\frak{g}_{2,\C} / \frak{k}_{2,\C})^* \otimes M_{\lambda^c,\C}^\vee
\]
as in (\ref{E:U(2) equiv.}) by replacing $\lambda$ by $\lambda^c = (\lambda_{\overline{\iota}},\lambda_{\iota})$.
We claim that
\begin{align*}
\overline{\xi}_{*}^{(2)} = (-1)^{\lambda_{1}-\lambda_{2}}\cdot {\xi}_{*}^{c,(2)}.
\end{align*}
Suppose the claim holds. Note that $\itPi_v = \itPi_v^c$ if and only if $\lambda=\lambda^c$. In this case, we have $\overline{\xi}_*^{(2)} = \xi_*^{(2)}$ and 
$(-1)^{\lambda_1-\lambda_2} = (-1)^{{\sf w}}$.
The assertion thus follows from the claim.
To prove the claim, it suffices to show that
\[
\overline{\xi}_{*}^{(2)}(x^{\kappa_1-\kappa_2}) = (-1)^{\lambda_{1}-\lambda_{2}}\cdot {\xi}_{*}^{c,(2)}(x^{\kappa_1-\kappa_2}).
\]
Note that
\[
\rho_{\underline{\kappa}^\vee}\left( \bp 0 & 1 \\ -1 & 0\ep\right)\cdot x^{\kappa_{1}-\kappa_{2}} = (-1)^{\kappa_1-\kappa_2}\cdot y^{\kappa_{1}-\kappa_{2}} = y^{\kappa_{1}-\kappa_{2}}
\]
and
\begin{align*}
&
{\rm Ad}\left( \bp 0 & 1 \\ -1 & 0\ep\right)Y_{(-1,1)}^*\otimes \rho_{\lambda^\vee}\left( \bp 0 & 1 \\ -1 & 0\ep\right)(x^{\lambda_{1}-\lambda_{2}}\otimes y^{\lambda_{1}-\lambda_{2}}) \\
&= (-1)^{\lambda_{1}-\lambda_{2}}\cdot Y_{(1,-1)}^* \otimes (y^{\lambda_{1}-\lambda_{2}}\otimes x^{\lambda_{1}-\lambda_{2}}),\\
&
{\rm Ad}\left( \bp 0 & 1 \\ -1 & 0\ep\right)(Y_{(0,0)}^*\wedge Y_{(-1,1)}^*)\otimes \rho_{\lambda^\vee}\left( \bp 0 & 1 \\ -1 & 0\ep\right)(x^{\lambda_{1}-\lambda_{2}}\otimes y^{\lambda_{1}-\lambda_{2}}) \\
&= (-1)^{1+\lambda_{1}-\lambda_{2}}\cdot Y_{(0,0)}^*\wedge Y_{(1,-1)}^* \otimes (y^{\lambda_{1}-\lambda_{2}}\otimes x^{\lambda_{1}-\lambda_{2}})
\end{align*}
Therefore,
\begin{align*}
\xi_{b}^{c,(2)} (x^{\kappa_{1}-\kappa_{2}}) &= (-1)^{\lambda_{1}-\lambda_{2}}\cdot (c_{(1,-1)}^{(2)})^*(Y_{(1,-1)}^*)\otimes {\bf c}_{\lambda^\vee}^{(2)}(y^{\lambda_{1}-\lambda_{2}}\otimes x^{\lambda_{1}-\lambda_{2}})\\
& = (-1)^{\lambda_{1}-\lambda_{2}}\cdot Y_{(-1,1)}^* \otimes (x^{\lambda_{1}-\lambda_{2}}\otimes y^{\lambda_{1}-\lambda_{2}}),\\
\xi_{t}^{c,(2)} (x^{\kappa_{1}-\kappa_{2}}) &= (-1)^{1+\lambda_{1}-\lambda_{2}}\cdot (c_{(1,-1)}^{(2)})^*(Y_{(0,0)}^*)\wedge(c_{(1,-1)}^{(2)})^*(Y_{(1,-1)}^*)\otimes {\bf c}_{\lambda^\vee}^{(2)}(y^{\lambda_{1}-\lambda_{2}}\otimes x^{\lambda_{1}-\lambda_{2}})\\
& = (-1)^{\lambda_{1}-\lambda_{2}}\cdot Y_{(0,0)}^*\wedge Y_{(-1,1)}^* \otimes (x^{\lambda_{1}-\lambda_{2}}\otimes y^{\lambda_{1}-\lambda_{2}})
\end{align*}
by (\ref{E:conjugate equiv. 1}) and (\ref{E:conjugate equiv. Lie 1}).
This completes the proof.
\end{proof}

\subsubsection{Whittaker periods for $\GL_3$}\label{SS:Whittaker period GL_3}
Let $\itSigma$ be a cohomological irreducible cuspidal automorphic representation of $\GL_3(\A_\K)$ with central character $\omega_\itSigma$.
We have $|\omega_\itSigma| = |\mbox{ }|_{\A_\E}^{3{\sf w}/2}$ for some ${\sf w} \in \Z$.
For $v \in S_\infty$, we have
\[
\itSigma_v = {\rm Ind}_{B_3(\C)}^{\GL_3(\C)}(\chi_{\ell_{1,v}} \boxtimes \chi_{\ell_{2,v}} \boxtimes \chi_{\ell_{3,v}})\otimes |\mbox{ }|_\C^{{\sf w}/2}
\]
for some $\underline{\ell}_v = (\ell_{1,v},\ell_{2,v},\ell_{3,v}) \in \Z^3$ such that
\begin{align*}
&\ell_{1,v} > \ell_{2,v} > \ell_{3,v},\quad \ell_{1,v} \equiv \ell_{2,v} \equiv \ell_{3,v} \equiv {\sf w} \,({\rm mod}\,2).
\end{align*}
Note that $V_{\underline{\ell}_v}^{(3)}$ is the minimal ${\rm U}(3)$-type of $\itSigma_v$.
Then $\itSigma$ contributes to the cuspidal cohomology of $\GL_3(\A_\K)$ with coefficients in $M_{\mu,\C}^\vee$, where $\mu = \prod_{v \in S_\infty}\mu_v$ with 
\[
\mu_{\iota_v} = \left( \tfrac{\ell_{1,v}-2+{\sf w}}{2}, \tfrac{\ell_{2,v}+{\sf w}}{2},\tfrac{\ell_{3,v}+2+{\sf w}}{2}\right)
,\quad \mu_{\overline{\iota}_v} = \mu_{\iota_v}^\vee + ({\sf w},{\sf w},{\sf w}).
\]

We have $*_3=3$ (resp.\,$*_3=5$) if $*=b$ (resp.\,$*=t$). 
Let $v \in S_\infty$.
Note that $V_{\underline{\ell}_v^\vee}^{(3)}$ appears in 
$
\wedge^{*_3}(\frak{g}_{3,\C} / \frak{k}_{3,\C})^* \otimes M_{\mu_v,\C}^\vee
$
with multiplicity one.
Let
\[
\xi_{*,v}^{(3)} : V_{\underline{\ell}_v^\vee}^{(3)}\longrightarrow \wedge^{*_3}(\frak{g}_{3,\C} / \frak{k}_{3,\C})^* \otimes M_{\mu_v,\C}^\vee
\]
be the ${\rm U}(3)$-equivariant homomorphism 
such that
\begin{align}\label{E:U(3) equiv.}
\xi_{*,v}^{(3)}({\bf v}_{(\underline{\ell}_v^\vee;\,(-\ell_{3,v},-\ell_{2,v}),0)}) =
\begin{cases}
X_{(-1,1,0)}^*\wedge X_{(0,-1,1)}^* \wedge X_{(-1,0,1)}^*\\
\otimes \left({\bf v}_{(\mu_{\iota_v}^\vee;\,(-\mu_{3,v},-\mu_{2,v}),0)}\otimes {\bf v}_{(\mu_{\overline{\iota}_v}^\vee;\,(\mu_{2,v}-{\sf w},\mu_{3,v}-{\sf w}),\mu_{2,v}-\mu_{3,v})}\right)
 & \mbox{ if $*=b$},\\
Z_{12}^*\wedge Z_{23}^* \wedge X_{(-1,1,0)}^*\wedge X_{(0,-1,1)}^* \wedge X_{(-1,0,1)}^*\\
\otimes \left({\bf v}_{(\mu_{\iota_v}^\vee;\,(-\mu_{3,v},-\mu_{2,v}),0)}\otimes {\bf v}_{(\mu_{\overline{\iota}_v}^\vee;\,(\mu_{2,v}-{\sf w},\mu_{3,v}-{\sf w}),\mu_{2,v}-\mu_{3,v})}\right)
 & \mbox{ if $*=t$}.
\end{cases}
\end{align}
Let
\begin{align}\label{E:U(3) equiv. 2}
\xi_{*,v}^{c,(3)} : V_{\underline{\ell}_v}^{(3)}\longrightarrow \wedge^{*_3}(\frak{g}_{3,\C} / \frak{k}_{3,\C})^* \otimes M_{\mu_v^c,\C}^\vee
\end{align}
be the ${\rm U}(3)$-equivariant homomorphism defined so that the diagrams
\[
\begin{tikzcd}
V_{\underline{\ell}_v}^{(3)} \arrow[r, "\xi_{*,v}^{c,(3)}"] \arrow[d, "c_{\underline{\ell}_v}^{(3)}"] & \wedge^{*_3}(\frak{g}_{3,\C} / \frak{k}_{3,\C})^* \otimes M_{\mu_v^c,\C}^\vee\\
V_{\underline{\ell}_v^\vee}^{(3)} \arrow[r, "\xi_{*,v}^{(3)}"]  & \wedge^{*_3}(\frak{g}_{3,\C} / \frak{k}_{3,\C})^* \otimes M_{\mu_v,\C}^\vee \arrow["\wedge^{*_3}(c_{(1,0,-1)}^{(3)})^* \otimes \,{\bf c}_{\mu_v^\vee}^{(3)}"', u]
\end{tikzcd}
\]
commute. Here $c_{\underline{\ell}_v}^{(3)}$, $c_{(1,0,-1)}^{(3)}$, and ${\bf c}_{\mu_v^\vee}^{(3)}$ are the conjugate-equivariant isomorphisms defined in \S\,\ref{SS:rep GLn}.
The following theorem of Hirano and Oda \cite[Theorem 7.7]{HO2009} (see also \cite[Theorem 4.1]{HIM2012}) is crucial to our computation of archimedean local integrals.
\begin{thm}[Hirano--Oda]\label{T:HO}
There exists a unique ${\rm U}(3)$-equivariant homomorphism 
\begin{align*}
V_{\underline{\ell}_v}^{(3)}\longrightarrow \mathcal{W}(\itSigma_v,\psi_{3,\C}),\quad {\bf v}_{(\underline{\ell}_v;\,\underline{n})} \longmapsto W_{(\underline{\ell}_v,{\sf w};\,\underline{n})}
\end{align*}
such that
\begin{align*}
&W_{(\underline{\ell}_v,{\sf w};\,\underline{n})}({\rm diag}(a_1a_2a_3,a_2a_3,a_3))\\
&= (\sqrt{-1})^{\ell_{1,v}+\ell_{2,v}+\ell_{3,v}+n_1-n_3+n_{12}-n_{23}}(a_1a_2^2a_3^3)^{{\sf w}}\\
&\times \int_{L_1}\frac{ds_1}{2\pi\sqrt{-1}}\int_{L_2}\frac{ds_2}{2\pi\sqrt{-1}}\,a_1^{-s_1+2}a_2^{-s_2+2}\\
&\times\frac{\Gamma_\C(\tfrac{s_1+n_2+n_3+n_{23}}{2})\Gamma_\C(\tfrac{s_1+n_1+n_{23}}{2})\Gamma_\C(\tfrac{s_1+n_1+n_{12}+n_{13}}{2})\Gamma_\C(\tfrac{s_2+n_1+n_2+n_{12}}{2})\Gamma_\C(\tfrac{s_2+n_3+n_{12}}{2})\Gamma_\C(\tfrac{s_2+n_3+n_{13}+n_{23}}{2})}{\Gamma_\C(\tfrac{s_1+s_2+n_1+n_3+n_{12}+n_{23}}{2})}
\end{align*}
for all $a_1,a_2,a_3>0$ and all $\underline{n} = (n_1,n_2,n_3,n_{23},n_{13},n_{12}) \in \Z_{\geq 0}^6$ with 
\[
n_1+n_2+n_3 = \mu_1-\mu_2,\quad n_{12}+n_{23}+n_{13} = \mu_2-\mu_3.
\]
\end{thm}
For $\underline{\ell}_v \succ \underline{\kappa}_v' = (\kappa_{1,v}',\kappa_{2,v}')$ and $0 \leq i \leq \kappa_{1,v}'-\kappa_{2,v}'$, let
\[
W_{(\underline{\ell}_v,{\sf w};\,\underline{\kappa}_v',i)} \in \mathcal{W}(\itSigma_v,\psi_{3,\C})
\]
be the image of ${\bf v}_{(\underline{\ell}_v;\,\underline{\kappa}_v',i)}$ under the homomorphism in Theorem \ref{T:HO}. 
Following (\ref{E:generator 2}), let the classes $[\itSigma_v]_*$ and $[\itSigma_v^c]_*$ in (\ref{E:generators}) be defined by
\begin{align}\label{E:generator GL_3}
\begin{split}
[\itSigma_v]_* & = \sum_{\underline{\ell}_v \succ \underline{\kappa}_v'}C_{(\underline{\ell}_v;\,\underline{\kappa}_v')}\sum_{i=0}^{\kappa_{1,v}'-\kappa_{2,v}'}(-1)^i{\kappa_{1,v}'-\kappa_{2,v}' \choose i}\cdot W_{(\underline{\ell}_v,{\sf w};\,\underline{\kappa}_v',i)}\otimes \xi_{*,v}^{(3)}({\bf v}_{(\underline{\ell}_v^\vee;\,(\underline{\kappa}_v')^\vee,\kappa_{1,v}'-\kappa_{2,v}'-i)}),\\
[\itSigma_v^c]_* & = \sum_{\underline{\ell}_v \succ \underline{\kappa}_v'}C_{(\underline{\ell}_v;\,\underline{\kappa}_v')}\sum_{i=0}^{\kappa_{1,v}'-\kappa_{2,v}'}(-1)^i{\kappa_{1,v}'-\kappa_{2,v}' \choose i}\cdot W_{(\underline{\ell}_v^\vee,{\sf w};\,(\underline{\kappa}_v')^\vee,i)}\\
&\quad\otimes  \begin{cases}
(-1)^{1+{\sf w}}\cdot\xi_{*,v}^{c,(3)}({\bf v}_{(\underline{\ell}_v;\,\underline{\kappa}_v',\kappa_{1,v}'-\kappa_{2,v}'-i)}) & \mbox{ if $*=b$},\\
(-1)^{{\sf w}}\cdot\xi_{*,v}^{c,(3)}({\bf v}_{(\underline{\ell}_v;\,\underline{\kappa}_v',\kappa_{1,v}'-\kappa_{2,v}'-i)}) & \mbox{ if $*=t$}.
\end{cases}
\end{split}
\end{align}
In the following lemma, we verify that the classes are well-defined.
\begin{lemma}
Let $v \in S_\infty$. Suppose $\itSigma_v = \itSigma_v^c$, then $[\itSigma_v]_* = [\itSigma_v^c]_*$.
\end{lemma}

\begin{proof}
We drop the subscript $v$ for brevity.
Define the ${\rm U}(3)$-equivariant homomorphism
\[
\overline{\xi}_{*}^{(3)} : V_{\underline{\ell}}^{(3)}\longrightarrow \wedge^{*_3}(\frak{g}_{3,\C} / \frak{k}_{3,\C})^* \otimes M_{\mu^c,\C}^\vee
\]
as in (\ref{E:U(3) equiv.}) by replacing $\mu$ by $\mu^c = (\mu_{\overline{\iota}},\mu_{\iota})$.
We claim that
\[
\overline{\xi}_{*}^{(3)} =
\begin{cases}
(-1)^{1+\mu_{1}-\mu_{3}}\cdot \xi_{*}^{c,(3)} & \mbox{ if $*=b$},\\
(-1)^{\mu_1-\mu_3} \cdot \xi_*^{c,(3)} & \mbox{ if $*=t$}.
\end{cases}
\]
Suppose the claim holds. Note that $\itSigma_v = \itSigma_v^c$ if and only if $\mu=\mu^c$. In this case, we have $\overline{\xi}_*^{(3)} = \xi_*^{(3)}$ and 
$(-1)^{\mu_1-\mu_3} = (-1)^{{\sf w}}$.
The assertion thus follows from the claim.
To prove the claim, it suffices to show that
\[
\overline{\xi}_{*}^{(3)}({\bf v}_{(\underline{\ell};\,(\ell_{1},\ell_{2}),0)}) = 
\begin{cases}
(-1)^{1+\mu_{1}-\mu_{3}}\cdot \xi_{*}^{c,(3)}({\bf v}_{(\underline{\ell};\,(\ell_{1},\ell_{2}),0)}) & \mbox{ if $*=b$},\\
(-1)^{\mu_1-\mu_3} \cdot \xi_*^{c,(3)}({\bf v}_{(\underline{\ell};\,(\ell_{1},\ell_{2}),0)}) & \mbox{ if $*=t$}.
\end{cases}
\]
Note that
\begin{align*}
c_{\underline{\ell}}^{(3)}({\bf v}_{(\underline{\ell};\,(\ell_{1},\ell_{2}),0)}) 
& = {\bf v}_{(\underline{\ell}^\vee;\,(-\ell_2,-\ell_1),\ell_1-\ell_2)}\\
& = (-1)^{\ell_{1}}\cdot\rho_{\underline{\ell}^\vee}\left( \bp 0&0&1\\0&1&0\\1&0&0\ep\right)\cdot {\bf v}_{(\underline{\ell}^\vee;\,(-\ell_{3},-\ell_{2}),0)}\\
& = (-1)^{{\sf w}}\cdot\rho_{\underline{\ell}^\vee}\left( \bp 0&0&1\\0&1&0\\1&0&0\ep\right)\cdot {\bf v}_{(\underline{\ell}^\vee;\,(-\ell_{3},-\ell_{2}),0)}
\end{align*}
by (\ref{E:conjugate equiv. 3}), and
\begin{align*}
&{\rm Ad}\left( \bp 0&0&1\\0&1&0\\1&0&0\ep \right) \left(X_{(-1,1,0)}^*\wedge X_{(0,-1,1)}^* \wedge X_{(-1,0,1)}^*\right)\\
&\otimes \rho_{\lambda^\vee}\left( \bp 0&0&1\\0&1&0\\1&0&0\ep \right)\left({\bf v}_{(\mu_{\iota}^\vee;\,(-\mu_{3},-\mu_{2}),0)}\otimes {\bf v}_{(\mu_{\overline{\iota}}^\vee;\,(\mu_{2}-{\sf w},\mu_{3}-{\sf w}),\mu_{2}-\mu_{3})}\right)\\
&= (-1)^{1+{\sf w}+\mu_{1}-\mu_{3}}\cdot X_{(0,1,-1)}^*\wedge X_{(1,-1,0)}^* \wedge X_{(1,0,-1)}^* \\
&\otimes \left({\bf v}_{(\mu_{\iota}^\vee;\,(-\mu_{2},-\mu_{1}),\mu_{1}-\mu_{2})}\otimes {\bf v}_{(\mu_{\overline{\iota}}^\vee;\,(\mu_{1}-{\sf w},\mu_{2}-{\sf w}),0)}\right),\\
&{\rm Ad}\left( \bp 0&0&1\\0&1&0\\1&0&0\ep \right) \left(Z_{12}^* \wedge Z_{23}^* \wedge X_{(-1,1,0)}^*\wedge X_{(0,-1,1)}^* \wedge X_{(-1,0,1)}^*\right)\\
&\otimes \rho_{\lambda^\vee}\left( \bp 0&0&1\\0&1&0\\1&0&0\ep \right)\left({\bf v}_{(\mu_{\iota}^\vee;\,(-\mu_{3},-\mu_{2}),0)}\otimes {\bf v}_{(\mu_{\overline{\iota}}^\vee;\,(\mu_{2}-{\sf w},\mu_{3}-{\sf w}),\mu_{2}-\mu_{3})}\right)\\
&= (-1)^{{\sf w}+\mu_{1}-\mu_{3}}\cdot Z_{12}^* \wedge Z_{23}^* \wedge X_{(0,1,-1)}^*\wedge X_{(1,-1,0)}^* \wedge X_{(1,0,-1)}^*\\
&\otimes \left({\bf v}_{(\mu_{\iota}^\vee;\,(-\mu_{2},-\mu_{1}),\mu_{1}-\mu_{2})}\otimes {\bf v}_{(\mu_{\overline{\iota}}^\vee;\,(\mu_{1}-{\sf w},\mu_{2}-{\sf w}),0)}\right).
\end{align*}
Thus 
\begin{align*}
\overline{\xi}_{b}^{(3)}({\bf v}_{(\underline{\ell};\,(\ell_{1},\ell_{2}),0)}) &= (-1)^{1+\mu_{1}-\mu_{3}}\cdot \wedge^3(c_{(1,0,-1)}^{(3)})^*\left(X_{(0,1,-1)}^*\wedge X_{(1,-1,0)}^* \wedge X_{(1,0,-1)}^*\right) \\
&\otimes {\bf c}_{\mu^\vee}^{(3)}\left({\bf v}_{(\mu_{\iota}^\vee;\,(-\mu_{2},-\mu_{1}),\mu_{1}-\mu_{2})}\otimes {\bf v}_{(\mu_{\overline{\iota}}^\vee;\,(\mu_{1}-{\sf w},\mu_{2}-{\sf w}),0)}\right)\\
& = (-1)^{1+\mu_{1}-\mu_{3}}\cdot X_{(-1,1,0)}^*\wedge X_{(0,-1,1)}^* \wedge X_{(-1,0,1)}^*\\
&\otimes \left({\bf v}_{(\mu_{\overline{\iota}}^\vee;\,(\mu_{1}-{\sf w},\mu_{2}-{\sf w}),0)}\otimes{\bf v}_{(\mu_{\iota}^\vee;\,(-\mu_{2},-\mu_{1}),\mu_{1}-\mu_{2})}\right),\\
\overline{\xi}_{t}^{(3)}({\bf v}_{(\underline{\ell};\,(\ell_{1},\ell_{2}),0)}) &= (-1)^{\mu_{1}-\mu_{3}}\cdot \wedge^5(c_{(1,0,-1)}^{(3)})^*\left(Z_{12}^* \wedge Z_{23}^*\wedge X_{(0,1,-1)}^*\wedge X_{(1,-1,0)}^* \wedge X_{(1,0,-1)}^*\right) \\
&\otimes {\bf c}_{\mu^\vee}^{(3)}\left({\bf v}_{(\mu_{\iota}^\vee;\,(-\mu_{2},-\mu_{1}),\mu_{1}-\mu_{2})}\otimes {\bf v}_{(\mu_{\overline{\iota}}^\vee;\,(\mu_{1}-{\sf w},\mu_{2}-{\sf w}),0)}\right)\\
& = (-1)^{\mu_{1}-\mu_{3}}\cdot Z_{12}^*\wedge Z_{23}^*\wedge X_{(-1,1,0)}^*\wedge X_{(0,-1,1)}^* \wedge X_{(-1,0,1)}^*\\
&\otimes \left({\bf v}_{(\mu_{\overline{\iota}}^\vee;\,(\mu_{1}-{\sf w},\mu_{2}-{\sf w}),0)}\otimes{\bf v}_{(\mu_{\iota}^\vee;\,(-\mu_{2},-\mu_{1}),\mu_{1}-\mu_{2})}\right)
\end{align*}
by (\ref{E:conjugate equiv. Lie 2}) and (\ref{E:conjugate equiv. Lie 3}).
This completes the proof.
\end{proof}

\section{Algebraicity of the Rankin--Selberg $L$-functions for $\GL_3 \times \GL_2$}\label{S:RS}

Let $\E$ be a totally imaginary quadratic extension of $\F$.

\subsection{Algebraicity for $\GL_3 \times \GL_2$}

Let $\itPi$ and $\itSigma$ be cohomological irreducible cuspidal automorphic representations of $\GL_2(\A_\K)$ and $\GL_3(\A_\K)$ with central characters $\omega_\itPi$ and $\omega_\itSigma$, respectively.
We have $|\omega_\itPi| = |\mbox{ }|_{\A_\E}^{{\sf w}(\itPi)}$ and $|\omega_\itSigma| = |\mbox{ }|_{\A_\E}^{3{\sf w}(\itSigma)/2}$ for some ${\sf w}(\itPi),{\sf w}(\itSigma) \in \Z$.
For $v \in S_\infty$, we have
\[
\itPi_v = {\rm Ind}_{B_2(\C)}^{\GL_2(\C)}(\chi_{\kappa_{1,v}}\boxtimes \chi_{\kappa_{2,v}})\otimes |\mbox{ }|_\C^{{\sf w}(\itPi)/2},\quad \itSigma_v = {\rm Ind}_{B_3(\C)}^{\GL_3(\C)}(\chi_{\ell_{1,v}} \boxtimes \chi_{\ell_{2,v}} \boxtimes \chi_{\ell_{3,v}})\otimes |\mbox{ }|_\C^{{\sf w}(\itSigma)/2}
\]
for some $\underline{\kappa}_v = (\kappa_{1,v},\kappa_{2,v}) \in \Z^2$ and $\underline{\ell}_v = (\ell_{1,v},\ell_{2,v},\ell_{3,v}) \in \Z^3$ such that
\begin{align*}
\kappa_{1,v} > \kappa_{2,v},\quad &\kappa_{1,v} \equiv \kappa_{2,v}\equiv 1+ {\sf w}(\itPi)\,({\rm mod}\,2),\\
\ell_{1,v} > \ell_{2,v} > \ell_{3,v},\quad &\ell_{1,v} \equiv \ell_{2,v} \equiv \ell_{3,v} \equiv {\sf w}(\itSigma) \,({\rm mod}\,2).
\end{align*}
Then $\itPi$ and $\itSigma$ contribute to the cuspidal cohomology of $\GL_2(\A_\K)$ and $\GL_3(\A_\K)$ with coefficients in $M_{\lambda,\C}^\vee$ and $M_{\mu,\C}^\vee$, respectively, where $\lambda = \prod_{v \in S_\infty}\lambda_v $ and $\mu = \prod_{v \in S_\infty}\mu_v$ with 
\begin{align*}
\lambda_{\iota_v} = \left(\tfrac{\kappa_{1,v}-1+{\sf w}(\itPi)}{2}, \tfrac{\kappa_{2,v}+1+{\sf w}(\itPi)}{2}\right),\quad
\mu_{\iota_v} = \left( \tfrac{\ell_{1,v}-2+{\sf w}(\itSigma)}{2}, \tfrac{\ell_{2,v}+{\sf w}(\itSigma)}{2},\tfrac{\ell_{3,v}+2+{\sf w}(\itSigma)}{2}\right)
\end{align*}
and 
\[
\lambda_{\overline{\iota}_v} = \lambda_{\iota_v}^\vee + ({\sf w}(\itPi),{\sf w}(\itPi)),\quad \mu_{\overline{\iota}_v} = \mu_{\iota_v}^\vee + ({\sf w}(\itSigma),{\sf w}(\itSigma),{\sf w}(\itSigma)).
\]

Let $\delta : \GL_2 \rightarrow \GL_3$ be the embedding defined by $\delta(g) = \bp g & 0 \\ 0 & 1\ep$. 
Let 
\[
s: \wedge^3 (\frak{g}_{3,\C} / \frak{k}_{3,\C})^* \times (\frak{g}_{2,\C} / \frak{k}_{2,\C})^* \longrightarrow \C
\]
be the bilinear homomorphism defined by
\[
\delta^*X_1^* \wedge \delta^* X_2^* \wedge \delta^*X_3^* \wedge {\rm pr}(Y^*) = s(X_1^* \wedge X_2^* \wedge X_3^*, \, Y^*)\cdot e_{11}^*\wedge e_{22}^* \wedge e_{12}^* \wedge \sqrt{-1}\,e_{12}^*.
\]
Here ${\rm pr}: (\frak{g}_{2,\C} / \frak{k}_{2,\C})^* \rightarrow (\frak{g}_{2,\C} / \frak{u}(2)_\C)^*$ is the natural surjection and $\delta^* : (\frak{g}_{3,\C} / \frak{k}_{3,\C})^*  \rightarrow (\frak{g}_{2,\C} / \frak{u}(2)_\C)^*$ is the homomorphism induced by the embedding $\delta$.
Let $v \in S_\infty$.
For $W_1 \in \mathcal{W}(\itSigma_v, \psi_{3,\C})$ and $W_2 \in \mathcal{W}(\itPi_v, \psi_{2,\C})$, we define the local zeta integral
\[
Z_v(s,W_1,W_2) = \int_{N_2(\C)\backslash \GL_2(\C)}W_1(\delta(g))W_2({\rm diag}(-1,1)g)|\det(g)|_\C^{s-1/2}\,dg.
\]
The integral converges absolutely for ${\rm Re}(s)$ sufficiently large and admits meromorphic continuation to $s \in \C$.
Moreover, the ratio
\[
\frac{Z_v(s,W_1,W_2)}{L(s,\itSigma_v \times \itPi_v)}
\]
is entire (cf.\,\cite{JS1990} and \cite{Jacquet2009}).
Let $m+\tfrac{1}{2} \in \Z+\tfrac{1}{2}$. 
Under the assumption that
\[
\ell_{1,v} > -\kappa_{2,v} > \ell_{2,v} > -\kappa_{1,v} > \ell_{3,v},
\]
it is easy to verify that $m+\tfrac{1}{2}$ is not a pole of neither $L(s,\itSigma_v \times \itPi_v)$ nor $L(1-s,\itSigma_v^\vee \times \itPi_v^\vee)$ if and only if 
\[
{\rm Hom}_{\GL_2}(M_{\lambda_{\iota_v}+m}^\vee,M_{\mu_{\iota_v}})\neq 0,\quad {\rm Hom}_{\GL_2}(M_{\lambda_{\overline{\iota}_v}+m}^\vee,M_{\mu_{\overline{\iota}_v}})\neq 0.
\]
In this case, we fix non-zero homomorphisms
\[
\iota_{m,v} \in {\rm Hom}_{\GL_2(\C)}(M_{\lambda_{\iota_v}+m,\C}^\vee,M_{\mu_{\iota_v},\C}),\quad \overline{\iota}_{m,v} \in {\rm Hom}_{\GL_2(\C)}(M_{\lambda_{\overline{\iota}_v+m},\C}^\vee,M_{\mu_{\overline{\iota}_v},\C})
\]
defined over $\Q$.
We also fix non-zero $\GL_3(\C)$-equivariant bilinear pairings 
\[
\<\cdot,\cdot\>_{\mu_{\iota_v},\C} : M_{\mu_{\iota_v},\C}^\vee \times M_{\mu_{\iota_v},\C} \longrightarrow \C,\quad \<\cdot,\cdot\>_{\mu_{\overline{\iota}_v},\C} : M_{\mu_{\overline{\iota}_v},\C}^\vee \times M_{\mu_{\overline{\iota}_v},\C} \longrightarrow \C
\]
defined over $\Q$.
Note that $\rho_{\lambda_{\iota_v}+m}^\vee = \rho_{\lambda_{\iota_v}}^\vee \otimes \det^{-m}$, thus $M_{\lambda_{\iota_v}+m}^\vee = M_{\lambda_{\iota_v}}^\vee$ as vector spaces by our convention in \S\,\ref{SS:rep GL2}.
Similarly $M_{\lambda_{\overline{\iota}_v}+m}^\vee = M_{\lambda_{\overline{\iota}_v}}^\vee$ as vector spaces.
Let 
\begin{align*}
\<\cdot,\cdot\>_{m,v}: \left(\mathcal{W}(\itSigma_v,\psi_{3,\C}) \otimes \wedge^3 (\frak{g}_{3,\C} / \frak{k}_{3,\C})^* \otimes M_{\mu_v,\C}^\vee\right)\times\left(\mathcal{W}(\itPi_v,\psi_{2,\C}) \otimes  (\frak{g}_{2,\C} / \frak{k}_{2,\C})^* \otimes M_{\lambda_v,\C}^\vee\right) &\longrightarrow \C,\\
\<\cdot,\cdot\>_{m,v}^c: \left(\mathcal{W}(\itSigma_v^c,\psi_{3,\C}) \otimes \wedge^3 (\frak{g}_{3,\C} / \frak{k}_{3,\C})^* \otimes M_{\mu_v^c,\C}^\vee\right)\times\left(\mathcal{W}(\itPi_v^c,\psi_{2,\C}) \otimes  (\frak{g}_{2,\C} / \frak{k}_{2,\C})^* \otimes M_{\lambda_v^c,\C}^\vee\right) &\longrightarrow \C
\end{align*}
be the bilinear homomorphisms defined by
\begin{align*}
&\<W_1 \otimes X_1^*\wedge X_2^* \wedge X_3^* \otimes ({\bf v}_{\iota_v} \otimes {\bf v}_{\overline{\iota}_v}),\, W_2 \otimes Y^* \otimes ({\bf u}_{\iota_v} \otimes {\bf u}_{\overline{\iota}_v})\>_{m,v}\\
& = Z_v(m+\tfrac{1}{2},W_1,W_2) \cdot s(X_1^* \wedge X_2^* \wedge X_3^*, \, Y^*) \cdot \<{\bf v}_{\iota_v},\iota_{m,v}({\bf u}_{\iota_v})\>_{\mu_{\iota_v},\C}\cdot \<{\bf v}_{\overline{\iota}_v},\overline{\iota}_{m,v}({\bf u}_{\overline{\iota}_v})\>_{\mu_{\overline{\iota}_v},\C},\\
&\<W_1 \otimes X_1^*\wedge X_2^* \wedge X_3^* \otimes ({\bf v}_{\overline{\iota}_v}\otimes {\bf v}_{\iota_v}),\, W_2 \otimes Y^* \otimes ({\bf u}_{\overline{\iota}_v}\otimes {\bf u}_{\iota_v})\>_{m,v}^c\\
& = Z_v(m+\tfrac{1}{2},W_1,W_2) \cdot s(X_1^* \wedge X_2^* \wedge X_3^*, \, Y^*) \cdot \<{\bf v}_{\overline{\iota}_v},\overline{\iota}_{m,v}({\bf u}_{\overline{\iota}_v})\>_{\mu_{\overline{\iota}_v},\C}\cdot\<{\bf v}_{\iota_v},\iota_{m,v}({\bf u}_{\iota_v})\>_{\mu_{\iota_v},\C}.
\end{align*}
Assume that
\[
\ell_{1,v} > -\kappa_{2,v} > \ell_{2,v} > -\kappa_{1,v} > \ell_{3,v}
\]
for all $v \in S_\infty$ and $m+\tfrac{1}{2}$ is critical for $L(s,\itSigma \times \itPi)$.
For $\sigma \in {\rm Aut}(\C)$, we then have the bilinear homomorphism
\[
\<\cdot,\cdot\>_{\sigma,m,\infty} = \left(\bigotimes_{v \in S_\infty,\,\sigma^{-1}\circ \iota_v = \iota_{\sigma^{-1}\circ v}}\<\cdot,\cdot\>_{m,\sigma^{-1} \circ v} \right) \otimes \left(\bigotimes_{v \in S_\infty,\,\sigma^{-1}\circ \iota_v = \overline{\iota}_{\sigma^{-1}\circ v}}\<\cdot,\cdot\>_{m,\sigma^{-1} \circ v}^c \right)
\]
from 
\begin{align*}
\left(\mathcal{W}({}^\sigma\!\itSigma_\infty,\psi_{3,\E_\infty}) \otimes \wedge^{3d} (\frak{g}_{3,\infty,\C} / \frak{k}_{3,\infty,\C})^* \otimes M_{{}^\sigma\!\mu,\C}^\vee\right)\times\left(\mathcal{W}({}^\sigma\!\itPi_\infty,\psi_{2,\E_\infty}) \otimes  \wedge^d(\frak{g}_{2,\infty,\C} / \frak{k}_{2,\infty,\C})^* \otimes M_{{}^\sigma\!\lambda,\C}^\vee\right)
\end{align*}
to $\C$.
We simply write $\<\cdot,\cdot\>_{\sigma,m,\infty} = \<\cdot,\cdot\>_{m,\infty}$ when $\sigma$ is the identity map.

We have the following result due to Raghuram \cite{Raghuram2016} and Sun \cite{Sun2017} on the algebraicity of the critical values of $L(s,\itSigma \times \itPi)$ in terms of the bottom degree Whittaker periods of $\itSigma$ and $\itPi$.
Note that we have incorporated the period relation proved by Raghuram--Shahidi \cite{RS2008} into the formula.
\begin{thm}\label{T:Raghuram}
Assume $\ell_{1,v} > -\kappa_{2,v} > \ell_{2,v} > -\kappa_{1,v} > \ell_{3,v}$ for all $v \in S_\infty$.
Let $m+\tfrac{1}{2}$ be critical for $L(s,\itSigma \times \itPi)$. 
\begin{itemize}
\item[(1)]{\rm (Raghuram)}
We have
\begin{align*}
&\sigma\left(\frac{L^{(\infty)}(m+\tfrac{1}{2},\itSigma \times \itPi)}{|D_\K|^{1/2}\cdot G(\omega_\itPi)\cdot p^b(\itSigma)\cdot p^b(\itPi)}\cdot \<[\itSigma_\infty]_b , [\itPi_\infty]_b\>_{m,\infty}\right)\\
& = \frac{L^{(\infty)}(m+\tfrac{1}{2},{}^\sigma\!\itSigma \times {}^\sigma\!\itPi)}{|D_\K|^{1/2}\cdot G({}^\sigma\!\omega_\itPi)\cdot p^b({}^\sigma\!\itSigma)\cdot p^b({}^\sigma\!\itPi)}\cdot \<[{}^\sigma\!\itSigma_\infty]_b , [{}^\sigma\!\itPi_\infty]_b\>_{\sigma,m,\infty}
\end{align*}
for all $\sigma \in {\rm Aut}(\C)$.
\item[(2)]{\rm (Sun)} The archimedean factor $\<[\itSigma_\infty]_b , [\itPi_\infty]_b\>_{m,\infty}$ is non-zero.
\end{itemize}
\end{thm}

\begin{rmk}
The factor $|D_\E|^{1/2}$ is due to the comparison between Haar measures on $\A_\E$ that
\[
dx = |D_\E|^{1/2}\cdot dx^{\rm Tam},
\]
where $dx = \prod_v dx_v$ is the standard measure in \S\,\ref{SS:measure} and $dx^{\rm Tam}$ is the Tamagawa measure.
\end{rmk}
In the following lemma, we compute the archimedean local zeta integrals appearing in the archimedean factor $\<[\itSigma_\infty]_b , [\itPi_\infty]_b\>_{m,\infty}$.
The result was announced in \cite[Theorem 4.1]{HIM2016} where the choice of Whittaker functions are given in \S\,4.2 therein. We give a proof for the sake of completeness.
\begin{lemma}\label{L:RS}
Let $v \in S_\infty$.
Assume $\ell_{1,v} \geq -\kappa_{2,v} \geq \ell_{2,v} \geq -\kappa_{1,v} \geq \ell_{3,v}$. We have
\begin{align*}
&Z_v(s,W_{(\underline{\ell}_v,{\sf w}(\itSigma);\,\underline{\kappa}_v^\vee,0)},W_{(\underline{\kappa}_v,{\sf w}(\itPi);\,\kappa_{1,v}-\kappa_{2,v})}) \\
&= Z_v(s,W_{(\underline{\ell}_v^\vee,{\sf w}(\itSigma);\,\underline{\kappa}_v,\kappa_{1,v}-\kappa_{2,v},)},W_{(\underline{\kappa}_v^\vee,{\sf w}(\itPi);\,0)})\\
&= 4(\kappa_{1,v}-\kappa_{2,v}+1)^{-1}(\sqrt{-1})^{\kappa_{1,v}+\kappa_{2,v}}\cdot L(s,\itSigma_v\times\itPi_v).
\end{align*}
\end{lemma}

\begin{proof}
We drop the subscript $v$ for brevity. 
We may assume ${\sf w}(\itSigma) = {\sf w}(\itPi)=0$ after replacing $\itSigma_v$ and $\itPi_v$ by $\itSigma_v \otimes |\mbox{ }|_\C^{-{\sf w}(\itSigma)/2}$ and $\itPi_v \otimes |\mbox{ }|_\C^{-{\sf w}(\itPi)/2}$, respectively. 
Write $W_{(\underline{\ell};\,n_1,n_2,n_3,n_{23},n_{13},n_{12})}=W_{(\underline{\ell},{\sf w}(\itSigma);\,n_1,n_2,n_3,n_{23},n_{13},n_{12})}$ and $W_{(\underline{\kappa};\,i)} = W_{(\underline{\kappa},{\sf w}(\itPi);\,i)}$.
Since the assertion is an equality for meromorphic functions in $s \in \C$, it suffices to prove the equality for ${\rm Re}(s)$ sufficiently large.
We assume
\[
{\rm Re}(s) > \max\{\tfrac{2\ell_2+\ell_3+2\kappa_1+\kappa_2}{2},0\} - \min\{\tfrac{2\ell_2+\kappa_1+\kappa_2}{2},0\}.
\]
Let $L_1$ be the vertical path from $c_1-\sqrt{-1}\,\infty$ to $c_1+\sqrt{-1}\,\infty$ for some
\[
\max\{2\ell_2+\ell_3+2\kappa_1+\kappa_2,0\} < c_1 < 2{\rm Re}(s) + \min\{2\ell_2+\kappa_1+\kappa_2,0\}.
\]
Let $L_2=L_3$ be the vertical path from $1-\sqrt{-1}\,\infty$ to $1+\sqrt{-1}\,\infty$.
For $k = \bp \alpha & \beta \\ -\overline{\beta} & \overline{\alpha}\ep \in {\rm SU}(2)$, we have
\begin{align*}
\rho\left(\bp k & 0 \\ 0 & 1\ep\right)W_{(\underline{\ell};\,\underline{\kappa}^\vee,0)} &= \sum_{i=0}^{-\ell_2-\kappa_2}\sum_{j=0}^{\ell_2+\kappa_1}(-1)^{i+j}{-\ell_2-\kappa_2 \choose i}{\ell_2+\kappa_1 \choose j}\alpha^{\kappa_1-\kappa_2-i-j}\overline{\beta}^{i+j}\\
&\quad\quad\quad\quad\quad\quad\quad\quad\quad\times W_{(\underline{\ell};\,-\ell_2-\kappa_2-i,i,\ell_1+\kappa_2,j,\ell_2+\kappa_1-j,-\ell_3-\kappa_1)},\\
\rho({\rm diag}(-1,1)k)W_{(\underline{\kappa};\,\kappa_1-\kappa_2)} & = \sum_{m=0}^{\kappa_1-\kappa_2}(-1)^m{\kappa_1-\kappa_2 \choose m}\overline{\alpha}^{\kappa_1-\kappa_2-m}\beta^m \cdot W_{(\underline{\kappa};\,\kappa_1-\kappa_2-m)}.
\end{align*}
Note that (cf.\,\cite[Lemma 6.6]{Ichino2005})
\[
\int_{{\rm SU}(2)}\alpha^{\kappa_1-\kappa_2-i-j}\overline{\beta}^{i+j}\overline{\alpha}^{\kappa_1-\kappa_2-m}\beta^m\,dk = \begin{cases}\displaystyle{
(\kappa_1-\kappa_2+1)^{-1}{\kappa_1-\kappa_2 \choose m}^{-1}} & \mbox{ if $i+i=m$},\\
0 & \mbox{ if $i+j \neq m$}.
\end{cases}
\]
Hence we have
\begin{align*}
&Z(s,W_{(\underline{\ell};\,\underline{\kappa}^\vee,0)},W_{(\underline{\kappa};\,\kappa_1-\kappa_2)})\\
& = \sum_{i=0}^{-\ell_2-\kappa_2}\sum_{j=0}^{\ell_2+\kappa_1}\sum_{m=0}^{\kappa_1-\kappa_2}(-1)^{i+j+m}(\sqrt{-1})^{-\kappa_1-\kappa_2-i-j+m}{-\ell_2-\kappa_2 \choose i}{\ell_2+\kappa_1 \choose j}{\kappa_1-\kappa_2 \choose m}\\
&\quad\times \int_{{\rm SU}(2)}\alpha^{\kappa_1-\kappa_2-i-j}\overline{\beta}^{i+j}\overline{\alpha}^{\kappa_1-\kappa_2-m}\beta^m\,dk\\
&\quad\times \int_{0}^\infty d^\times a_1 \int_0^\infty d^\times a_2\int_{L_1}\frac{ds_1}{2\pi\sqrt{-1}}\int_{L_2}\frac{ds_2}{2\pi\sqrt{-1}}\int_{L_3}\frac{ds_3}{2\pi\sqrt{-1}}\,a_1^{-s_1-s_3+2s}a_2^{-s_2+4s}\\
&\quad\quad\quad\quad\quad\quad\quad\quad\times\Gamma_\C(\tfrac{s_1+\ell_1+\kappa_2+i+j}{2})\Gamma_\C(\tfrac{s_1-\ell_2-\kappa_2-i+j}{2})\Gamma_\C(\tfrac{s_1-\ell_3-\kappa_2-i-j}{2})\\
&\quad\quad\quad\quad\quad\quad\quad\quad\times\Gamma_\C(\tfrac{s_2-\ell_2-\ell_3-\kappa_1-\kappa_2}{2})\Gamma_\C(\tfrac{s_2+\ell_1-\ell_3-\kappa_1+\kappa_2}{2})\Gamma_\C(\tfrac{s_2+\ell_1+\ell_2+\kappa_1+\kappa_2}{2})\\
&\quad\quad\quad\quad\quad\quad\quad\quad\times\Gamma_\C(\tfrac{s_3+m}{2})\Gamma_\C(\tfrac{s_3+\kappa_1-\kappa_2-m}{2})\Gamma_\C(\tfrac{s_1+s_2+\ell_1-\ell_2-\ell_3-\kappa_1-i+j}{2})^{-1}\\
& = (\kappa_1-\kappa_2+1)^{-1}(\sqrt{-1})^{-\kappa_1-\kappa_2}\Gamma_\C(2s+\tfrac{-\ell_2-\ell_3-\kappa_1-\kappa_2}{2})\Gamma_\C(2s+\tfrac{\ell_1-\ell_3-\kappa_1+\kappa_2}{2})\Gamma_\C(2s+\tfrac{\ell_1+\ell_2+\kappa_1+\kappa_2}{2})\\
&\times\sum_{i=0}^{-\ell_2-\kappa_2}\sum_{j=0}^{\ell_2+\kappa_1}{-\ell_2-\kappa_2 \choose i}{\ell_2+\kappa_1 \choose j}\int_{L_1}\frac{ds_1}{2\pi\sqrt{-1}}\,\Gamma_\C(\tfrac{s_1+\ell_1+\kappa_2+i+j}{2})\Gamma_\C(\tfrac{s_1-\ell_2-\kappa_2-i+j}{2})\Gamma_\C(\tfrac{s_1-\ell_3-\kappa_2-i-j}{2})\\
&\quad\quad\quad\quad\quad\quad\quad\quad\quad\quad\quad\quad\quad\quad\times\Gamma_\C(\tfrac{-s_1+2s+i+j}{2})\Gamma_\C(\tfrac{-s_1+2s+\kappa_1-\kappa_2-i-j}{2})\Gamma_\C(\tfrac{s_1+4s+\ell_1-\ell_2-\ell_3-\kappa_1-i+j}{2})^{-1}.
\end{align*}
Here the last equality follows from the Mellin inversion formula.
Now we make a change of variable from $s_1$ to $s_1+i-j$.
By \cite[Lemma 12.7-(i)]{Ichino2005}, we have
\begin{align*}
&\sum_{i=0}^{-\ell_2-\kappa_2}{-\ell_2-\kappa_2 \choose i}\Gamma_\C(\tfrac{s_1+\ell_1+\kappa_2}{2}+i)\Gamma_\C(\tfrac{-s_1+2s+\kappa_1-\kappa_2}{2}-i) \\
&= \Gamma_\C(\tfrac{s_1+\ell_1+\kappa_2}{2})\Gamma_\C(s+\tfrac{\ell_1+\kappa_1}{2})\Gamma_\C(\tfrac{-s_1+2s+2\ell_2+\kappa_1+\kappa_2}{2})\Gamma_\C(s+\tfrac{\ell_1+2\ell_2+\kappa_1+2\kappa_2}{2})^{-1},\\
&\sum_{j=0}^{\ell_2+\kappa_1}{\ell_2+\kappa_1 \choose j}\Gamma_\C(\tfrac{-s_1+2s}{2}+j)\Gamma_\C(\tfrac{s_1-\ell_3-\kappa_2}{2}-j)\\
&= \Gamma_\C(\tfrac{-s_1+2s}{2})\Gamma_\C(s+\tfrac{-\ell_3-\kappa_2}{2})\Gamma_\C(\tfrac{s_1-2\ell_2-\ell_3-2\kappa_1-\kappa_2}{2})\Gamma_\C(s+\tfrac{-2\ell_2-\ell_3-2\kappa_1-\kappa_2}{2})^{-1}.
\end{align*}
Therefore, by our choice of the path $L_1$, we have
\begin{align*}
&Z(s,W_{(\underline{\ell};\,\underline{\kappa}^\vee,0)},W_{(\underline{\kappa};\,\kappa_1-\kappa_2)})\\
& = (\kappa_1-\kappa_2+1)^{-1}(\sqrt{-1})^{-\kappa_1-\kappa_2}\cdot\Gamma_\C(2s+\tfrac{-\ell_2-\ell_3-\kappa_1-\kappa_2}{2})\Gamma_\C(2s+\tfrac{\ell_1-\ell_3-\kappa_1+\kappa_2}{2})\Gamma_\C(2s+\tfrac{\ell_1+\ell_2+\kappa_1+\kappa_2}{2})\\
&\times\Gamma_\C(s+\tfrac{\ell_1+\kappa_1}{2})\Gamma_\C(s+\tfrac{-\ell_3-\kappa_2}{2})\Gamma_\C(s+\tfrac{\ell_1+2\ell_2+\kappa_1+2\kappa_2}{2})^{-1}\Gamma_\C(s+\tfrac{-2\ell_2-\ell_3-2\kappa_1-\kappa_2}{2})^{-1}\\
&\times \int_{L_1}\frac{ds_1}{2\pi\sqrt{-1}}\,\frac{\Gamma_\C(\tfrac{s_1-\ell_2-\kappa_2}{2})\Gamma_\C(\tfrac{s_1+\ell_1+\kappa_2}{2})\Gamma_\C(\tfrac{s_1-2\ell_2-\ell_3-2\kappa_1-\kappa_2}{2})\Gamma_\C(\tfrac{-s_1+2s}{2})\Gamma_\C(\tfrac{-s_1+2s+2\ell_2+\kappa_1+\kappa_2}{2})}{\Gamma_\C(\tfrac{s_1+4s+\ell_1-\ell_2-\ell_3-\kappa_1}{2})}.
\end{align*}
By the second Barnes lemma, the above integration in $s_1$ is equal to
\begin{align*}
&4\cdot\Gamma_\C(2s+\tfrac{\ell_1-\ell_3-\kappa_1+\kappa_2}{2})^{-1}\Gamma_\C(2s+\tfrac{-\ell_2-\ell_3-\kappa_1-\kappa_2}{2})^{-1}\Gamma_\C(2s+\tfrac{\ell_1+\ell_2+\kappa_1+\kappa_2}{2})^{-1}\\
&\times\Gamma_\C(s+\tfrac{-\ell_2-\kappa_2}{2})\Gamma_\C(s+\tfrac{\ell_2+\kappa_1}{2})\Gamma_\C(s+\tfrac{\ell_1+\kappa_2}{2})\Gamma_\C(s+\tfrac{-\ell_3-\kappa_1}{2})\\
&\times\Gamma_\C(s+\tfrac{\ell_1+2\ell_2+\kappa_1+2\kappa_2}{2})\Gamma_\C(s+\tfrac{-2\ell_2-\ell_3-2\kappa_1-\kappa_2}{2}).
\end{align*}
We conclude that
\begin{align*}
&Z(s,W_{(\underline{\ell};\,\underline{\kappa}^\vee,0)},W_{(\underline{\kappa};\,\kappa_1-\kappa_2)})\\
& = 4(\kappa_1-\kappa_2+1)^{-1}(\sqrt{-1})^{-\kappa_1-\kappa_2}\\
&\times\Gamma_\C(s+\tfrac{\ell_1+\kappa_1}{2})\Gamma_\C(s+\tfrac{\ell_1+\kappa_2}{2})\Gamma_\C(s+\tfrac{\ell_2+\kappa_1}{2})\Gamma_\C(s+\tfrac{-\ell_2-\kappa_2}{2})\Gamma_\C(s+\tfrac{-\ell_3-\kappa_1}{2})\Gamma_\C(s+\tfrac{-\ell_3-\kappa_2}{2})\\
& = 4(\kappa_1-\kappa_2+1)^{-1}(\sqrt{-1})^{-\kappa_1-\kappa_2}\cdot L(s,\itSigma_v \times \itPi_v)\\
& = 4(\kappa_1-\kappa_2+1)^{-1}(\sqrt{-1})^{\kappa_1+\kappa_2}\cdot L(s,\itSigma_v \times \itPi_v).
\end{align*}
Replac $\itSigma_v$ and $\itPi_v$ by $\itSigma_v^c$ and $\itPi_v^c$, respectively, we have
\begin{align*}
Z(s,W_{(\underline{\ell}^\vee,{\sf w}(\itSigma);\,\underline{\kappa},0)},W_{(\underline{\kappa}^\vee,{\sf w}(\itPi);\,\kappa_1-\kappa_2)})
& = 4(\kappa_1-\kappa_2+1)^{-1}(\sqrt{-1})^{\kappa_1+\kappa_2}\cdot L(s,\itSigma_v^c \times \itPi_v^c)\\
& = 4(\kappa_1-\kappa_2+1)^{-1}(\sqrt{-1})^{\kappa_1+\kappa_2}\cdot L(s,\itSigma_v \times \itPi_v).
\end{align*}
Note that
\begin{align*}
W_{(\underline{\ell}^\vee,{\sf w}(\itSigma);\,\underline{\kappa},\kappa_1-\kappa_2)} = \rho\left(\bp 0 & 1 & 0 \\ -1 & 0 & 0 \\ 0&0&1 \ep\right)W_{(\underline{\ell}^\vee,{\sf w}(\itSigma);\,\underline{\kappa},0)},\quad 
W_{(\underline{\kappa}^\vee,{\sf w}(\itPi);\,0)} = \rho\left(\bp 0 & 1 \\ -1 & 0 \ep\right) W_{(\underline{\kappa}^\vee,{\sf w}(\itPi);\,\kappa_1-\kappa_2)}. 
\end{align*}
Hence
\begin{align*}
&Z(s,W_{(\underline{\ell}^\vee,{\sf w}(\itSigma);\,\underline{\kappa},\kappa_1-\kappa_2)},W_{(\underline{\kappa}^\vee,{\sf w}(\itPi);\,0)})= Z(s,W_{(\underline{\ell}^\vee,{\sf w}(\itSigma);\,\underline{\kappa},0)},W_{(\underline{\kappa}^\vee,{\sf w}(\itPi);\,\kappa_1-\kappa_2)}).
\end{align*}
This completes the proof.
\end{proof}

\begin{thm}\label{T:RS}
Assume $\ell_{1,v} > -\kappa_{2,v} > \ell_{2,v} > -\kappa_{1,v} > \ell_{3,v}$ for all $v \in S_\infty$.
Let $m+\tfrac{1}{2}$ be critical for $L(s,\itSigma \times \itPi)$.
We have
\begingroup
\begin{align*}
&\sigma\left(\frac{L^{(\infty)}(m+\tfrac{1}{2},\itSigma \times \itPi)}{(2\pi\sqrt{-1})^{6dm}\cdot G(\omega_{\E/\F})\cdot G(\omega_\itPi)\cdot\pi^{3d(1+{\sf w}(\itSigma)+{\sf w}(\itPi))+\sum_{v \in S_\infty}(\ell_{1,v}-\ell_{3,v}+(\kappa_{1,v}-\kappa_{2,v})/2)}\cdot  p^b(\itSigma)\cdot p^b(\itPi)} \right)\\
&= \frac{L^{(\infty)}(m+\tfrac{1}{2},{}^\sigma\!\itSigma \times {}^\sigma\!\itPi)}{(2\pi\sqrt{-1})^{6dm}\cdot G(\omega_{\E/\F})\cdot G({}^\sigma\!\omega_\itPi)\cdot\pi^{3d(1+{\sf w}(\itSigma)+{\sf w}(\itPi))+\sum_{v \in S_\infty}(\ell_{1,v}-\ell_{3,v}+(\kappa_{1,v}-\kappa_{2,v})/2)}\cdot p^b({}^\sigma\!\itSigma)\cdot p^b({}^\sigma\!\itPi)}
\end{align*}
\endgroup
for all $\sigma \in {\rm Aut}(\C)$. 
\end{thm}

\begin{proof}
By Lemma \ref{L:Gauss sum} and Theorem \ref{T:Raghuram}, it suffices to show that
\begin{align}
\<[\itSigma_v]_b , [\itPi_v]_b\>_{m,v} &\in \sqrt{-1}\cdot\pi^{-6m-3-3{\sf w}(\itSigma)-3{\sf w}(\itPi) -\ell_{1,v}+\ell_{3,v}-(\kappa_{1,v}-\kappa_{2,v})/2}\cdot {\Q},\label{E:RS proof 1}\\
\<[\itSigma_v^c]_b , [\itPi_v^c]_b\>_{m,v}^c &= \<[\itSigma_v]_b , [\itPi_v]_b\>_{m,v} \label{E:RS proof 2}
\end{align}
for all $v \in S_\infty$.
Fix $v \in S_\infty$. In the rest of the proof, we drop the subscript $v$ for brevity. 
It is clear that the local zeta integral defines an element in 
\[
{\rm Hom}_{\GL_2(\C)}(\itSigma_v \otimes \itPi_v, \,|\mbox{ }|_\C^{-s+1/2}).
\]
Hence 
\[
Z(m+\tfrac{1}{2},W_{(\underline{\ell},{\sf w}(\itSigma);\,\underline{\kappa}',i)},W_{(\underline{\kappa},{\sf w}(\itPi);\,j)}) = 0
\]
unless $\underline{\kappa}' = \underline{\kappa}^\vee$ and $i+j = \kappa_{1}-\kappa_{2}$. 
We deduce that
\begin{align*}\label{E:RS proof 3}
\begin{split}
\<[\itSigma_v]_b , [\itPi_v]_b\>_{m} &= C_{(\underline{\ell};\,\underline{\kappa}^\vee)}\sum_{i=0}^{\kappa_{1}-\kappa_{2}}{\kappa_{1}-\kappa_{2}\choose i}^2\cdot Z(m+\tfrac{1}{2},W_{(\underline{\ell},{\sf w}(\itSigma);\,\underline{\kappa}^\vee,i)},W_{(\underline{\kappa},{\sf w}(\itPi);\,\kappa_{1}-\kappa_{2}-i)})\\
&\quad\times \left(s\otimes \<\cdot,\iota_{m}(\cdot)\>_{\mu_{\iota},\C}\otimes \<\cdot,\overline{\iota}_{m}(\cdot)\>_{\mu_{\overline{\iota}},\C}\right)\left(\xi_{b}^{(3)}({\bf v}_{(\underline{\ell}^\vee;\,\underline{\kappa},\kappa_{1}-\kappa_{2}-i)}),\,\xi_{b}^{(2)}(x^{\kappa_{1}-\kappa_{2}-i}y^i)\right).
\end{split}
\end{align*}
By the $\GL_2(\C)$-equivariance of the local zeta integral and Lemma \ref{L:RS}, we have
\begin{align*}
&Z(m+\tfrac{1}{2},W_{(\underline{\ell},{\sf w}(\itSigma);\,\underline{\kappa}^\vee,i)},W_{(\underline{\kappa},{\sf w}(\itPi);\,\kappa_{1}-\kappa_{2}-i)})\\
& = (-1)^i{\kappa_{1}-\kappa_{2} \choose i}^{-1}\cdot Z(m+\tfrac{1}{2},W_{(\underline{\ell},{\sf w}(\itSigma);\,\underline{\kappa}^\vee,0)},W_{(\underline{\kappa},{\sf w}(\itPi);\,\kappa_{1}-\kappa_{2})})\\
& \in L(m+\tfrac{1}{2}, \itSigma_v \times \itPi_v) \cdot \Q^\times = \pi^{-6m-3-3{\sf w}(\itSigma)-3{\sf w}(\itPi)-\ell_{1}+\ell_{3}-(\kappa_{1}-\kappa_{2})/2}\cdot \Q^\times.
\end{align*}
On the other hand, by (\ref{E:Lie algebra action 1}) and (\ref{E:Lie algebra action 2}), we have
\begin{align*}
\xi_{b}^{(2)}(x^jy^{\kappa_{1}-\kappa_{2}-j})\in  (\frak{g}_{2,\C} / \frak{k}_{2,\C})^*_\Q \,\otimes M_{\lambda_{\iota}}^\vee\otimes M_{\lambda_{\overline{\iota}}}^\vee,\quad 
\xi_{b}^{(3)}({\bf v}_{(\underline{\ell}^\vee;\,(\underline{\kappa}')^\vee,i)}) \in \wedge^3(\frak{g}_{3,\C} / \frak{k}_{3,\C})^*_\Q\otimes M_{\mu_{\iota}}^\vee\otimes M_{\mu_{\overline{\iota}}}^\vee
\end{align*}
for all $\underline{\ell} \succ \underline{\kappa}'$, $0 \leq i \leq \kappa_{1}'-\kappa_{2}'$, and $0 \leq j \leq \kappa_{1}-\kappa_{2}$.
It is clear that
\[
\<{\bf v}_{\iota},\iota_{m}({\bf u}_{\iota})\>_{\mu_{\iota},\C}\cdot \<{\bf v}_{\overline{\iota}},\overline{\iota}_{m}({\bf u}_{\overline{\iota}})\>_{\mu_{\overline{\iota}},\C} \in \Q
\]
for all ${\bf v}_{\iota} \otimes {\bf v}_{\overline{\iota}} \in M_{\mu_{\iota}}^\vee\otimes M_{\mu_{\overline{\iota}}}^\vee$ and ${\bf u}_{\iota} \otimes {\bf u}_{\overline{\iota}} \in M_{\lambda_{\iota}}^\vee\otimes M_{\lambda_{\overline{\iota}}}^\vee$.
To prove (\ref{E:RS proof 1}) holds, it remains to show that
\[
s\left(\wedge^3(\frak{g}_{3,\C} / \frak{k}_{3,\C})^*_\Q, \,(\frak{g}_{2,\C} / \frak{k}_{2,\C})^*_\Q\right) \subset \sqrt{-1}\cdot\Q.
\]
Indeed, the image is contained in a scalar multiple of $\Q$. To determine the scalar, note that the map $\frak{g}_{2,\C} / \frak{u}(2)_\C \rightarrow \frak{g}_{3,\C} / \frak{k}_{3,\C}$ induced by $\delta$ is given by
\[
e_{11}+e_{22} \longmapsto \frac{1}{3}\cdot Z_{12}+\frac{2}{3}\cdot Z_{23},\quad Y_{(0,0)} \longmapsto Z_{12},\quad Y_{\pm(1,-1)} \longmapsto X_{\pm(1,-1,0)}.
\]
Hence we have
\begin{align*}
\delta^*Z_{12}^* &= \frac{1}{3}\cdot (e_{11}^*+e_{22}^*) + {\rm pr}(Y_{(0,0)}^*),\quad \delta^*Z_{23}^* = \frac{2}{3}\cdot (e_{11}^*+e_{22}^*),\\ \delta^*X_{\pm(1,-1,0)}^* &= {\rm pr}(Y_{\pm(1,-1)}^*),\quad
\delta^*X_{\pm(0,1,-1)} = \delta^*X_{\pm(1,0,-1)} = 0.
\end{align*}
In particular, we have
\begin{align*}
&\delta^*Z_{23}^*\wedge \delta^*X_{(1,-1,0)}^* \wedge \delta^*X_{(-1,1,0)}^* \wedge {\rm pr}(Y_{(0,0)}^*)\\
 &= \frac{2}{3}\cdot (e_{11}^*+e_{22}^*)\wedge (-e_{12}^*-e_{21}^*)\wedge(\sqrt{-1}\,e_{12}^*-\sqrt{-1}\,e_{21}^*) \wedge (e_{11}^*-e_{22}^*) \otimes \sqrt{-1}\\
& = \frac{16}{3}\cdot e_{11}^*\wedge e_{22}^* \wedge e_{12}^* \wedge \sqrt{-1}\,e_{12}^* \otimes \sqrt{-1}.
\end{align*}
In other words, 
\[
s(Z_{23}\wedge X_{(1,-1,0)}^* \wedge X_{(-1,1,0)}^*,\,Y_{(0,0)}^*) = \frac{16}{3}\cdot \sqrt{-1}.
\]
Now we show that (\ref{E:RS proof 2}) holds.
By the $\GL_2(\C)$-equivariance of the local zeta integral and Lemma \ref{L:RS}, we have
\begin{align*}
& Z(m+\tfrac{1}{2},W_{(\underline{\ell},{\sf w}(\itSigma);\,\underline{\kappa}^\vee,i)},W_{(\underline{\kappa},{\sf w}(\itPi);\,\kappa_{1}-\kappa_{2}-i)}) = Z(m+\tfrac{1}{2},W_{(\underline{\ell}^\vee,{\sf w}(\itSigma);\,\underline{\kappa},\kappa_{1}-\kappa_{2}-i)},W_{(\underline{\kappa}^\vee,{\sf w}(\itPi);\,i)})
\end{align*}
for $0 \leq i \leq \kappa_{1}-\kappa_{2}$.
Therefore, 
\begin{align*}
&\<[\itSigma_v^c]_b , [\itPi_v^c]_b\>_{m}^c \\
&= C_{(\underline{\ell}^\vee;\,\underline{\kappa})}\sum_{i=0}^{\kappa_{1}-\kappa_{2}}{\kappa_{1}-\kappa_{2}\choose i}^2\cdot Z(m+\tfrac{1}{2},W_{(\underline{\ell}^\vee,{\sf w}(\itSigma);\,\underline{\kappa},i)},W_{(\underline{\kappa}^\vee,{\sf w}(\itPi);\,\kappa_{1}-\kappa_{2}-i)})\\
&\quad\times (-1)^{1+{\sf w}(\itPi)+{\sf w}(\itSigma)}\cdot(s\otimes \<\cdot,\overline{\iota}_{m}(\cdot)\>_{\mu_{\overline{\iota}},\C}\otimes \<\cdot,\iota_{m}(\cdot)\>_{\mu_{\iota},\C})\left(\xi_{b}^{c,(3)}({\bf v}_{(\underline{\ell};\,\underline{\kappa}^\vee,\kappa_{1}-\kappa_{2}-i)}),\,\xi_{b}^{c,(2)}(x^{\kappa_{1}-\kappa_{2}-i}y^i)\right)\\
&= C_{(\underline{\ell}^\vee;\,\underline{\kappa})}\sum_{i=0}^{\kappa_{1}-\kappa_{2}}{\kappa_{1}-\kappa_{2}\choose i}^2\cdot Z(m+\tfrac{1}{2},W_{(\underline{\ell},{\sf w}(\itSigma);\,\underline{\kappa}^\vee,\kappa_1-\kappa_2-i)},W_{(\underline{\kappa},{\sf w}(\itPi);\,i)})\\
&\quad\times (-1)^{1+{\sf w}(\itPi)+{\sf w}(\itSigma)}\cdot(s\otimes \<\cdot,\overline{\iota}_{m}(\cdot)\>_{\mu_{\overline{\iota}},\C}\otimes \<\cdot,\iota_{m}(\cdot)\>_{\mu_{\iota},\C})\left(\xi_{b}^{c,(3)}({\bf v}_{(\underline{\ell};\,\underline{\kappa}^\vee,\kappa_{1}-\kappa_{2}-i)}),\,\xi_{b}^{c,(2)}(x^{\kappa_{1}-\kappa_{2}-i}y^i)\right).
\end{align*}
Note that $C_{(\underline{\ell};\,\underline{\kappa}^\vee)}=C_{(\underline{\ell}^\vee;\,\underline{\kappa})}$ by (\ref{E:scalar relation}).
Therefore, to prove (\ref{E:RS proof 2}), it suffices to show that
\begin{align*}
&(s\otimes \<\cdot,\iota_{m}(\cdot)\>_{\mu_{\iota},\C}\otimes \<\cdot,\overline{\iota}_{m}(\cdot)\>_{\mu_{\overline{\iota}},\C})\left(\xi_{b}^{(3)}({\bf v}_{(\underline{\ell}^\vee;\,\underline{\kappa},i)}),\,\xi_{b}^{(2)}(x^iy^{\kappa_{1}-\kappa_{2}-i})\right)\\
&=(-1)^{1+{\sf w}(\itPi)+{\sf w}(\itSigma)}\cdot(s\otimes \<\cdot,\overline{\iota}_{m}(\cdot)\>_{\mu_{\overline{\iota}},\C}\otimes\<\cdot,\iota_{m}(\cdot)\>_{\mu_{\iota},\C})\left(\xi_{b}^{c,(3)}({\bf v}_{(\underline{\ell};\,\underline{\kappa}^\vee,\kappa_{1}-\kappa_{2}-i)}),\,\xi_{b}^{c,(2)}(x^{\kappa_{1}-\kappa_{2}-i}y^i)\right)
\end{align*}
for $0 \leq i \leq \kappa_1-\kappa_2$.
By the ${\rm U}(2)$-equivariance of the bilinear pairing, we need only to prove the euqlity for $i=0$.
Note that $c_{\underline{\kappa}}^{(2)}(x^{\kappa_{1}-\kappa_{2}}) = y^{\kappa_{1}-\kappa_{2}}$ and
\begin{align*}
c_{\underline{\ell}}^{(3)}({\bf v}_{(\underline{\ell};\,\underline{\kappa}^\vee,\kappa_{1}-\kappa_{2})})
& = (-1)^{\ell_{2}-\kappa_{2}}\cdot {\bf v}_{(\underline{\ell}^\vee;\,\underline{\kappa},0)} = (-1)^{1+{\sf w}(\itPi)+{\sf w}(\itSigma) }\cdot {\bf v}_{(\underline{\ell}^\vee;\,\underline{\kappa},0)}
\end{align*}
by (\ref{E:conjugate equiv. 3}).
It follows from definition that
\[
\left(\<\cdot,\iota_{m}(\cdot)\>_{\mu_{\iota},\C}\otimes \<\cdot,\overline{\iota}_{m}(\cdot)\>_{\mu_{\overline{\iota}},\C} \right)({\bf v},{\bf u}) = \left(\<\cdot,\overline{\iota}_{m}(\cdot)\>_{\mu_{\overline{\iota}},\C} \otimes \<\cdot,\iota_{m}(\cdot)\>_{\mu_{\iota},\C} \right)({\bf c}_{\mu^\vee}^{(3)}({\bf v}),{\bf c}_{\lambda^\vee}^{(2)}({\bf u}))
\]
for ${\bf v} \in M_{\mu,\C}^\vee$ and ${\bf u} \in M_{\lambda,\C}^\vee$. Therefore, it remains to show that
\[
s \circ \left(\wedge^3(c_{(1,0,-1)}^{(3)})^* , \, (c_{(1,-1)}^{(2)})^* \right) = s.
\]
Indeed, the two bilinear pairings on $\wedge^3 (\frak{g}_{3,\C} / \frak{k}_{3,\C})^* \times (\frak{g}_{2,\C} / \frak{k}_{2,\C})^*$ differ by a scalar. By (\ref{E:conjugate equiv. Lie 1}), (\ref{E:conjugate equiv. Lie 2}), and (\ref{E:conjugate equiv. Lie 3}) we have
\begin{align*}
&\delta^*\circ (c_{(1,0,-1)}^{(3)})^* (Z_{23}^*)\wedge \delta^*\circ (c_{(1,0,-1)}^{(3)})^* (X_{(1,-1,0)}^*) \wedge \delta^*\circ (c_{(1,0,-1)}^{(3)})^* (X_{(-1,1,0)}^*) \wedge {\rm pr}\circ (c_{(1,-1)}^{(2)})^* (Y_{(0,0)}^*)\\
&=\delta^*Z_{23}^* \wedge \delta^* (-X_{(-1,1,0)}^*) \wedge \delta^*(-X_{(1,-1,0)}^*) \wedge {\rm pr}(-Y_{(0,0)}^*),\\
& = \delta^*Z_{23}^* \wedge \delta^* X_{(1,-1,0)}^*\wedge \delta^* X_{(-1,1,0)}^* \wedge {\rm pr}(Y_{(0,0)}^*).
\end{align*}
Thus the scalar between the pairings is equal to $1$.
This completes the proof.
\end{proof}

\begin{rmk}
Under the normalization of the generators in \cite[\S\,1.5.3]{GL2020}, similar result was proved by Grobner and Lin \cite[Theorem A]{GL2020} for $\GL_n(\A_\E) \times \GL_{n-1}(\A_\E)$ with additional assumption that $\itPi_\infty^\vee = \itPi_\infty^c$, $\itSigma_\infty^\vee = \itSigma_\infty^c$, and the Galois-equivariance was proved for $\sigma \in {\rm Aut}(\C/\E^{\Gal})$.
For the central critical point, they also need to assume certain non-vanishing hypotheses on central critical values.
See also the result of Januszewski \cite[Theorem A]{Januszewski2019}.
\end{rmk}

As a consequence of Theorem \ref{T:RS}, we obtain the algebraicity of the ratios of critical values.
\begin{corollary}\label{C:ratio}
Assume $\ell_{1,v} > -\kappa_{2,v} > \ell_{2,v} > -\kappa_{1,v} > \ell_{3,v}$ for all $v \in S_\infty$.
Let $m_1+\tfrac{1}{2}$ and $m_2+\tfrac{1}{2}$ be critical for $L(s,\itSigma \times \itPi)$ with $m_2 \neq 0$.
We have
\begin{align*}
\sigma\left((2\pi\sqrt{-1})^{6d(m_2-m_1)}\cdot\frac{L^{(\infty)}(m_1+\tfrac{1}{2},\itSigma \times \itPi)}{L^{(\infty)}(m_2+\tfrac{1}{2},\itSigma \times \itPi)} \right) = (2\pi\sqrt{-1})^{6d(m_2-m_1)}\cdot\frac{L^{(\infty)}(m_1+\tfrac{1}{2},{}^\sigma\!\itSigma \times {}^\sigma\!\itPi)}{L^{(\infty)}(m_2+\tfrac{1}{2},{}^\sigma\!\itSigma \times {}^\sigma\!\itPi)}
\end{align*}
for all $\sigma \in {\rm Aut}(\C)$.
\end{corollary}

\begin{rmk}
We refer to \cite{Raghuram2020} for more general result over totally imaginary number fields.
\end{rmk}

\subsection{Period relations for base change lifting}

In this section, we establish period relations for Whittaker periods of $\GL_2$ and $\GL_3$ under base change. 
Let $\itPi$ and $\itSigma$ be cohomological irreducible cuspidal automorphic representations of $\GL_2(\A_\F)$ and $\GL_3(\A_\F)$ with central characters $\omega_\itSigma$ and $\omega_\itSigma$, respectively. We have $|\omega_\itPi| = |\mbox{ }|_{\A_\F}^{{\sf w}(\itPi)}$ and $|\omega_\itSigma| = |\mbox{ }|_{\A_\F}^{3{\sf w}(\itSigma)/2}$ for some ${\sf w}(\itPi) \in \Z$ and some even integer ${\sf w}(\itSigma)$. For each $v \in S_\infty$, we have
\[
\itPi_v = D_{\kappa_v} \otimes |\mbox{ }|_\R^{{\sf w}(\itPi)/2},\quad \itSigma_v = {\rm Ind}^{\GL_3(\R)}_{P_{2,1}(\R)}(D_{\ell_v} \boxtimes {\rm sgn}^\delta) \otimes |\mbox{ }|_\R^{{\sf w}(\itSigma)/2}
\]
for some $\kappa_v \in \Z_{\geq 2}$ such that $\kappa_v \equiv {\sf w}(\itPi) \,({\rm mod}\,2)$, some odd integer $\ell_v \in \Z_{\geq 3}$, and some $\delta \in \{0,1\}$ independent of $v$.
Here $P_{2,1}$ is the standard parabolic subgroup of $\GL_3$ with Levi part $\GL_2 \times \GL_1$.
Let $p^b(\itSigma)$ and $p(\itPi,\underline{\varepsilon})$ be the bottom degree Whittaker periods of $\itSigma$ and $\itPi$ defined as in \cite[Lemma 2.4 and \S\,4.1]{Chen2020}.

\begin{thm}\label{T:period relation}
Let ${\rm BC}_\E(\itPi)$ and ${\rm BC}_\E(\itSigma)$ be the base change lifts of $\itPi$ and $\itSigma$ to $\GL_2(\A_\E)$ and $\GL_3(\A_\E)$, respectively. 
\begin{itemize}
\item[(1)] Assume ${\rm BC}_\E(\itPi)$ is cuspidal. We have
\begin{align*}
\sigma \left( \frac{G(\omega_{\E/\F})\cdot p(\itPi,\underline{\varepsilon})\cdot p(\itPi,-\underline{\varepsilon})}{p^b({\rm BC}_\E(\itPi))}\right) = \frac{G(\omega_{\E/\F})\cdot p({}^\sigma\!\itPi,\underline{\varepsilon})\cdot p({}^\sigma\!\itPi,-\underline{\varepsilon})}{p^b({\rm BC}_\E({}^\sigma\!\itPi))} 
\end{align*}
for all $\sigma \in {\rm Aut}(\C)$ and $\underline{\varepsilon} \in \{\pm1\}^{S_\infty}$.
\item[(2)]
Assume ${\rm BC}_\E(\itSigma)$ is cuspidal and $\ell_v \geq 5$ for all $v \in S_\infty$. We have
\begin{align*}
\sigma \left( \frac{G(\omega_{\E/\F})\cdot p^b(\itSigma)^2}{p^b({\rm BC}_\E(\itSigma))} \right) = \frac{G(\omega_{\E/\F})\cdot p^b({}^\sigma\!\itSigma)^2}{p^b({\rm BC}_\E({}^\sigma\!\itSigma))}
\end{align*}
for all $\sigma \in {\rm Aut}(\C)$.
\end{itemize}
\end{thm}

\begin{proof}
Note that ${\rm BC}_\E(\itPi)$ and ${\rm BC}_\E(\itSigma)$ are cohomological. Indeed, we have
\[
{\rm BC}_\E(\itPi)_v = {\rm Ind}_{B_2(\C)}^{\GL_2(\C)}(\chi_{\kappa_{v}-1}\boxtimes \chi_{1-\kappa_v})\otimes |\mbox{ }|_\C^{{\sf w}(\itPi)/2},\quad {\rm BC}_\E(\itSigma)_v = {\rm Ind}_{B_3(\C)}^{\GL_3(\C)}(\chi_{\ell_v-1} \boxtimes {\bf 1} \boxtimes \chi_{1-\ell_v})\otimes |\mbox{ }|_\C^{{\sf w}(\itSigma)/2}
\]
for $v \in S_\infty$. Moreover, we have the factorization of $L$-functions:
\begin{align}
L(s,{\rm BC}_{\E}(\itPi),{\rm As}^+\otimes \omega_\itPi^{-1}) &= L(s,\itPi,{\rm Sym}^2\otimes \omega_\itPi^{-1})\cdot L(s,\omega_{\E/\F}),\label{E:period proof 1}\\
L(s,{\rm BC}_{\E}(\itSigma)\times{\rm BC}_{\E}(\itPi)) &= L(s,\itSigma \times \itPi)\cdot L(s,\itSigma \times \itPi\otimes\omega_{\E/\F})\label{E:period proof 2}.
\end{align}
By (\ref{E:period proof 1}), the assertion for $\itPi$ is a consequence of Theorem \ref{T:algebraicity results}-(2)-(ii), (3), (4), and Theorem \ref{T:Asai 2}-(1) below on the algebraicity of twisted Asai $L$-functinos for $\GL_2$.
To prove the assertion for $\itSigma$, we choose $\itPi$ such that 
\begin{align}\label{E:period proof 3}
\ell_v -1 > \kappa_v >2
\end{align}
for all $v \in S_\infty$.
The existence of such $\itPi$ is guaranteed by the assumption on $\ell_v$.
By (\ref{E:period proof 2}) and our result \cite[Theorem 4.11]{Chen2020} on the algebraicity of Rankin--Selberg $L$-functions for $\GL_3(\A_\F)\times\GL_2(\A_\F)$, we have
\begingroup
\begin{align}\label{E:period proof 4}
\begin{split}
&\sigma\left( \frac{L^{(\infty)}(m+\tfrac{1}{2},{\rm BC}_{\E}(\itSigma)\times{\rm BC}_{\E}(\itPi))}{(2\pi\sqrt{-1})^{6dm+\sum_{v \in S_\infty}(2\ell_v+\kappa_v+3{\sf w}(\itPi)+3{\sf w}(\itSigma))}\cdot G(\omega_{\E/\F}\cdot\omega_\itPi^2)\cdot p^b(\itSigma)^2\cdot p(\itPi,\varepsilon_m(\itSigma))\cdot p(\itPi,-\varepsilon_m(\itSigma))}\right) \\
& = \frac{L^{(\infty)}(m+\tfrac{1}{2},{\rm BC}_{\E}({}^\sigma\!\itSigma)\times{\rm BC}_{\E}({}^\sigma\!\itPi))}{(2\pi\sqrt{-1})^{6dm+\sum_{v \in S_\infty}(2\ell_v+\kappa_v+3{\sf w}(\itPi)+3{\sf w}(\itSigma))}\cdot G(\omega_{\E/\F}\cdot{}^\sigma\!\omega_\itPi^2)\cdot p^b({}^\sigma\!\itSigma)^2\cdot p({}^\sigma\!\itPi,\varepsilon_m(\itSigma))\cdot p({}^\sigma\!\itPi,-\varepsilon_m(\itSigma))}
\end{split}
\end{align}
\endgroup
for all $\sigma \in {\rm Aut}(\C)$ and all critical points $m+\tfrac{1}{2}$ of $L(s,\itSigma \times \itPi)$, where $\varepsilon_m(\itSigma) = (-1)^{m+{\sf w}(\itSigma)/2+\delta}$.
Note that condition (\ref{E:period proof 3}) implies that $L(s,\itSigma \times \itPi)$ has at least two critical points. In particular, there is a non-zero critical value for $L(s,{\rm BC}_{\E}(\itSigma)\times{\rm BC}_{\E}(\itPi))$.
The period relation then follows from comparing Theorem \ref{T:RS} with (\ref{E:period proof 4}).
This complete the proof.
\end{proof}

\section{Algebraicity of the adjoint $L$-functions for $\GL_2$ and $\GL_3$}\label{S:adjoint}

Let $\E$ be a totally imaginary quadratic extension of $\F$.

\subsection{Algebraicity of the adjoint $L$-functions for $\GL_n$}
Let $\itPi$ be a cohomological irreducible cuspidal automorphic representation of $\GL_n(\A_\E)$. 
Let 
\[
\mu = \prod_{v \in S_\infty} \mu_v \in \prod_{v \in S_\infty} (X^+(T_n) \times X^+(T_n))
\]
with $\mu_v = (\mu_{\iota_v},\mu_{\overline{\iota}_v})$ such that $\itPi$ contributes to the cuspidal cohomology of $\GL_n(\A_\E)$ with coefficients in $M^\vee_{\mu,\C}$.
We keep the notation of \S\,\ref{SS:Whittaker GL_n}.
With respect to the choice of the generators (\ref{E:generators}) in the cuspidal cohomology, let $p^b({}^\sigma\!\itPi)$ and $p^t({}^\sigma\!\itPi)$ be the Whittaker periods of ${}^\sigma\!\itPi$ for each $\sigma\in{\rm Aut}(\C)$ defined as in Lemma \ref{L:Whittaker periods}.
Note that the periods are normalized so that the diagram in Lemma \ref{L:Whittaker periods} commutes.

Let 
\begin{align*}
s^{(n)}: \wedge^{b_n} (\frak{g}_{n,\C} / \frak{k}_{n,\C})^* \times \wedge^{t_n} (\frak{g}_{n,\C} / \frak{k}_{n,\C})^* &\longrightarrow \C
\end{align*}
be the ${\rm U}(n)$-equivariant bilinear pairing defined by
\begin{align*}
(\wedge_{i=1}^{b_n}X_i) \wedge (\wedge_{j=1}^{t_n}Y_j) &= s^{(n)}\left(\wedge_{i=1}^{b_n}X_i, \wedge_{j=1}^{t_n}Y_j\right)\\
&\times(\wedge_{i=1}^{n-1} \,e_{ii}^*) \wedge (\wedge_{1\leq i <j \leq n}\,e_{ij}^* )\wedge (\wedge_{1\leq i <j \leq n}\,\sqrt{-1}\,e_{ij}^* ).
\end{align*}
Here $\{e_{kk},e_{i,j},\sqrt{-1}\,e_{i,j}\,\vert\, 1\leq k \leq n-1,\,1\leq i<j \leq n \} \subset \frak{g}_{n,\C}/\frak{k}_{n,\C}$ is ordered lexicographically.
Note that $(\wedge_{i=1}^{n-1} \,e_{ii}^*) \wedge (\wedge_{1\leq i <j \leq n}\,e_{ij}^* )\wedge (\wedge_{1\leq i <j \leq n}\,\sqrt{-1}\,e_{ij}^* )$ corresponds to the standard measure on $\GL_n(\C)/K_n$ defined in \S\,\ref{SS:measure}.
Let $v \in S_\infty$.
Let 
\begin{align*}
\<\cdot,\cdot\>_v^{(n)} : \mathcal{W}(\itPi_v,\psi_{n,\C}) \times \mathcal{W}(\itPi_v^\vee,\psi_{n,\C}) &\longrightarrow \C
\end{align*}
be the non-zero $\GL_n(\C)$-equivariant bilinear pairing defined by the local integrals
\begin{align}\label{E:local bilinear}
\begin{split}
\<W_1,W_2\>_v^{(n)} 
&= \int_{N_{n-1}(\C)\backslash \GL_{n-1}(\C)} W_1\left(\bp g & 0 \\ 0 & 1\ep\right)W_2\left({\rm diag}((-1)^{n-1},(-1)^{n-2},\cdots,1)\bp g & 0 \\ 0 & 1\ep\right)\,dg.
\end{split}
\end{align}
Note that the integrals converge absolutely (cf.\,\cite[(3.17)]{JS1981}).
We define
\begin{align*}
\<\cdot,\cdot\>_v^{c,(n)} : \mathcal{W}(\itPi_v^c,\psi_{n,\C}) \times \mathcal{W}((\itPi_v^c)^\vee,\psi_{n,\C}) &\longrightarrow \C
\end{align*}
in a similar way.
We also fix non-zero $\GL_n(\C)$-equivariant bilinear pairings
\begin{align*}
&\<\cdot,\cdot\>_{\mu_{\iota_v},\C} : M_{\mu_{\iota_v},\C}^\vee \times M_{\mu_{\iota_v},\C} \longrightarrow \C,\quad \<\cdot,\cdot\>_{\mu_{\overline{\iota}_v,\C}} : M_{\mu_{\overline{\iota}_v,\C}}^\vee \times M_{\mu_{\overline{\iota}_v,\C}} \longrightarrow \C
\end{align*}
defined over $\Q$.
Let
\begin{align*}
B(\cdot,\cdot)_v^{(n)}:
&\left(\mathcal{W}(\itPi_v,\psi_{n,\C}) \otimes \wedge^{b_n}(\frak{g}_{n,\C} / \frak{k}_{n,\C})^* \otimes M_{\lambda_v,\C}^\vee\right) \times \left(\mathcal{W}(\itPi_v^\vee,\psi_{n,\C}) \otimes \wedge^{t_n} (\frak{g}_{n,\C} / \frak{k}_{n,\C})^* \otimes M_{\lambda_v,\C}\right)\longrightarrow \C
\end{align*}
be the bilinear pairing defined by
\begin{align*}
&B\left(W_1 \otimes \wedge_{i=1}^{b_n}X_i^*\otimes ({\bf v}_{\iota_v}'\otimes{\bf v}_{\overline{\iota}_v}') ,\, W_2 \otimes \wedge_{j=1}^{t_n}Y_j^* \otimes ({\bf u}_{\iota_v}'\otimes{\bf u}_{\overline{\iota}_v}')\right)_v^{(n)}\\
&= \<W_1,W_2\>_v^{(n)}\cdot s^{(n)}\left(\wedge_{i=1}^{b_n}X_i^*,\wedge_{j=1}^{t_n}Y_j^*\right)\cdot \<{\bf v}_{{\iota}_v}',{\bf u}_{\iota_v}'\>_{\mu_{\iota_v},\C}\cdot \<{\bf v}_{\overline{\iota}_v}',{\bf u}_{\overline{\iota}_v}'\>_{\mu_{\overline{\iota}_v},\C}.
\end{align*}
We define bilinear pairing $B(\cdot,\cdot)_v^{c,(n)}$ in a similar way by replacing $\itPi_v$ and $\mu_v$ by $\itPi_v^c$ and $\mu_v^c$, respectively.
For $\sigma \in {\rm Aut}(\C)$, we then have the bilinear pairing
\begin{align*}
B(\cdot,\cdot)_{\sigma,\infty}^{(n)} &= \left( \bigotimes_{v \in S_\infty,\,\sigma^{-1}\circ\iota_v = \iota_{\sigma^{-1}\circ v}}B(\cdot,\cdot)_{\sigma^{-1}\circ v}^{(n)} \right) \otimes \left( \bigotimes_{v \in S_\infty,\,\sigma^{-1}\circ\iota_v = \overline{\iota}_{\sigma^{-1}\circ v}}B(\cdot,\cdot)_{\sigma^{-1}\circ v}^{c,(n)} \right)
\end{align*}
on
\begin{align*}
&\left(\mathcal{W}({}^\sigma\!\itPi_\infty,\psi_{n,\E_\infty}) \otimes \wedge^{b_nd}(\frak{g}_{n,\infty,\C} / \frak{k}_{n,\infty,\C})^* \otimes M_{{}^\sigma\!\mu,\C}^\vee\right)\times \left(\mathcal{W}({}^\sigma\!\itPi_\infty^\vee,\psi_{n,\E_\infty}) \otimes \wedge^{t_nd} (\frak{g}_{n,\infty,\C} / \frak{k}_{n,\infty,\C})^* \otimes M_{{}^\sigma\!\mu,\C}\right).
\end{align*}
We simply write $B(\cdot,\cdot)_{\sigma,\infty}^{(n)} = B(\cdot,\cdot)_{\infty}^{(n)}$ when $\sigma$ is the identity map.
Let 
\[
L(s,\itPi \times \itPi^\vee)
\]
be the Rankin--Selberg $L$-function of $\itPi \times \itPi^\vee$. The corresponding $L$-function without archimedean $L$-factors is denoted by $L^{(\infty)}(s,\itPi \times \itPi^\vee)$.
We have the following result of Grobner--Harris--Lapid and Balasubramanyam--Raghuram in \cite[Theorem 5.3]{GHL2016} and \cite[Theorem 3.3.11]{BR2017} on the algebraicity of 
\[
{\rm Res}_{s=1}L^{(\infty)}(s,\itPi \times \itPi^\vee)
\]
in terms of the product of the bottom degree and top degree Whittaker periods. 

\begin{thm}[Grobner--Harris--Lapid,\,Balasubramanyam--Raghuram]\label{T:BR}
\noindent
\begin{itemize}
\item[(1)] We have
\begin{align*}
&\sigma\left(\frac{{\rm Res}_{s=1}L^{(\infty)}(s,\itPi \times \itPi^\vee)}
{|D_\K|^{n(n+1)/4}\cdot{\rm Reg}_\E\cdot\pi^{nd}\cdot p^b(\itPi)\cdot p^t(\itPi^\vee)}\cdot B([\itPi_\infty]_b,[\itPi_\infty^\vee]_t)_\infty^{(n)}\right) \\
&=\frac{{\rm Res}_{s=1}L^{(\infty)}(s,{}^\sigma\!\itPi \times {}^\sigma\!\itPi^\vee)}
{|D_\K|^{n(n+1)/4}\cdot{\rm Reg}_\E\cdot\pi^{nd}\cdot p^b({}^\sigma\!\itPi)\cdot p^t({}^\sigma\!\itPi^\vee)}\cdot B([{}^\sigma\!\itPi_\infty]_b,[{}^\sigma\!\itPi_\infty^\vee]_t)_{,\sigma,\infty}^{(n)}
\end{align*}
for all $\sigma \in {\rm Aut}(\C)$.
\item[(2)] The archimedean factor $B([\itPi_\infty]_b,[\itPi_\infty^\vee]_t)_\infty^{(n)}$ is non-zero.
\end{itemize}
\end{thm}

\begin{rmk}
The factor $|D_\E|^{n(n+1)/4}$ is due to the comparison between Haar measures on $\A_\E^\times \backslash \GL_n(\A_\E)$ and $N_n(\A_\E)$ (cf.\,\cite[(5.2)]{GHL2016}). The factor ${\rm Reg}_\E$ appears in the volume of $\E^\times \backslash \A_\E^1$ with respect to the standard measure $\prod_v d^\times x_v$ on $\A_\E^\times$.
As for the factor $\pi^{nd}$, note that we have
\begin{align*}
\int_{\GL_n(\C)}f(g)\,dg
& = \Gamma_\C(n)\cdot\int_{\C^{n-1}}dx\int_{\C^{n-1}}dy\int_{\C^\times}\frac{d a}{|a|_\C} \int_{\GL_{n-1}(\C)}\frac{dg}{|\det(g)|_\C}\\
&\quad\quad\quad\quad\quad\quad\quad\quad f\left(\bp {\bf 1}_{n-1} & x \\ 0 & 1\ep \bp g & 0 \\ 0 & 1\ep \bp {\bf 1}_{n-1} & 0 \\ {}^ty & 1\ep a \right)
\end{align*}
for $f \in L^1(\GL_n(\C))$, where the measures are the standard measures defined in \S\,\ref{SS:measure}.
This formula is used to relate the Rankin--Selberg local zeta integrals for $\itPi_v \times \itPi_v^\vee$ at $s=1$ to the bilinear pairing (\ref{E:local bilinear}) (cf.\,\cite[(3.9)]{Zhang2014}) for $v \in S_\infty$.
\end{rmk}

\subsection{Algebraicity of the adjoint $L$-functions for $\GL_2$ and $\GL_3$}

In this section, we refine Theorem \ref{T:BR} for $n=2$ and $n=3$ in Theorem \ref{T:adjoint} below by explicitly determine the archimedean factors.
We keep the notation of \S\,\ref{S:RS}.

In the following lemmas, we compute the archimedean local pairings appearing in the archimedean factors $B([\itPi_\infty]_b,[\itPi_\infty^\vee]_t)_\infty^{(2)}$ and $B([\itSigma_\infty]_b,[\itSigma_\infty^\vee]_t)_\infty^{(3)}$.
\begin{lemma}\label{L:BR 1}
Let $v \in S_\infty$.
We have
\begin{align*}
&\< W_{(\underline{\kappa}_v,{\sf w}(\itPi);\,0)}, W_{(\underline{\kappa}_v^\vee,-{\sf w}(\itPi);\,\kappa_{1,v}-\kappa_{2,v})}\>_v^{(2)} \\
&= \< W_{(\underline{\kappa}_v^\vee,{\sf w}(\itPi);\,\kappa_{1,v}-\kappa_{2,v})}, W_{(\underline{\kappa}_v,-{\sf w}(\itPi);\,0)}\>_v^{c,(2)}\\
&= 2^2\cdot \frac{\Gamma_\C(1+\kappa_{1,v}-\kappa_{2,v})}{\Gamma_\C(1)\Gamma_\C(2+\kappa_{1,v}-\kappa_{2,v})}\cdot L(1,\itPi_v \times \itPi_v^\vee).
\end{align*}
\end{lemma}

\begin{proof}
We drop the subscript $v$ for brevity.
Note that $W_{(\underline{\kappa},{\sf w};\,i)} = W_{(\underline{\kappa},0;\,i)} \otimes |\mbox{ }|_\C^{{\sf w}/2}$ and $W_{(\underline{\kappa}^\vee,{\sf w};\,i)} = W_{(\underline{\kappa}^\vee,0;\,i)} \otimes |\mbox{ }|_\C^{{\sf w}/2}$ for ${\sf w} \in \Z$ and $0 \leq i \leq \kappa_1-\kappa_2$.
Hecne
\[
\< W_{(\underline{\kappa},{\sf w}(\itPi);\,0)}, W_{(\underline{\kappa}^\vee,-{\sf w}(\itPi);\,\kappa_1-\kappa_2)}\>^{(2)} = \< W_{(\underline{\kappa}^\vee,{\sf w}(\itPi);\,\kappa_1-\kappa_2)}, W_{(\underline{\kappa},-{\sf w}(\itPi);\,0)}\>^{c,(2)}.
\]
By the Mellin inversion formula and the first Barnes lemma, we have
\begin{align*}
&\< W_{(\underline{\kappa},{\sf w}(\itPi);\,0)}, W_{(\underline{\kappa}^\vee,-{\sf w}(\itPi);\,\kappa_1-\kappa_2)}\>^{(2)} \\
&= \int_{\C^\times} W_{(\underline{\kappa},{\sf w}(\itPi);\,0)}({\rm diag}(z,1))W_{(\underline{\kappa}^\vee,-{\sf w}(\itPi);\,\kappa_1-\kappa_2)}({\rm diag}(-z,1))\,d^\times z\\
&= \int_{0}^\infty
d^\times a \int_{L}\frac{ds}{2\pi\sqrt{-1}}\int_{L}\frac{ds'}{2\pi\sqrt{-1}}\,a^{-s-s'+2}\Gamma_\C(\tfrac{s}{2})\Gamma_\C(\tfrac{s+\kappa_1-\kappa_2}{2})\Gamma_\C(\tfrac{s'}{2})\Gamma_\C(\tfrac{s'+\kappa_1-\kappa_2}{2})\\
& = \int_{L}\frac{ds}{2\pi\sqrt{-1}}\,\Gamma_\C(\tfrac{s}{2})\Gamma_\C(\tfrac{s+\kappa_1-\kappa_2}{2})\Gamma_\C(1-\tfrac{s}{2})\Gamma_\C(1+\tfrac{-s+\kappa_1-\kappa_2}{2})\\
& = 4\cdot\frac{\Gamma_\C(1)\Gamma_\C(1+\tfrac{\kappa_1-\kappa_2}{2})^2\Gamma_\C(1+\kappa_1-\kappa_2)}{\Gamma_\C(2+\kappa_1-\kappa_2)}\\
& =  4\cdot\frac{\Gamma_\C(1+\kappa_1-\kappa_2)}{\Gamma_\C(1)\Gamma_\C(2+\kappa_1-\kappa_2)}\cdot L(1,\itPi_v \times \itPi_v^\vee).
\end{align*} 
This completes the proof.
\end{proof}

\begin{lemma}\label{L:BR 2}
Let $v \in S_\infty$.
We have
\begin{align*}
& \<W_{(\underline{\ell}_v,{\sf w}(\itSigma);\,(\ell_{2,v},\ell_{2,v}),0)},W_{(\underline{\ell}_v^\vee,-{\sf w}(\itSigma);\,(-\ell_{2,v},-\ell_{2,v}),0)}\>_v^{(3)}\\
& =  \<W_{(\underline{\ell}_v^\vee,{\sf w}(\itSigma);\,(-\ell_{2,v},-\ell_{2,v}),0)},W_{(\underline{\ell}_v,-{\sf w}(\itSigma);\,(\ell_{2,v},\ell_{2,v}),0)}\>_v^{c,(3)}\\
& = 2^4\cdot\frac{\Gamma_\C(1+\ell_{1,v}-\ell_{2,v})\Gamma_\C(1+\ell_{2,v}-\ell_{3,v})}{\Gamma_\C(1)^2\Gamma_\C(3+\ell_{1,v}-\ell_{3,v})}\cdot L(1,\itSigma_v \times \itSigma_v^\vee).
\end{align*}
\end{lemma}

\begin{proof}
We drop the subscript $v$ for brevity.
Note that
\[
W_{(\underline{\ell},{\sf w};\,n_1,n_2,n_3,n_{23},n_{13},n_{12})} = W_{(\underline{\ell},0;\,n_1,n_2,n_3,n_{23},n_{13},n_{12})} \otimes |\mbox{ }|_\C^{{\sf w}/2}
\]
for ${\sf w} \in \Z$ and $(n_1,n_2,n_3,n_{23},n_{13},n_{12}) \in \Z^6$ with $n_1+n_2+n_3 = \ell_1-\ell_2$ and $n_{23}+n_{13}+n_{12}=\ell_2-\ell_3$. Similarly for $\underline{\ell}^\vee$. 
Hence
\begin{align*}
& \<W_{(\underline{\ell},{\sf w}(\itSigma);\,(\ell_2,\ell_2),0)},W_{(\underline{\ell}^\vee,-{\sf w}(\itSigma);\,(-\ell_2,-\ell_2),0)}\>^{(3)} =  \<W_{(\underline{\ell}^\vee,{\sf w}(\itSigma);\,(-\ell_2,-\ell_2),0)},W_{(\underline{\ell},-{\sf w}(\itSigma);\,(\ell_2,\ell_2),0)}\>^{c,(3)}.
\end{align*}
Let $L_1$ and $L_2$ be vertical paths from $c_1-\sqrt{-1}\,\infty$ to $c_1+\sqrt{-1}\,\infty$ and $c_2-\sqrt{-1}\,\infty$ to $c_2+\sqrt{-1}\,\infty$, respectively, for some
$0<c_1<\min\{2+\ell_1-\ell_2,2+\ell_2-\ell_3\}$ and $0<c_2<\min\{4+\ell_1-\ell_2,4+\ell_2-\ell_3\}$.
We have
\begin{align*}
& \<W_{(\underline{\ell},{\sf w}(\itSigma);\,(\ell_2,\ell_2),0)},W_{(\underline{\ell}^\vee,-{\sf w}(\itSigma);\,(-\ell_2,-\ell_2),0)}\>^{(3)}\\
&= \int_{0}^\infty d^\times a_1 \int_0^\infty d^\times a_2 \int_{L_1}\frac{ds_1}{2\pi\sqrt{-1}}\int_{L_2}\frac{ds_2}{2\pi\sqrt{-1}}\int_{L_1}\frac{ds_1'}{2\pi\sqrt{-1}}\int_{L_2}\frac{ds_2'}{2\pi\sqrt{-1}}\,a_1^{-s_1-s_1'+2}a_2^{-s_2-s_2'+4}\\
&\quad\quad\quad\times \Gamma_\C(\tfrac{s_1+\ell_1-\ell_2}{2})\Gamma_\C(\tfrac{s_1}{2})\Gamma_\C(\tfrac{s_1+\ell_2-\ell_3}{2})\Gamma_\C(\tfrac{s_2+\ell_2-\ell_3}{2})\Gamma_\C(\tfrac{s_2+\ell_1-\ell_3}{2})\Gamma_\C(\tfrac{s_2+\ell_1-\ell_2}{2})\Gamma_\C(\tfrac{s_1+s_2+\ell_1-\ell_3}{2})^{-1}\\
&\quad\quad\quad\times \Gamma_\C(\tfrac{s_1'+\ell_1-\ell_2}{2})\Gamma_\C(\tfrac{s_1'}{2})\Gamma_\C(\tfrac{s_1'+\ell_2-\ell_3}{2})\Gamma_\C(\tfrac{s_2'+\ell_2-\ell_3}{2})\Gamma_\C(\tfrac{s_2'+\ell_1-\ell_3}{2})\Gamma_\C(\tfrac{s_2'+\ell_1-\ell_2}{2})\Gamma_\C(\tfrac{s_1'+s_2'+\ell_1-\ell_3}{2})^{-1}\\
&= \int_{L_2}\frac{ds_2}{2\pi\sqrt{-1}}\, \Gamma_\C(\tfrac{s_2+\ell_2-\ell_3}{2})\Gamma_\C(\tfrac{s_2+\ell_1-\ell_3}{2})\Gamma_\C(\tfrac{s_2+\ell_1-\ell_2}{2})\\
&\quad\times\Gamma_\C(\tfrac{-s_2+4+\ell_2-\ell_3}{2})\Gamma_\C(\tfrac{-s_2+4+\ell_1-\ell_3}{2})\Gamma_\C(\tfrac{-s_2+4+\ell_1-\ell_2}{2})\\
&\quad\times \int_{L_1}\frac{ds_1}{2\pi\sqrt{-1}}\,\Gamma_\C(\tfrac{s_1+\ell_1-\ell_2}{2})\Gamma_\C(\tfrac{s_1}{2})\Gamma_\C(\tfrac{s_1+\ell_2-\ell_3}{2})\Gamma_\C(\tfrac{s_1+s_2+\ell_1-\ell_3}{2})^{-1}\\
&\quad\times \Gamma_\C(\tfrac{-s_1+2+\ell_1-\ell_2}{2})\Gamma_\C(\tfrac{-s_1+2}{2})\Gamma_\C(\tfrac{-s_1+2+\ell_2-\ell_3}{2})\Gamma_\C(\tfrac{-s_1-s_2+6+\ell_1-\ell_3}{2})^{-1}.
\end{align*}
Here the last equality follows from Mellin inversion formula.
By the first Barnes lemma, we have
\begin{align*}
&\int_{L_1'}\frac{dt_1}{2\pi\sqrt{-1}}\,\Gamma_\C(t_1+\tfrac{s_1}{2})\Gamma_\C(t_1+\tfrac{s_2}{2})\Gamma_\C(-t_1+\tfrac{\ell_1-\ell_2}{2})\Gamma_\C(-t_1+\tfrac{\ell_2-\ell_3}{2})\\
& = 2\cdot \Gamma_\C(\tfrac{s_1+\ell_1-\ell_2}{2}) \Gamma_\C(\tfrac{s_1+\ell_2-\ell_3}{2})\Gamma_\C(\tfrac{s_2+\ell_1-\ell_2}{2}) \Gamma_\C(\tfrac{s_2+\ell_2-\ell_3}{2})\Gamma_\C(\tfrac{s_1+s_2+\ell_1-\ell_3}{2})^{-1},\\
&\int_{L_2'}\frac{dt_2}{2\pi\sqrt{-1}}\,\Gamma_\C(t_2-\tfrac{s_1}{2})\Gamma_\C(t_2+\tfrac{-s_2+2}{2})\Gamma_\C(-t_2+\tfrac{2+\ell_1-\ell_2}{2})\Gamma_\C(-t_2+\tfrac{2+\ell_2-\ell_3}{2})\\
& = 2 \cdot \Gamma_\C(\tfrac{-s_1+2+\ell_1-\ell_2}{2})\Gamma_\C(\tfrac{-s_1+2+\ell_2-\ell_3}{2})\Gamma_\C(\tfrac{-s_2+4+\ell_1-\ell_2}{2})\Gamma_\C(\tfrac{-s_2+4+\ell_2-\ell_3}{2})\Gamma_\C(\tfrac{-s_1-s_2+6+\ell_1-\ell_3}{2})^{-1}.
\end{align*}
Here $L_1'$ and $L_2'$ are vertical paths from $c_1'-\sqrt{-1}\,\infty$ to $c_1'+\sqrt{-1}\,\infty$ and $c_2'-\sqrt{-1}\,\infty$ to $c_2'+\sqrt{-1}\,\infty$, respectively, for some
$0<c_1'<\min\{\tfrac{\ell_1-\ell_2}{2},\tfrac{\ell_2-\ell_3}{2}\}$ and $\max\{\tfrac{c_1}{2},\tfrac{c_2-2}{2}\}< c_2' < \min\{\tfrac{2+\ell_1-\ell_2}{2},\tfrac{2+\ell_2-\ell_3}{2}\}$.
Therefore we have
\begin{align*}
&\int_{L_1}\frac{ds_1}{2\pi\sqrt{-1}}\,\Gamma_\C(\tfrac{s_1+\ell_1-\ell_2}{2})\Gamma_\C(\tfrac{s_1}{2})\Gamma_\C(\tfrac{s_1+\ell_2-\ell_3}{2})\Gamma_\C(\tfrac{s_1+s_2+\ell_1-\ell_3}{2})^{-1}\\
&\quad\times \Gamma_\C(\tfrac{-s_1+2+\ell_1-\ell_2}{2})\Gamma_\C(\tfrac{-s_1+2}{2})\Gamma_\C(\tfrac{-s_1+2+\ell_2-\ell_3}{2})\Gamma_\C(\tfrac{-s_1-s_2+6+\ell_1-\ell_3}{2})^{-1}\\
& = 2^{-2}\cdot \Gamma_\C(\tfrac{s_2+\ell_1-\ell_2}{2})^{-1} \Gamma_\C(\tfrac{s_2+\ell_2-\ell_3}{2})^{-1}\Gamma_\C(\tfrac{-s_2+4+\ell_1-\ell_2}{2})^{-1}\Gamma_\C(\tfrac{-s_2+4+\ell_2-\ell_3}{2})^{-1}\\
&\times \int_{L_1'}\frac{dt_1}{2\pi\sqrt{-1}}\,\Gamma_\C(t_1+\tfrac{s_2}{2})\Gamma_\C(-t_1+\tfrac{\ell_1-\ell_2}{2})\Gamma_\C(-t_1+\tfrac{\ell_2-\ell_3}{2})\\
&\times \int_{L_2'}\frac{dt_2}{2\pi\sqrt{-1}}\,\Gamma_\C(t_2+\tfrac{-s_2+2}{2})\Gamma_\C(-t_2+\tfrac{2+\ell_1-\ell_2}{2})\Gamma_\C(-t_2+\tfrac{2+\ell_2-\ell_3}{2})\\
&\times\int_{L_1}\frac{ds_1}{2\pi\sqrt{-1}} \Gamma_\C(\tfrac{s_1}{2})\Gamma_\C(\tfrac{s_1}{2}+t_1)\Gamma_\C(\tfrac{-s_1+2}{2})\Gamma_\C(\tfrac{-s_1}{2}+t_2)\\
& = \Gamma_\C(1)\Gamma_\C(\tfrac{s_2+\ell_1-\ell_2}{2})^{-1} \Gamma_\C(\tfrac{s_2+\ell_2-\ell_3}{2})^{-1}\Gamma_\C(\tfrac{-s_2+4+\ell_1-\ell_2}{2})^{-1}\Gamma_\C(\tfrac{-s_2+4+\ell_2-\ell_3}{2})^{-1}\\
&\times \int_{L_1'}\frac{dt_1}{2\pi\sqrt{-1}}\,\Gamma_\C(t_1+1)\Gamma_\C(t_1+\tfrac{s_2}{2})\Gamma_\C(-t_1+\tfrac{\ell_1-\ell_2}{2})\Gamma_\C(-t_1+\tfrac{\ell_2-\ell_3}{2})\\
&\times \int_{L_2'}\frac{dt_2}{2\pi\sqrt{-1}}\,\Gamma_\C(t_2)\Gamma_\C(t_1+t_2)\Gamma_\C(t_2+\tfrac{-s_2+2}{2})\Gamma_\C(-t_2+\tfrac{2+\ell_1-\ell_2}{2})\Gamma_\C(-t_2+\tfrac{2+\ell_2-\ell_3}{2})\Gamma_\C(t_1+t_2+1)^{-1}.
\end{align*}
Here the last equality follows from first Barnes lemma.
Hence $
\<W_{(\underline{\ell},{\sf w}(\itSigma);\,(\ell_2,\ell_2),0)},W_{(\underline{\ell}^\vee,-{\sf w}(\itSigma);\,(-\ell_2,-\ell_2),0)}\>^{(3)}$
is equal to 
\begin{align*}
&\Gamma_\C(1) \int_{L_1'}\frac{dt_1}{2\pi\sqrt{-1}}\,\Gamma_\C(t_1+1)\Gamma_\C(-t_1+\tfrac{\ell_1-\ell_2}{2})\Gamma_\C(-t_1+\tfrac{\ell_2-\ell_3}{2})\\
&\quad\quad\quad\times \int_{L_2'}\frac{dt_2}{2\pi\sqrt{-1}}\,\frac{\Gamma_\C(t_2)\Gamma_\C(t_1+t_2)\Gamma_\C(-t_2+\tfrac{2+\ell_1-\ell_2}{2})\Gamma_\C(-t_2+\tfrac{2+\ell_2-\ell_3}{2})}{\Gamma_\C(t_1+t_2+1)}\\
&\quad\quad\quad\times \int_{L_2}\frac{ds_2}{2\pi\sqrt{-1}}\,\Gamma_\C(\tfrac{s_2+\ell_1-\ell_3}{2})\Gamma_\C(\tfrac{s_2}{2}+t_1)\Gamma_\C(\tfrac{-s_2+4+\ell_1-\ell_3}{2})\Gamma_\C(\tfrac{-s_2}{2}+1+t_2)\\
& = \Gamma_\C(1)\Gamma_\C(2+\ell_1-\ell_3)\\
&\times\int_{L_1'}\frac{dt_1}{2\pi\sqrt{-1}}\,\Gamma_\C(t_1+1)\Gamma_\C(t_1+\tfrac{4+\ell_1-\ell_3}{2})\Gamma_\C(-t_1+\tfrac{\ell_1-\ell_2}{2})\Gamma_\C(-t_1+\tfrac{\ell_2-\ell_3}{2})\\
&\times \int_{L_2'}\frac{dt_2}{2\pi\sqrt{-1}}\,\frac{\Gamma_\C(t_2)\Gamma_\C(t_2+t_1)\Gamma_\C(t_2+\tfrac{2+\ell_1-\ell_3}{2})\Gamma_\C(-t_2+\tfrac{2+\ell_1-\ell_2}{2})\Gamma_\C(-t_2+\tfrac{2+\ell_2-\ell_3}{2})}{\Gamma_\C(t_1+t_2+3+\ell_1-\ell_3)}\\
& = 8\cdot\frac{\Gamma_\C(1)\Gamma_\C(2+\ell_1-\ell_3)\Gamma_\C(1+\tfrac{\ell_1-\ell_2}{2})\Gamma_\C(1+\tfrac{\ell_2-\ell_3}{2})\Gamma_\C(2+\tfrac{2\ell_1-\ell_2-\ell_3}{2})\Gamma_\C(2+\tfrac{\ell_1+\ell_2-2\ell_3}{2})}{\Gamma_\C(3+\ell_1-\ell_3)}\\
&\times\int_{L_1'}\frac{dt_1}{2\pi\sqrt{-1}}\,\frac{\Gamma_\C(t_1+1)\Gamma_\C(t_1+\tfrac{2+\ell_1-\ell_2}{2})\Gamma_\C(t_1+\tfrac{2+\ell_2-\ell_3}{2})\Gamma_\C(-t_1+\tfrac{\ell_1-\ell_2}{2})\Gamma_\C(-t_1+\tfrac{\ell_2-\ell_3}{2})}{\Gamma_\C(t_1+3+\ell_1-\ell_3)}\\
& = 16\cdot \frac{\Gamma_\C(1)\Gamma_\C(1+\ell_1-\ell_2)\Gamma_\C(1+\ell_2-\ell_3)\Gamma_\C(1+\tfrac{\ell_1-\ell_2}{2})^2\Gamma_\C(1+\tfrac{\ell_2-\ell_3}{2})^2\Gamma_\C(1+\tfrac{\ell_1-\ell_3}{2})^2}{\Gamma_\C(3+\ell_1-\ell_3)}.
\end{align*}
Here the first equality  and the last two equalities follow from Barnes' first and second lemmas, respectively.
Finally, note that
\begin{align*}
L(s,\itSigma_v \times \itSigma_v^\vee) = \Gamma_\C(s)^3\Gamma_\C(s+\tfrac{\ell_1-\ell_2}{2})^2\Gamma_\C(s+\tfrac{\ell_2-\ell_3}{2})^2\Gamma_\C(s+\tfrac{\ell_1-\ell_3}{2})^2.
\end{align*}
This completes the proof.
\end{proof}

\begin{thm}\label{T:adjoint}
\noindent
\begin{itemize}
\item[(1)]
We have
\begin{align*}
&\sigma\left(\frac{{\rm Res}_{s=1}L^{(\infty)}(s,\itPi \times \itPi^\vee)}{G(\omega_{\E/\F})\cdot{\rm Reg}_\E\cdot\pi^{4d+\sum_{v\in S_\infty}(\kappa_{1,v}-\kappa_{2,v})}\cdot p^b(\itPi)\cdot p^t(\itPi^\vee)}\right)\\
& = \frac{{\rm Res}_{s=1}L^{(\infty)}(s,{}^\sigma\!\itPi \times {}^\sigma\!\itPi^\vee)}{G(\omega_{\E/\F})\cdot{\rm Reg}_\E\cdot\pi^{4d+\sum_{v\in S_\infty}(\kappa_{1,v}-\kappa_{2,v})}\cdot p^b({}^\sigma\!\itPi)\cdot p^t({}^\sigma\!\itPi^\vee)}
\end{align*}
for all $\sigma\in{\rm Aut}(\C)$.
\item[(2)]
We have
\begin{align*}
&\sigma\left(\frac{{\rm Res}_{s=1}L^{(\infty)}(s,\itSigma \times \itSigma^\vee)}{(\sqrt{-1})^d\cdot{\rm Reg}_\E\cdot\pi^{9d+\sum_{v \in S_\infty}(2\ell_{1,v}-2\ell_{3,v})}\cdot p^b(\itSigma)\cdot p^t(\itSigma^\vee)}\right) \\
& = \frac{{\rm Res}_{s=1}L^{(\infty)}(s,\itSigma \times \itSigma^\vee)}{(\sqrt{-1})^d\cdot{\rm Reg}_\E\cdot\pi^{9d+\sum_{v \in S_\infty}(2\ell_{1,v}-2\ell_{3,v})}\cdot p^b({}^\sigma\!\itSigma)\cdot p^t({}^\sigma\!\itSigma^\vee)}
\end{align*}
for all $\sigma \in {\rm Aut}(\C)$.
\end{itemize}
\end{thm}

\begin{proof}
By Lemma \ref{L:Gauss sum} and Theorem \ref{T:BR}, to prove the first assertion, it suffices to show that
\begin{align*}
B([\itPi_v]_b,[\itPi_v^\vee]_t)_v^{(2)}
=B([\itPi_v^c]_b,[(\itPi_v^c)^\vee]_t)_v^{c,(2)} \in \sqrt{-1}\cdot \pi^{-2-\kappa_{1,v}+\kappa_{2,v}}\cdot \Q^\times
\end{align*}
for all $v \in S_\infty$. Fix $v \in S_\infty$.
We drop the subscript $v$ for brevity.
Let 
\[
(\xi_t^{(2)})^\vee : V_{\underline{\kappa}}^{(2)}\longrightarrow \wedge^2(\frak{g}_{2,\C} / \frak{k}_{2,\C})^* \otimes M_{\lambda,\C},\quad (\xi_t^{c,(2)})^\vee : V_{\underline{\kappa}^\vee}^{(2)}\longrightarrow \wedge^2(\frak{g}_{2,\C} / \frak{k}_{2,\C})^* \otimes M_{\lambda^c,\C}
\]
be the ${\rm U}(2)$-equivariant homomorphisms defined similarly to (\ref{E:U(2) equiv.}) and (\ref{E:U(2) equiv. 2}).  
Note that
\[
\xi_b^{(2)}(y^{\kappa_1-\kappa_2}) = \xi_b^{(2)}\left(\rho_{\underline{\kappa}^\vee}\left(\bp0 &1 \\ -1 &0 \ep \right)\cdot x^{\kappa_1-\kappa_2}\right)=
(-1)^{\lambda_1-\lambda_2}\cdot Y_{(1,-1)}^* \otimes (y^{\lambda_1-\lambda_2}\otimes x^{\lambda_1-\lambda_2}).
\]
Thus
\begin{align*}
&(s^{(2)} \otimes \<\cdot,\cdot\>_{\lambda_{\iota},\C}\otimes\<\cdot,\cdot\>_{\lambda_{\overline{\iota}},\C}) \left( \xi_{b}^{(2)}(y^{\kappa_{1}-\kappa_{2}}),\,(\xi_{t}^{(2)})^\vee(x^{\kappa_{1}-\kappa_{2}})\right)\\
& = (-1)^{\lambda_1-\lambda_2}\cdot s^{(2)}\left(Y_{(1,-1)}^*,\, Y_{(0,0)}^*\wedge Y_{(-1,1)}^*\right)\cdot\<y^{\lambda_1-\lambda_2},x^{\lambda_1-\lambda_2}\>_{\lambda_\iota,\C}\cdot\<x^{\lambda_1-\lambda_2},y^{\lambda_1-\lambda_2}\>_{\lambda_{\overline{\iota}},\C}\\
& = (-1)^{\lambda_1-\lambda_2}\cdot2^3\cdot \sqrt{-1} \cdot\<y^{\lambda_1-\lambda_2},x^{\lambda_1-\lambda_2}\>_{\lambda_\iota,\C}\cdot\<x^{\lambda_1-\lambda_2},y^{\lambda_1-\lambda_2}\>_{\lambda_{\overline{\iota}},\C}.
\end{align*}
Therefore, by the ${\rm U}(2)$-equivariance of the pairing and Lemma \ref{L:BR 1}, we have
\begin{align*}
B([\itPi_v]_b,[\itPi_v^\vee]_t)^{(2)}
&= (\kappa_{1}-\kappa_{2}+1)\cdot \< W_{(\underline{\kappa},{\sf w}(\itPi);\,0)}, W_{(\underline{\kappa}^\vee,-{\sf w}(\itPi);\,\kappa_{1}-\kappa_{2})}\>^{(2)} \\
&\times (s^{(2)} \otimes \<\cdot,\cdot\>_{\lambda_{\iota},\C}\otimes\<\cdot,\cdot\>_{\lambda_{\overline{\iota}},\C}) \left( \xi_{b}^{(2)}(y^{\kappa_{1}-\kappa_{2}}),\,\xi_{t}^{(2)}(x^{\kappa_{1}-\kappa_{2}})\right)\\
& = (-1)^{\lambda_1-\lambda_2}\cdot2^5 \cdot\sqrt{-1}\cdot\<x^{\lambda_1-\lambda_2},y^{\lambda_1-\lambda_2}\>_{\lambda_\iota,\C}\cdot\<y^{\lambda_1-\lambda_2},x^{\lambda_1-\lambda_2}\>_{\lambda_{\overline{\iota}},\C}\\
&\times\frac{\Gamma_\C(1+\kappa_1-\kappa_2)}{\Gamma_\C(1)\Gamma_\C(2+\kappa_1-\kappa_2)}\cdot L(1,\itPi_v \times \itPi_v^\vee)\\
& \in \sqrt{-1}\cdot \pi^{-2-\kappa_{1}+\kappa_{2}}\cdot \Q^\times.
\end{align*}
Moreover, it is clear that
\[
\<x^{\lambda_1-\lambda_2},y^{\lambda_1-\lambda_2}\>_{\lambda_\iota,\C}\cdot\<y^{\lambda_1-\lambda_2},x^{\lambda_1-\lambda_2}\>_{\lambda_{\overline{\iota}},\C} = \<y^{\lambda_1-\lambda_2},x^{\lambda_1-\lambda_2}\>_{\lambda_\iota,\C}\cdot\<x^{\lambda_1-\lambda_2},y^{\lambda_1-\lambda_2}\>_{\lambda_{\overline{\iota}},\C}.
\]
It then follows from Lemma \ref{L:BR 1} again that
\[
B([\itPi_v]_b,[\itPi_v^\vee]_t)^{(2)}= B([\itPi_v^c]_b,[(\itPi_v^c)^\vee]_t)^{c,(2)}.
\]
This completes the proof of the first assertion.

By Theorem \ref{T:BR}, to prove the second assertion, it suffices to show that
\begin{align*}
B([\itSigma_v]_b,[\itSigma_v^\vee]_t)_v^{(3)}
=B([\itSigma_v^c]_b,[(\itSigma_v^c)^\vee]_t)_v^{c,(3)} \in \sqrt{-1}\cdot \pi^{-6-2\ell_{1,v}+2\ell_{3,v}}\cdot \Q^\times
\end{align*}
for all $v \in S_\infty$. Fix $v \in S_\infty$.
We drop the subscript $v$ for brevity.
First note that for any ${\rm U}(3)$-equivariant bilinear pairing 
\[
\<\cdot,\cdot\> : V_{\underline{\ell}}^{(3)} \times V_{\underline{\ell}^\vee}^{(3)} \longrightarrow \C,
\]
our normalization in \S\,\ref{SS:rep GL3} implies that
\begin{align}\label{E:adjoint proof 1}
\begin{split}
\<{\bf v}_{(\underline{\ell};\,\underline{\kappa}',i)},\,{\bf v}_{(\underline{\ell}^\vee;\,(\underline{\kappa}')^\vee,\kappa_1'-\kappa_2'-i)}\> 
&= C_{(\underline{\ell};\, \underline{\kappa}')}^{-1}\cdot(-1)^i\cdot{\kappa_1'-\kappa_2'\choose i}^{-1}\cdot \<{\bf v}_{(\underline{\ell};\,(\ell_2,\ell_2),0)},\,{\bf v}_{(\underline{\ell}^\vee;\,(-\ell_2,-\ell_2),0)}\>
\end{split}
\end{align}
for all $\underline{\ell} \succ \underline{\kappa}'$ and $0 \leq i \leq \kappa_{1}'-\kappa_{2}'$.
Let 
\[
(\xi_t^{(3)})^\vee : V_{\underline{\ell}}^{(3)}\longrightarrow \wedge^5(\frak{g}_{3,\C} / \frak{k}_{3,\C})^* \otimes M_{\mu,\C},\quad (\xi_t^{c,(3)})^\vee : V_{\underline{\ell}^\vee}^{(3)}\longrightarrow \wedge^5(\frak{g}_{3,\C} / \frak{k}_{3,\C})^* \otimes M_{\mu^c,\C}
\]
be the ${\rm U}(3)$-equivariant homomorphisms defined similarly to (\ref{E:U(3) equiv.}) and (\ref{E:U(3) equiv. 2}).
Note that
\begin{align*}
\xi_b^{(3)}({\bf v}_{(\underline{\ell}^\vee;\,(-\ell_2,-\ell_1),\ell_1-\ell_2)})) &= (-1)^{\ell_1}\cdot \xi_b^{(3)} \left(\rho_{\underline{\ell}^\vee}\left(\bp0 &0&1 \\0& 1 &0 \\1&0&0\ep \right)\cdot {\bf v}_{(\underline{\ell}^\vee;\,(-\ell_3,-\ell_2),0)}\right)\\
&=(-1)^{1+\mu_{1}-\mu_{3}}\cdot X_{(0,1,-1)}^*\wedge X_{(1,-1,0)}^* \wedge X_{(1,0,-1)}^* \\
&\otimes \left({\bf v}_{(\mu_{\iota}^\vee;\,(-\mu_{2},-\mu_{1}),\mu_{1}-\mu_{2})}\otimes {\bf v}_{(\mu_{\overline{\iota}}^\vee;\,(\mu_{1}-{\sf w}(\itSigma),\mu_{2}-{\sf w}(\itSigma)),0)}\right).
\end{align*}
Thus
\begin{align*}
&(s^{(3)} \otimes \<\cdot,\cdot\>_{\mu_{\iota},\C}\otimes\<\cdot,\cdot\>_{\mu_{\overline{\iota}},\C}) \left( \xi_{b}^{(3)}({\bf v}_{(\underline{\ell}^\vee;\,(-\ell_2,-\ell_1),\ell_1-\ell_2)}),\,\xi_{t}^{(3)}({\bf v}_{(\underline{\ell};\,(\ell_1,\ell_2),0)})\right)\\
& = (-1)^{1+\mu_{1}-\mu_{3}}\cdot s^{(3)}\left( X_{(0,1,-1)}^*\wedge X_{(1,-1,0)}^* \wedge X_{(1,0,-1)}^*,\,Z_{12}^*\wedge Z_{23}^* \wedge X_{(-1,1,0)}^*\wedge X_{(0,-1,1)}^* \wedge X_{(-1,0,1)}^*\right)\\
&\times\<{\bf v}_{(\mu_{\iota}^\vee;\,(-\mu_{2},-\mu_{1}),\mu_{1}-\mu_{2})},\,{\bf v}_{(\mu_{\iota};\,(\mu_1,\mu_2),0)}\>_{\mu_\iota,\C}\\
&\times\<{\bf v}_{(\mu_{\overline{\iota}}^\vee;\,(\mu_{1}-{\sf w}(\itSigma),\mu_{2}-{\sf w}(\itSigma)),0)},\,{\bf v}_{(\mu_{\overline{\iota}};\,(-\mu_{2}+{\sf w}(\itSigma),-\mu_{1}+{\sf w}(\itSigma)),\mu_{1}-\mu_{2})}\>_{\mu_{\overline{\iota}},\C}\\
& = (-1)^{1+\mu_1-\mu_3}\cdot2^8\cdot\sqrt{-1}\cdot\<{\bf v}_{(\mu_{\iota}^\vee;\,(-\mu_{2},-\mu_{1}),\mu_{1}-\mu_{2})},\,{\bf v}_{(\mu_{\iota};\,(\mu_1,\mu_2),0)}\>_{\mu_\iota,\C}\\
&\times\<{\bf v}_{(\mu_{\overline{\iota}}^\vee;\,(\mu_{1}-{\sf w}(\itSigma),\mu_{2}-{\sf w}(\itSigma)),0)},\,{\bf v}_{(\mu_{\overline{\iota}};\,(-\mu_{2}+{\sf w}(\itSigma),-\mu_{1}+{\sf w}(\itSigma)),\mu_{1}-\mu_{2})}\>_{\mu_{\overline{\iota}},\C}.
\end{align*}
Therefore, by the ${\rm U}(3)$-equivariance of the pairing, Lemma \ref{L:BR 2}, and (\ref{E:adjoint proof 1}), we have
\begin{align*}
&B([\itSigma_v]_b,[\itSigma_v^\vee]_t)^{(3)} \\
& = {\rm dim}_\C\,V_{\underline{\ell}}\cdot C_{(\underline{\ell};\,(\ell_1,\ell_2))}\cdot \<W_{(\underline{\ell},{\sf w}(\itSigma);\,(\ell_2,\ell_2),0)},W_{(\underline{\ell}^\vee,-{\sf w}(\itSigma);\,(-\ell_2,-\ell_2),0)}\>^{(3)}\\
&\times(s^{(3)} \otimes \<\cdot,\cdot\>_{\mu_{\iota},\C}\otimes\<\cdot,\cdot\>_{\mu_{\overline{\iota}},\C}) \left( \xi_{b}^{(3)}({\bf v}_{(\underline{\ell}^\vee;\,(-\ell_2,-\ell_1),\ell_1-\ell_2)}),\,\xi_{t}^{(3)}({\bf v}_{(\underline{\ell};\,(\ell_1,\ell_2),0)})\right)\\
&= {\rm dim}_\C\,V_{\underline{\ell}}\cdot C_{(\underline{\ell};\,(\ell_1,\ell_2))}\cdot (-1)^{1+\mu_1-\mu_3}\cdot2^{12}\cdot\sqrt{-1}\cdot\<{\bf v}_{(\mu_{\iota}^\vee;\,(-\mu_{2},-\mu_{1}),\mu_{1}-\mu_{2})},\,{\bf v}_{(\mu_{\iota};\,(\mu_1,\mu_2),0)}\>_{\mu_\iota,\C}\\
&\times\<{\bf v}_{(\mu_{\overline{\iota}}^\vee;\,(\mu_{1}-{\sf w}(\itSigma),\mu_{2}-{\sf w}(\itSigma)),0)},\,{\bf v}_{(\mu_{\overline{\iota}};\,(-\mu_{2}+{\sf w}(\itSigma),-\mu_{1}+{\sf w}(\itSigma)),\mu_{1}-\mu_{2})}\>_{\mu_{\overline{\iota}},\C}\\
&\times\frac{\Gamma_\C(1+\ell_{1}-\ell_{2})\Gamma_\C(1+\ell_{2}-\ell_{3})}{\Gamma_\C(1)^2\Gamma_\C(3+\ell_{1}-\ell_{3})}\cdot L(1,\itSigma_v \times \itSigma_v^\vee)\\
& \in \sqrt{-1}\cdot \pi^{-6-2\ell_1+2\ell_3}\cdot \Q^\times.
\end{align*}
Finally, since
\begin{align*}
\rho_{\underline{\ell}}\left(\bp0 &0&1 \\0& 1 &0 \\1&0&0\ep \right)\cdot {\bf v}_{(\underline{\ell};\,(\ell_1,\ell_2),0)} &= (-1)^{\ell_3}\cdot {\bf v}_{(\underline{\ell};\,(\ell_2,\ell_3),\ell_2-\ell_3)},\\
\rho_{\underline{\ell}^\vee}\left(\bp0 &0&1 \\0& 1 &0 \\1&0&0\ep \right)\cdot {\bf v}_{(\underline{\ell}^\vee;\,(-\ell_2,-\ell_1),\ell_1-\ell_2)} &= (-1)^{\ell_1}\cdot {\bf v}_{(\underline{\ell}^\vee;\,(-\ell_3,-\ell_2),0)},
\end{align*}
by (\ref{E:adjoint proof 1}) we have $C_{(\underline{\ell};\,(\ell_1,\ell_2))} = C_{(\underline{\ell};\,(\ell_2,\ell_3))}$.
Thus $C_{(\underline{\ell};\,(\ell_1,\ell_2))}=C_{(\underline{\ell}^\vee;\,(-\ell_3,-\ell_2))}$ by (\ref{E:scalar relation}).
Similarly, we have 
\begin{align*}
&\<{\bf v}_{(\mu_{\iota}^\vee;\,(-\mu_{2},-\mu_{1}),\mu_{1}-\mu_{2})},\,{\bf v}_{(\mu_{\iota};\,(\mu_1,\mu_2),0)}\>_{\mu_\iota,\C}\\
& = (-1)^{\mu_1-\mu_3}\cdot\<{\bf v}_{(\mu_{\iota}^\vee;\,(-\mu_{3},-\mu_{2}),0)},\,{\bf v}_{(\mu_{\iota};\,(\mu_2,\mu_3),\mu_{2}-\mu_{3})}\>_{\mu_\iota,\C},\\
&\<{\bf v}_{(\mu_{\overline{\iota}}^\vee;\,(\mu_{1}-{\sf w}(\itSigma),\mu_{2}-{\sf w}(\itSigma)),0)},\,{\bf v}_{(\mu_{\overline{\iota}};\,(-\mu_{2}+{\sf w}(\itSigma),-\mu_{1}+{\sf w}(\itSigma)),\mu_{1}-\mu_{2})}\>_{\mu_{\overline{\iota}},\C}\\
& = (-1)^{\mu_1-\mu_3}\cdot\<{\bf v}_{(\mu_{\overline{\iota}}^\vee;\,(\mu_{2}-{\sf w}(\itSigma),\mu_{3}-{\sf w}(\itSigma)),\mu_{2}-\mu_{3})},\,{\bf v}_{(\mu_{\overline{\iota}};\,(-\mu_{3}+{\sf w}(\itSigma),-\mu_{2}+{\sf w}(\itSigma)),0)}\>_{\mu_{\overline{\iota}},\C}.
\end{align*}
It then follows from Lemma \ref{L:BR 2} again that
\[
B([\itSigma_v]_b,[\itSigma_v^\vee]_t)^{(3)}= B([\itSigma_v^c]_b,[(\itSigma_v^c)^\vee]_t)^{c,(3)}.
\]
This completes the proof.
\end{proof}

\section{Algebraicity of the Asai $L$-functions for $\GL_2$ and $\GL_3$}\label{S:Asai}

Let $\E$ be a totally imaginary quadratic extension of $\F$.
\subsection{Algebraicity of the Asai $L$-functions for $\GL_n$}

Let 
$
\delta_{(n),\R} : \GL_n(\R) \longrightarrow \GL_n(\C)
$
be the natural inclusion. 
Denote by $\frak{g}_n'$ and $\frak{k}_n'$ the Lie algebras of the real Lie groups $\GL_n(\R)$ and $K_n'=K_n \cap \GL_n(\R)$, respectively.
Let
\begin{align*}
s_\R^{(n)} : \wedge^{t_n} (\frak{g}_{n,\C}/\frak{k}_{n,\C})^*\longrightarrow \C
\end{align*}
be the linear functionals defined by
\begin{align*}
\wedge_{i=1}^{t_n}\delta_{(n),\R}^*X_i^* &= s_\R^{(n)}(\wedge_{i=1}^{t_n}X_i^*)\cdot (\wedge_{i=1}^{n-1} \,e_{ii}^*) \wedge (\wedge_{1\leq i <j \leq n}\,e_{ij}^* ).
\end{align*} 
Here $\{e_{kk},e_{i,j}\,\vert\, 1\leq k \leq n-1,\,1\leq i<j \leq n \} \subset \frak{g}_n'$ is ordered lexicographically and
\[
\delta_{(n),\R}^* :  (\frak{g}_{n,\C}/\frak{k}_{n,\C})^* \longrightarrow  (\frak{g}_{n,\C}'/\frak{k}_{n,\C}')^*
\]
is the homomorphism induced by $\delta_{(n),\R}$.
Note that $(\wedge_{i=1}^{n-1} \,e_{ii}^*) \wedge (\wedge_{1\leq i <j \leq n}\,e_{ij}^* )$ corresponds to the standard measure on $\GL_n(\R)/(K_n \cap \GL_n(\R))$ defined in \S\,\ref{SS:measure}.

Let $\itPi$ be a cohomological irreducible cuspidal automorphic representation of $\GL_n(\A_\E)$ with central character $\omega_\itPi$. We have $|\omega_\itPi| = |\mbox{ }|_{\A_\E}^{n{\sf w}/2}$ for some ${\sf w} \in \Z$.
Let 
\[
\mu = \prod_{v \in S_\infty} \mu_v \in \prod_{v \in S_\infty} (X^+(T_n) \times X^+(T_n))
\]
with $\mu_v = (\mu_{\iota_v},\mu_{\overline{\iota}_v})$ such that $\itPi$ contributes to the cuspidal cohomology of $\GL_n(\A_\E)$ with coefficients in $M^\vee_{\mu,\C}$. 
We keep the notation of \S\,\ref{SS:Whittaker GL_n}.
With respect to the choice of the generators (\ref{E:generators}) in the cuspidal cohomology, let $p^b({}^\sigma\!\itPi)$ and $p^t({}^\sigma\!\itPi)$ be the Whittaker periods of ${}^\sigma\!\itPi$ for each $\sigma\in{\rm Aut}(\C)$ defined as in Lemma \ref{L:Whittaker periods}.
Note that the periods are normalized so that the diagram in Lemma \ref{L:Whittaker periods} commutes.
We assume further that $\itPi$ is essentially conjugate self-dual, that is, 
\[
\itPi^\vee = \itPi^c \otimes \chi\circ {\rm N}_{\K/\F}
\]
for some algebraic Hecke character $\chi$ of $\A_\F^\times$ with parallel signature ${\rm sgn}(\chi)$.
After replacing $\chi$ by $\omega_{\E/\F}\cdot\chi$ if necessary, we may assume that
\begin{align}\label{E:sign normalization 1}
{\rm sgn}(\chi) = (-1)^{{\sf w}(\itPi)}.
\end{align}
In particular, we have
\[
\chi_v = {\rm sgn}^{\sf w}|\mbox{ }|_\R^{-{\sf w}}
\]
for all $v \in S_\infty$.
Let $v \in S_\infty$.
For $\alpha \in \C^\times$ with ${\rm tr}_{\C/\R}(\alpha)=0$, let 
\[
Z_v^{(n)}(\mbox{ };\alpha) \in {\rm Hom}_{(\frak{g}_n',K_n')}(\mathcal{W}(\itPi_v,\psi_{n,\C}),{\rm sgn}^{n-1}\cdot\chi_v^{-1})
\]
defined by (cf.\,\cite{Kemarsky2015})
\begin{align*}
Z_v^{(n)}(W;\alpha) &= \int_{N_{n-1}(\R)\backslash\GL_{n-1}(\R)}W\left({\rm diag}(1,\alpha,\cdots,\alpha^{n-1})\bp g & 0 \\ 0 & 1\ep \right){\rm sgn}^{n-1}(\det(g))\chi_v(\det(g))\,dg.
\end{align*}
Here $dg$ is the standard measure defined in \S\,\ref{SS:measure} and note that the integral converges absolutely.
We simply write $Z_v^{(n)}(\mbox{ };\alpha_v) = Z_v^{(n)}$ when $\alpha = \sqrt{-1}$.
Similarly we define
\[
Z_v^{c,(n)}(\mbox{ };\alpha) \in {\rm Hom}_{(\frak{g}_n',K_n')}(\mathcal{W}(\itPi_v^c,\psi_{n,\C}),{\rm sgn}^{n-1}\cdot\chi_v^{-1}).
\]
Fix a non-zero functional
\begin{align*}
f_{\mu_v,\C}\in{\rm Hom}_{\GL_n(\C)}(M_{\mu_{\iota_v},\C}^\vee \otimes M_{\mu_{\overline{\iota}_v},\C}^\vee,\,{\rm det}^{-{\sf w}})
\end{align*}
defined over $\Q$.
Define
$
f_{\mu_v^c,\C}\in{\rm Hom}_{\GL_n(\C)}(M_{\mu_{\overline{\iota}_v},\C}^\vee \otimes M_{\mu_{{\iota}_v},\C}^\vee,\,{\rm det}^{-{\sf w}})
$
by
\[
f_{\mu_v^c,\C} = f_{\mu_v,\C} \circ ({\bf c}_{\mu_v^\vee}^{(n)})^{-1}.
\]
Then we have
\begin{align*}
f_{\mu_v,\C} \vert_{\GL_n(\R)}\in{\rm Hom}_{\GL_n(\R)}(M_{\mu_v,\C}^\vee ,{\rm det}^{-{\sf w}}),\quad f_{\mu_v^c,\C} \vert_{\GL_n(\R)}\in{\rm Hom}_{\GL_n(\R)}(M_{\mu_v^c,\C}^\vee ,{\rm det}^{-{\sf w}}).
\end{align*}
Let
\begin{align*}
L_v^{(n)}(\mbox{ };\alpha): \mathcal{W}(\itPi_v,\psi_{n,\C}) \otimes \wedge^{t_n} (\frak{g}_{n,\C} / \frak{k}_{n,\C})^* \otimes M_{\mu_v,\C}^\vee &\longrightarrow \C,\\ 
L_v^{c,(n)}(\mbox{ };\alpha): \mathcal{W}(\itPi_v^c,\psi_{n,\C}) \otimes \wedge^{t_n} (\frak{g}_{n,\C} / \frak{k}_{n,\C})^* \otimes M_{\mu_v^c,\C}^\vee &\longrightarrow \C
\end{align*}
be the linear functionals defined by
\begin{align*}
L_v^{(n)}(W_1 \otimes  \wedge_{i=1}^{t_n} X_i^* \otimes ({\bf u}_{\iota_v}\otimes{\bf u}_{\overline{\iota}_v});\alpha) &= Z_v^{(n)}(W_1;\alpha)\cdot s_\R^{(n)}(\wedge_{i=1}^{t_n} X_i^*)\cdot f_{\mu_v,\C}({\bf u}_{\iota_v}\otimes{\bf u}_{\overline{\iota}_v}),\\
L_v^{c,(n)}(W_2 \otimes  \wedge_{i=1}^{t_n} Y_i^* \otimes ({\bf u}_{\overline{\iota}_v}\otimes{\bf u}_{\iota_v});\alpha) &= Z_v^{c,(n)}(W_2;\alpha)\cdot s_\R^{(n)}(\wedge_{i=1}^{t_n} Y_i^*)\cdot f_{\mu_v^c,\C}({\bf u}_{\overline{\iota}_v}\otimes{\bf u}_{\iota_v}).
\end{align*}
For $\sigma \in {\rm Aut}(\C)$ and $(\alpha_v)_{v \in S_\infty} \subset (\C^\times)^{S_\infty}$ with ${\rm tr}_{\C/\R}(\alpha_v)=0$, we then have the linear functional
\begin{align*}
L_{\sigma,\infty}^{(n)}(\mbox{ };(\alpha_v)_{v \in S_\infty})
= \left( \bigotimes_{v \in S_\infty,\,\sigma^{-1}\circ\iota_v = \iota_{\sigma^{-1}\circ v}}L_{\sigma^{-1}\circ v}^{(n)}(\mbox{ };\alpha_v) \right) \otimes \left( \bigotimes_{v \in S_\infty,\,\sigma^{-1}\circ\iota_v = \overline{\iota}_{\sigma^{-1}\circ v}}L_{\sigma^{-1}\circ v}^{c,(n)}(\mbox{ };\alpha_v) \right)
\end{align*}
on
\[
\mathcal{W}({}^\sigma\!\itPi_\infty,\psi_{n,\E_\infty}) \otimes \wedge^{t_nd} (\frak{g}_{n,\infty,\C} / \frak{k}_{n,\infty,\C})^* \otimes M_{{}^\sigma\!\mu,\C}^\vee.
\]
We simply write $L_{\sigma,\infty}^{(n)}(\mbox{ };(\alpha_v)_{v \in S_\infty}) = L_{\sigma,\infty}^{(n)}$ when $\alpha_v = \sqrt{-1}$ for all $v \in S_\infty$, and write $L_{\sigma,\infty}^{(n)}(\mbox{ };(\alpha_v)_{v \in S_\infty}) = L_{\infty}^{(n)}(\mbox{ };(\alpha_v)_{v \in S_\infty})$ when $\sigma$ is the identity map.

We have the following theorem on the algebraicity of 
\[
{\rm Res}_{s=1}L^{(\infty)}(s,\itPi,{\rm As}^{(-1)^{n-1}}\otimes\chi)
\]
in terms of the top degree Whittaker period of $\itPi$. The result is a generalization of the result of Grobner, Harris, and Lapid \cite[Theorem 6.4]{GHL2016} to essentially conjugate self-dual representations.
We also recall the non-vanishing result of Sun \cite[Theorem 4.5]{Sun2019}.

\begin{thm}\label{T:GHL}
\noindent
\begin{itemize}
\item[(1)]
We have
\begin{align*}
&\sigma\left(\frac{{\rm Res}_{s=1} L^{(\infty)}(s,\itPi,{\rm As}^{(-1)^{n-1}}\otimes\chi)}
{|D_\F|^{\tfrac{n(n+1)}{4}}\cdot |D_\E|^{\mbox{$\tfrac{n(n-1)(n-2)}{12}$}+\tfrac{{\sf w}\cdot n(n-1)}{4}}\cdot G(\omega_{\E/\F}^{n-1}\cdot \chi)^{\tfrac{n(n-1)}{2}}\cdot {\rm Reg}_\F\cdot \pi^{d\lfloor \tfrac{n}{2} \rfloor}\cdot p^t(\itPi)}\cdot L_\infty^{(n)}([\itPi_\infty]_t)\right)\\
& = \frac{{\rm Res}_{s=1}L^{(\infty)}(s,{}^\sigma\!\itPi,{\rm As}^{(-1)^{n-1}}\otimes{}^\sigma\!\chi)}{|D_\F|^{\tfrac{n(n+1)}{4}}\cdot |D_\E|^{\mbox{$\tfrac{n(n-1)(n-2)}{12}$}+\tfrac{{\sf w}\cdot n(n-1)}{4}}\cdot G(\omega_{\E/\F}^{n-1}\cdot{}^\sigma\!\chi)^{\tfrac{n(n-1)}{2}}\cdot {\rm Reg}_\F\cdot \pi^{d\lfloor \tfrac{n}{2} \rfloor}\cdot p^t({}^\sigma\!\itPi)}\cdot L_{\sigma,\infty}^{(n)}([{}^\sigma\!\itPi_\infty]_t)
\end{align*}
for all $\sigma \in {\rm Aut}(\C)$.
\item[(2)]{\rm (Sun)} The archimedean factor $L_\infty^{(n)}([\itPi_\infty]_t)$ is non-zero.
\end{itemize}
\end{thm}

\begin{proof}

For the second assertion, it is clear that $L_v^{(n)}$ defines a homomorphism
\[
H^{t_n}(\frak{g}_n,K_n;\itPi_v \otimes M_{\mu_v,\C}^\vee)\longrightarrow H^{t_n}(\frak{g}_n',K_n'; {\rm sgn}^{n-1})
\]
for $v \in S_\infty$. Moreover, $L_v^{(n)}([\itPi_v]_t) \neq 0$ if and only if the homomorphism is non-zero. Note that the $L$-factor 
$L(s,\itPi_v,{\rm As}^{(-1)^{n-1}}\otimes\chi_v)$ is holomorphic and non-vanishing at $s=1$. Therefore, by \cite[Theorem 1-(iii)]{Raphael2018} (see also \cite{CCI2020} for $n=2$), the functional $Z_v^{(n)}$ is non-zero. Also note that
\[
{\rm Hom}_{(\frak{g}_n',K_n')}(\itPi_v,{\rm sgn}^{n-1}\cdot\chi_v^{-1}) = {\rm Hom}_{(\frak{g}_n',K_n')}(\itPi_v\otimes\chi_{\sf w}|\mbox{ }|_\C^{-{\sf w}/2},{\rm sgn}^{n-1}),
\]
where $\chi_{\sf w}(z) = (z/\overline{z})^{{\sf w}/2}$. The non-vanishing of $L_v^{(n)}([\itPi_v]_t)$ then follows from the results \cite[Theorems 4.1 and A.3]{Sun2019} of Sun.

For each place $v$ of $\F$, let 
\[
\E_v = \prod_{w \mid v}\E_w,\quad \itPi_v = \bigotimes_{w \mid v}\itPi_w,\quad \psi_{n,\E_v} = \bigotimes_{w \mid v}\psi_{n,\E_w}.
\]
For each finite place $v$ of $\F$ and $\alpha_v \in \E_v^\times$ with ${\rm tr}_{\E_v/ \F_v}(\alpha_v)=0$, let 
\[
Z_v^{(n)}(\mbox{ };\alpha_v) \in {\rm Hom}_{\mathcal{P}_n(\F_v)}(\mathcal{W}(\itPi_v,\psi_{v,\E_v}), \omega_{\E_v/\F_v}^{n-1}\cdot\chi_v^{-1})
\]
defined by
\begin{align*}
Z_v^{(n)}(W_v;\alpha_v)& = \int_{N_{n-1}(\F_v)\backslash \GL_{n-1}(\F_v)}W_v\left({\rm diag}(1,\alpha_v,\cdots,\alpha_v^{n-1})\bp g_v & 0 \\ 0 & 1\ep\right)\omega_{\E_v/\F_v}^{n-1}(\det(g_v))\chi_v(\det(g_v))\,dg_v.
\end{align*}
Here $\mathcal{P}_n$ is the mirabolic subgroup of $\GL_n$ consisting of matrices whose last row is $(0,\cdots,0,1)$, and
$dg_v$ is the standard measure in \S\,\ref{SS:measure}. 
First note that the functional $Z_v^{(n)}(\mbox{ };\alpha_v)$ is also non-zero for all finite places $v$ by \cite[Lemma 4.1]{Offen2011}.
By the integral representation of the Asai $L$-functions \cite{Flicker1988}, \cite{Kable2004}, together with the corresponding local theory partially completed in \cite{Raphael2018}, we see that the order of the pole of $L(s,\itPi,{\rm As}^{(-1)^{n-1}}\otimes\chi)$ at $s=1$ is at most one.
Proceeding similarly as in the proof of \cite[Proposition 3.2]{Zhang2014}, we have the integral representation of the residue at $s=1$ as follows:
For $\varphi \in \itPi$ with $W_{\varphi,\psi_{n,\E}} = \bigotimes_v W_v \in \bigotimes_v \mathcal{W}(\itPi_v,\psi_{v,\E_v})$, we have
\begin{align}\label{E:proof A 4}
\begin{split}
&\int_{\A_\F^\times\GL_n(\F)\backslash \GL_n(\A_\F)}\varphi(g)\omega_{\E/\F}^{n-1}(\det(g))\chi(\det(g))\,dg\\
& = |D_\F|^{n(n+1)/4}\cdot\frac{n\cdot {\rm Res}_{s=1}L(s,\itPi,{\rm As}^{(-1)^{n-1}}\otimes\chi)}{{\rm vol}(\F^\times \backslash \A_\F^1)}\cdot\prod_{v}\zeta_{\F_v}(n)\cdot\frac{Z_v^{(n)}(W_v,;\alpha_v)}{L(1,\itPi_v,{\rm As}^{(-1)^{n-1}}\otimes\chi_v)}.
\end{split}
\end{align}
Here 
\begin{itemize}
\item $dg$ is the standard measure defined in \S\,\ref{SS:measure};
\item ${\rm vol}(\F^\times \backslash \A_\F^1)$ is the volume computed under the standard measure on $\A_\F^\times$  defined in \S\,\ref{SS:measure};
\item $\alpha \in \E^\times$ is any element such that ${\rm tr}_{\E/\F}(\alpha)=0$ and $\alpha_v$ is the image of $\alpha$ in $\E_v$ for each place $v$ of $\F$;
\item $\zeta_{\F_v}(s)$ is the local zeta function of $\F_v$;
\item we have
$
\zeta_{\F_v}(n)\cdot Z_v^{(n)}(W_v,;\alpha_v)=L(1,\itPi_v,{\rm As}^{(-1)^{n-1}}\otimes\chi_v)
$
for all but finitely many $v$.
\end{itemize}
Note that $\alpha_v = \iota_v(\alpha)$ for $v \in S_\infty$, and the local $L$-factor $L(s,\itPi_v,{\rm As}^{(-1)^{n-1}}\otimes\chi_v)$ is holomorphic and non-vanishing at $s=1$ for each place $v$ of $\F$ (cf.\,\cite[Proposition]{FZ1995} and \cite[Lemma 3.3.1]{Raphael2018}). 
We also remark that the extra factors $|D_\F|^{n(n+1)/4}$ and $\zeta_{\F_v}(n)$ comparing with \cite[Proposition 3.2]{Zhang2014} are due to different choices of Haar measures in \S\,\ref{SS:measure} and \cite[\S\,2.1]{Zhang2014}.
Now we consider the locally symmetric spaces
\[
\mathcal{S}_{n,\F} = \GL_n(\F)\backslash \GL_n(\A_\F)/K_{n,\infty}',\quad \mathcal{S}_{n,\E} = \GL_n(\E)\backslash \GL_n(\A_\E)/K_{n,\infty}.
\]
The natural inclusion $j : \GL_n(\F) \rightarrow \GL_n(\E)$ induces a homomorphism 
\[
j^*: H_c^{t_nd}(\mathcal{S}_{n,\E},\mathcal{M}_{\mu,\C}^\vee)\longrightarrow H_c^{t_nd}(\mathcal{S}_{n,\F},\mathcal{M}_{\underline{\sf w},\C}^\vee)
\]
for cohomology with compact support with coefficients in the locally constant sheaves $\mathcal{M}_{\mu,\C}^\vee$ and $\mathcal{M}_{\underline{\sf w},\C}^\vee$ on $\mathcal{S}_{n,\E}$ and $\mathcal{S}_{n,\F}$ induced by the algebraic representations ${M}_{\mu,\C}^\vee$ and ${M}_{\underline{\sf w},\C}^\vee$, respectively.
Here $\underline{\sf w} = \prod_{v \in S_\infty}({\sf w},\cdots,{\sf w}) \in \prod_{v \in S_\infty}X^+(T_n)$.
Moreover, the homomorphism is equivariant under ${\rm Aut}(\C)$, that is, the diagram
\[
\begin{tikzcd}
H_c^{t_nd}(\mathcal{S}_{n,\E},\mathcal{M}_{\mu,\C}^\vee) \arrow[r, "j^*"] \arrow[d, "T_\sigma"] &H_c^{t_nd}(\mathcal{S}_{n,\F},\mathcal{M}_{\underline{\sf w},\C}^\vee)\arrow[d, "T_\sigma"]\\
H_c^{t_nd}(\mathcal{S}_{n,\E},\mathcal{M}_{{}^\sigma\!\mu,\C}^\vee)\arrow[r, "j^*"] & H_c^{t_nd}(\mathcal{S}_{n,\F},\mathcal{M}_{\underline{\sf w},\C}^\vee)
\end{tikzcd}
\]
is commutative for all $\sigma \in {\rm Aut}(\C)$. 
Here the $\sigma$-linear isomorphisms $T_\sigma$ are induced by the natural $\sigma$-linear isomorphisms ${M}_{\mu,\C}^\vee \rightarrow {M}_{{}^\sigma\!\mu,\C}^\vee$ and ${M}_{\underline{\sf w},\C}^\vee \rightarrow {M}_{\underline{\sf w},\C}^\vee$ (cf.\,(\ref{E:sigma linear})).
With respect to the standard measures in \S\,\ref{SS:measure}, we have the Poincar\'e pairing
\[
H_c^{t_nd}(\mathcal{S}_{n,\F},\mathcal{M}_{\underline{\sf w},\C}^\vee) \otimes H^0(\mathcal{S}_{n,\F},\mathcal{M}_{\underline{\sf w},\C}) \longrightarrow \C,\quad c_1\otimes c_2\longmapsto \int_{\mathcal{S}_{n,\F}}c_1\wedge c_2
\]
which is equivariant under ${\rm Aut}(\C)$ (cf.\,\cite[p.\,124]{Clozel1990}).
Let $[\omega_{\E/\F}^{n-1}\cdot\chi] \in H^0(\mathcal{S}_{n,\F},\mathcal{M}_{\underline{\sf w},\C})$ be the class represented by $\omega_{\E/\F}^{n-1}\cdot\chi$. In particular, we have
\begin{align}\label{E:proof A 5}
\sigma\left( \int_{\mathcal{S}_{n,\F}} \frac{j^*\circ\mathcal{F}_t(W)}{p^t(\itPi)} \wedge [\omega_{\E/\F}^{n-1}\cdot\chi] \right) = \int_{\mathcal{S}_{n,\F}} \frac{j^*\circ\mathcal{F}_{\sigma,t}(t_\sigma W)}{p^t({}^\sigma\!\itPi)} \wedge [\omega_{\E/\F}^{n-1}\cdot{}^\sigma\!\chi]
\end{align}
for all $W \in \mathcal{W}(\itPi_f,\psi_{n,\E}^{(\infty)})$ and $\sigma \in {\rm Aut}(\C)$. Here $\mathcal{F}_{\sigma,t}$ is the $\GL_n(\A_{\E,f})$-equivariant isomorphism defined in (\ref{E:Galois equiv}). 
By the integral representation (\ref{E:proof A 4}) and proceed as in \cite[\S\,3.3.3]{BR2017}, we can interpret (\ref{E:proof A 5}) as follows:
We have
\begingroup
\small
\begin{align}\label{E:proof A 1}
\begin{split}
& \sigma \left( \frac{{\rm Res}_{s=1}L^{(\infty)}(s,\itPi,{\rm As}^{(-1)^{n-1}}\otimes\chi)}{|D_\F|^{n(n+1)/4}\cdot{\rm vol}(\F^\times \backslash \A_\F^1)\cdot \pi^{d\lfloor \tfrac{n}{2} \rfloor}\cdot p^t(\itPi)}\cdot\prod_{v\nmid \infty}\zeta_{\F_v}(n)\cdot\frac{Z_v^{(n)}(W_v,;\alpha_v)}{L(1,\itPi_v,{\rm As}^{(-1)^{n-1}}\otimes\chi_v)}\cdot L_{\infty}^{(n)}([\itPi_\infty]_t;(\iota_v(\alpha))_{v \in S_\infty})\right)\\
& =  \frac{{\rm Res}_{s=1}L^{(\infty)}(s,{}^\sigma\!\itPi,{\rm As}^{(-1)^{n-1}}\otimes{}^\sigma\!\chi)}{|D_\F|^{n(n+1)/4}\cdot{\rm vol}(\F^\times \backslash \A_\F^1)\cdot \pi^{d\lfloor \tfrac{n}{2} \rfloor}\cdot p^t({}^\sigma\!\itPi)}\cdot\prod_{v\nmid \infty}\zeta_{\F_v}(n)\cdot\frac{Z_v^{(n)}(t_\sigma W_v,;\alpha_v)}{L(1,{}^\sigma\!\itPi_v,{\rm As}^{(-1)^{n-1}}\otimes{}^\sigma\!\chi_v)}\cdot L_{\sigma,\infty}^{(n)}([{}^\sigma\!\itPi_\infty]_t;(\iota_v(\alpha))_{v \in S_\infty})
\end{split}
\end{align}
\endgroup
for all $\bigotimes_{v \nmid \infty}W_v \in \mathcal{W}(\itPi_f,\psi_{n,\E}^{(\infty)})$ and $\sigma\in{\rm Aut}(\C)$.
Here the factor $\pi^{d\lfloor \tfrac{n}{2} \rfloor}$ is due to $\Gamma_\R(n) \in \pi^{-\lfloor \tfrac{n}{2} \rfloor}\cdot \Q^\times$.
It is well-known that 
\[
{\rm vol}(\F^\times\backslash\A_\F^1) = \frac{2^d\cdot{\rm Reg}_\F\cdot h_\F}{w_\F},
\]
where $h_\F$ is the class number of $\F$ and $w_\F$ is the number of roots of unity contained in $\F$. In particular, we see that ${\rm vol}(\F^\times\backslash\A_\F^1) \in {\rm Reg}_\F\cdot\Q^\times$. 
Let $v$ be a finite place of $\F$ lying over a rational prime $p$. Recall the $\sigma$-linear isomorphism 
\[
t_\sigma : \mathcal{W}(\itPi_v,\psi_{n,\E_v}) \longrightarrow \mathcal{W}({}^\sigma\!\itPi_v,\psi_{n,\E_v})
\]
is defined by
\[
t_\sigma W(g) = \sigma \left( W({\rm diag}(u_{\sigma,p}^{-n+1},u_{\sigma,p}^{-n+2},\cdots,1)g)\right).
\]
Here $u_{\sigma,p} \in \Z_p^\times\subset \E_v^\times$ is the unique element such that $\sigma(\psi_{\Q_p}(x)) = \psi_{\Q_p}(u_{\sigma,p}x)$ for $x \in \Q_p$.
For $W \in \mathcal{W}(\itPi_v,\psi_{n,\E_v})$, we have
\begin{align*}
&\sigma(Z_v^{(n)}(W;\alpha_v))\\
&= \int_{N_{n-1}(\F_v)\backslash\GL_{n-1}(\F_v)}\sigma\left(W\left({\rm diag}(1,\alpha_v,\cdots,\alpha_v^{n-1})\bp g_v & 0 \\ 0 & 1\ep\right)\right)\omega_{\E_v/\F_v}^{n-1}(\det(g_v)){}^\sigma\!\chi_v(\det(g_v))\,dg_v\\
& = (\omega_{\E_v/\F_v}^{n-1}\cdot{}^\sigma\!\chi_v)^{-n(n-1)/2}(u_{\sigma,p})\cdot Z_v^{(n)}(t_\sigma W;\alpha_v).
\end{align*}
Here the first equality follows from \cite[Proposition A]{Grobner2018} and the asymptotics of Whittaker functions (cf.\,\cite[Proposition 2.2]{JPSS1979}). 
Similarly, by \cite[Theorem 5.3]{Matringe2011}, we have
\[
\sigma(L(1,\itPi_v,{\rm As}^{(-1)^{n-1}}\otimes\chi_v)) = L(1,{}^\sigma\!\itPi_v,{\rm As}^{(-1)^{n-1}}\otimes{}^\sigma\!\chi_v).
\]
Together with (\ref{E:Galois Gauss sum}), we conclude that
\begin{align}\label{E:proof A 2}
\begin{split}
&\sigma\left(G(\omega_{\E/\F}^{n-1}\cdot \chi)^{n(n-1)/2}\cdot \prod_{v \nmid \infty}\zeta_{\F_v}(n)\cdot\frac{Z_v^{(n)}(W_v;\alpha_v)}{L(1,\itPi_v,{\rm As}^{(-1)^{n-1}}\otimes\chi_v)}\right) \\
& = G(\omega_{\E/\F}^{n-1}\cdot {}^\sigma\!\chi)^{n(n-1)/2}\cdot \prod_{v \nmid \infty}\zeta_{\F_v}(n)\cdot\frac{Z_v^{(n)}(t_\sigma W_v;\alpha_v)}{L(1,{}^\sigma\!\itPi_v,{\rm As}^{(-1)^{n-1}}\otimes{}^\sigma\!\chi_v)}.
\end{split}
\end{align}
For $v \in S_\infty$ and $W \in \mathcal{W}(\itPi_v,\psi_{n,\C})$, we have
\begin{align*}
&Z_v^{(n)}(W;\iota_v(\alpha))\\
& = \int_{N_{n-1}(\R)\backslash \GL_{n-1}(\R)}W\left({\rm diag}(1,\iota_v(\alpha),\cdots,\iota_v(\alpha)^{n-1})\bp g & 0 \\ 0 & 1\ep\right){\rm sgn}^{n-1+{\sf w}}(\det(g))|\det(g)|^{\sf w}\,dg\\
& = |\iota_v(\alpha)|^{n(n-1)(n-2)/6}\cdot \left({\rm sgn}^{n-1+{\sf w}}(-\sqrt{-1}\cdot \iota_v(\alpha))|\iota_v(\alpha)|^{\sf w}\right)^{-n(n-1)/2}\cdot Z_v^{(n)}(W).
\end{align*}
Here the factor $|\iota_v(\alpha)|^{n(n-1)(n-2)/6}$ is due to the measure comparison under the change of variable from $g$ to ${\rm diag}((\sqrt{-1})^{1-n}\cdot \iota_v(\alpha)^{n-1},(\sqrt{-1})^{2-n}\cdot \iota_v(\alpha)^{n-2},\cdots,1)\cdot g$.
Therefore, we have
\begin{align}\label{E:proof A 3}
\begin{split}
&L_\infty^{(n)}([\itPi_\infty]_t,(\iota_v(\alpha))_{v \in S_\infty}) \\
&=\prod_{v \in S_\infty}|\iota_v(\alpha)|^{n(n-1)(n-2)/6}\cdot \left({\rm sgn}^{n-1+{\sf w}}(-\sqrt{-1}\cdot \iota_v(\alpha))|\iota_v(\alpha)|^{\sf w}\right)^{-n(n-1)/2} \cdot L_\infty^{(n)}([\itPi_\infty]_t).\\
\end{split}
\end{align}
The same factor appears when we compare $L_{\sigma,\infty}^{(n)}([{}^\sigma\!\itPi_\infty]_t,(\iota_v(\alpha))_{v \in S_\infty})$ with $L_{\sigma,\infty}^{(n)}([{}^\sigma\!\itPi_\infty]_t)$.
Note that 
\[
\prod_{v \in S_\infty}{\rm sgn}^{n-1+{\sf w}}(-\sqrt{-1}\cdot \iota_v(\alpha))^{-n(n-1)/2} \in \{\pm1\}.
\]
By (\ref{E:proof A 1})-(\ref{E:proof A 3}) and the non-vanishing of $Z_v^{(n)}(\mbox{ };\alpha_v)$ for all finite places $v$, to prove the theorem, it remains to show that
\begin{align*}
\sigma\left( \frac{\prod_{v \in S_\infty}|\iota_v(\alpha)|}{|D_\E|^{1/2}}\right) = \frac{\prod_{v \in S_\infty}|\iota_v(\alpha)|}{|D_\E|^{1/2}}.
\end{align*}
This equality is established in the proof of Lemma \ref{L:Gauss sum}.
This completes the proof.
\end{proof}

\begin{rmk}
The goal of \cite{Raphael2018} is to complete the local theory of Asai $L$-factors and $\varepsilon$-factors.
Since we concern only the behavior at $s=1$, we could avoid using \cite{Raphael2018} provided we can prove that $Z_v^{(n)}$ is non-vanishing for all $v \in S_\infty$.
We shall prove the cases $n=2,3$ by direct computations in Lemmas \ref{L:Asai 1} and \ref{L:Asai 2} below.
\end{rmk}

In the following theorem, we prove the algebraicity of the critical value 
$
L^{(\infty)}(1,\itPi,{\rm As}^{(-1)^n} \otimes \chi)
$
in terms of the bottom degree Whittaker period of $\itPi$.
The result is a generalization of the result of Grobner, Harris, and Lapid \cite[Theorem 7.1]{GHL2016} to essentially conjugate self-dual representations.

\begin{thm}\label{T:GHL 2}
We have
\begin{align*}
&\sigma\left(\frac{L^{(\infty)}(1,\itPi,{\rm As}^{(-1)^n}\otimes\chi)}
{|D_\F|^{\tfrac{n(n+1)}{4}}\cdot |D_\E|^{\mbox{$\tfrac{n(n^2+5)}{12}$}+\tfrac{{\sf w}\cdot n(n-1)}{4}}\cdot G(\omega_{\E/\F}^{n-1}\cdot \chi)^{\tfrac{n(n-1)}{2}}\cdot \pi^{d\lceil \tfrac{n}{2} \rceil}\cdot p^b(\itPi)}
\cdot \frac{L_{\infty}^{(n)}([\itPi_\infty^\vee]_t)}{B([\itPi_\infty]_b,[\itPi_\infty^\vee]_t)_\infty^{(n)}}\right)\\
&=\frac{L^{(\infty)}(1,{}^\sigma\!\itPi,{\rm As}^{(-1)^n}\otimes{}^\sigma\!\chi)}
{|D_\F|^{\tfrac{n(n+1)}{4}}\cdot |D_\E|^{\mbox{$\tfrac{n(n^2+5)}{12}$}+\tfrac{{\sf w}\cdot n(n-1)}{4}}\cdot G(\omega_{\E/\F}^{n-1}\cdot {}^\sigma\!\chi)^{\tfrac{n(n-1)}{2}}\cdot \pi^{d\lceil \tfrac{n}{2} \rceil}\cdot p^b({}^\sigma\!\itPi)}
\cdot \frac{L_{\sigma,\infty}^{(n)}([{}^\sigma\!\itPi_\infty^\vee]_t)}{B([{}^\sigma\!\itPi_\infty]_b,[{}^\sigma\!\itPi_\infty^\vee]_t)_{\sigma,\infty}^{(n)}}
\end{align*}
for all $\sigma \in {\rm Aut}(\C)$.
\end{thm}

\begin{proof}
The essentially conjugate self-dual condition implies that
\begin{align*}
L(s,\itPi \times \itPi^\vee) &= L(s,\itPi,{\rm As}^+ \otimes \chi) \cdot L(s,\itPi,{\rm As}^-\otimes \chi).
\end{align*}
Note that
\[
L(s,\itPi,{\rm As}^\pm \otimes \chi) = L(s,\itPi^\vee,{\rm As}^\pm \otimes \chi^{-1}).
\]
Also note that $L(s,\itPi,{\rm As}^{(-1)^{n-1}}\otimes \chi)$ has a simple pole at $s=1$.
Indeed, by (\ref{E:sign normalization 1}) and \cite[Lemma 4.1.4]{CHT2008}, there exists an algebraic Hecke character $\tilde{\chi}$ of $\A_\E^\times$ such that 
\[
\tilde{\chi} \vert_{\A_\F^\times} = \chi.
\]
This implies that $\itPi \otimes \tilde{\chi}$ is cohomological and conjugate self-dual.
Therefore, by \cite[Corollary 2.5.9]{Mok2015}, $L(s,\itPi\otimes\tilde{\chi},{\rm As}^{(-1)^{n-1}})$ has a pole at $s=1$.
Moreover, this pole must be simple as we explained in the proof of Theorem \ref{T:GHL}.
The assertions then follows immediately from Theorems \ref{T:adjoint}, \ref{T:Asai}, and the fact that (cf.\,\cite[Proposition 4.16]{Washingtonbook})
\[
\frac{{\rm Reg}_\E}{{\rm Reg}_\F} \in \Q^\times.
\]
This completes the proof.
\end{proof}

\subsection{Algebraicity of the Asai $L$-functions for $\GL_2$ and $\GL_3$}
In this section, we refine Theorem \ref{T:GHL} for $n=2$ and $n=3$ in Theorem \ref{T:Asai} below by explicitly determine the archimedean factors.
We keep the notation of \S\,\ref{S:RS}.
Assume $\itSigma$ and $\itPi$ are essentially conjugate self-dual with
\[
\itPi^\vee = \itPi^c \otimes \chi\circ {\rm N}_{\K/\F},\quad \itSigma^\vee = \itSigma^c \otimes \eta\circ {\rm N}_{\K/\F}
\]
for some algebraic Hecke characters $\chi$ and $\eta$ of $\A_\F^\times$ with parallel signatures ${\rm sgn}(\chi)$ and ${\rm sgn}(\eta)$, respectively.
After replacing $\chi$ by $\omega_{\E/\F}\cdot\chi$ and $\eta$ by $\omega_{\E/\F}\cdot\eta$ if necessary, we may assume that
\begin{align*}
{\rm sgn}(\chi) = (-1)^{{\sf w}(\itPi)},\quad {\rm sgn}(\eta) = (-1)^{{\sf w}(\itSigma)}.
\end{align*}
In the following lemmas, we compute the archimedean local zeta integrals appearing in the archimedean factor $L_\infty^{(2)}([\itPi_\infty]_t)$ and $L_\infty^{(3)}([\itSigma_\infty]_t)$.

\begin{lemma}\label{L:Asai 1}
Let 
\[
k_\infty = \displaystyle \frac{1}{\sqrt{2}}\bp 1 & \sqrt{-1} \\ \sqrt{-1} & 1\ep \in {\rm SU}(2).
\]
Let $v \in S_\infty$. 
\begin{itemize}
\item[(1)] For $i \neq \tfrac{\kappa_{1,v}-\kappa_{2,v}}{2}$, we have $Z_v^{(2)}(\rho(k_\infty)W_{
(\underline{\kappa}_v,{\sf w}(\itPi);\,i)})=0$.
\item[(2)] We have
\begin{align*}
Z_v^{(2)}\left(\rho(k_\infty)W_{\left(\underline{\kappa}_v, {\sf w}(\itPi);\,(\kappa_{1,v}-\kappa_{2,v})/2\right)}\right) 
&= (-1)^{(\kappa_{1,v}-\kappa_{2,v})/2}\cdot Z_v^{c,(2)}\left(\rho(\overline{k}_\infty)W_{\left(\underline{\kappa}_v^\vee, {\sf w}(\itPi);\,(\kappa_{1,v}-\kappa_{2,v})/2\right)}\right) \\
&= 2^{1-(\kappa_{1,v}-\kappa_{2,v})/2}(\sqrt{-1})^{(3\kappa_{1,v}+\kappa_{2,v})/2}\cdot \frac{\Gamma_\C(\tfrac{1}{2})^2}{\Gamma_\C(1)}\cdot L(1,\itPi_v,{\rm As}^-\otimes\chi_v).
\end{align*}
\end{itemize}
\end{lemma}

\begin{proof}
We drop the subscript $v$ for brevity.
Note that
\[
k_\infty^{-1}\bp \cos\theta & \sin\theta \\ -\sin\theta & \cos\theta\ep k_\infty = {\rm diag}( e^{\sqrt{-1}\,\theta},e^{-\sqrt{-1}\,\theta}).
\]
We see that $Z^{(2)}(\rho(k_\infty)W_{(\underline{\kappa};\,i)}) =0 $ unless $i=\tfrac{\kappa_1-\kappa_2}{2}$.
Also note that
\[
k_\infty^{-1}{\rm diag}(-1,1)k_\infty = \bp 0 & -\sqrt{-1} \\ \sqrt{-1} & 0 \ep
\]
and 
\[
\rho\left(\bp 0 & -\sqrt{-1} \\ \sqrt{-1} & 0 \ep\right)W_{\left(\underline{\kappa},{\sf w}(\itPi);\,(\kappa_1-\kappa_2)/2\right)} = (-1)^{\kappa_2}\cdot W_{\left(\underline{\kappa},{\sf w}(\itPi);\,(\kappa_1-\kappa_2)/2\right)} = (-1)^{1+{\sf w}(\itPi)}\cdot W_{\left(\underline{\kappa},{\sf w}(\itPi);\,(\kappa_1-\kappa_2)/2\right)}.
\]
Hence
\[
W_{\left(\underline{\kappa},{\sf w}(\itPi);\,(\kappa_1-\kappa_2)/2\right)}({\rm diag}(-a,\sqrt{-1})k_\infty) = (-1)^{1+{\sf w}(\itPi)}\cdot W_{\left(\underline{\kappa},{\sf w}(\itPi);\,(\kappa_1-\kappa_2)/2\right)}({\rm diag}(a,\sqrt{-1})k_\infty).
\]
Therefore, we have
\begin{align*}
&Z^{(2)}\left(\rho(k_\infty)W_{\left(\underline{\kappa},{\sf w}(\itPi);\,(\kappa_1-\kappa_2)/2\right)}\right) \\
&= \int_{\R^\times} W_{\left(\underline{\kappa},{\sf w}(\itPi);\,(\kappa_1-\kappa_2)/2\right)}\left( {\rm diag}(a,\sqrt{-1})k_\infty \right){\rm sgn}^{1+{\sf w}(\itPi)}(a)|a|_\R^{-{\sf w}(\itPi)}\,d^\times a\\
& = 2\int_{0}^{\infty} W_{\left(\underline{\kappa},{\sf w}(\itPi);\,(\kappa_1-\kappa_2)/2\right)}\left( {\rm diag}(a,\sqrt{-1})k_\infty \right)a^{-{\sf w}(\itPi)}\,d^\times a\\
& = 2^{1-(\kappa_1-\kappa_2)/2}(\sqrt{-1})^{(\kappa_1+\kappa_2)/2} \sum_{j=0}^{\tfrac{\kappa_1-\kappa_2}{2}}(-1)^j{\tfrac{\kappa_1-\kappa_2}{2} \choose j}\int_{0}^{\infty} W_{(\underline{\kappa},{\sf w}(\itPi);\,2j)}\left( {\rm diag}(a,1)\right)a^{-{\sf w}(\itPi)}\,d^\times a\\
& = 2^{1-(\kappa_1-\kappa_2)/2}(\sqrt{-1})^{(3\kappa_1+\kappa_2)/2} \sum_{j=0}^{\tfrac{\kappa_1-\kappa_2}{2}}{\tfrac{\kappa_1-\kappa_2}{2} \choose j}\int_0^\infty\int_{L}\frac{ds}{2\pi\sqrt{-1}}\,a^{-s+1}\Gamma_\C(\tfrac{s}{2}+j)\Gamma_\C(\tfrac{s+\kappa_1-\kappa_2}{2}-j) \,d^\times a\\
& = 2^{1-(\kappa_1-\kappa_2)/2}(\sqrt{-1})^{(3\kappa_1+\kappa_2)/2} \sum_{j=0}^{\tfrac{\kappa_1-\kappa_2}{2}}{\tfrac{\kappa_1-\kappa_2}{2} \choose j} \Gamma_\C(\tfrac{1}{2}+j)\Gamma_\C(\tfrac{1+\kappa_1-\kappa_2}{2}-j)\\
& = 2^{1-(\kappa_1-\kappa_2)/2}(\sqrt{-1})^{(3\kappa_1+\kappa_2)/2} \cdot \frac{\Gamma_\C(\tfrac{1}{2})^2\Gamma_\C(1+\tfrac{\kappa_1-\kappa_2}{2})}{\Gamma_\C(1)}.
\end{align*}
Here the last equality follows from \cite[Lemma 12.7-(i)]{Ichino2005}.
Note that 
\[
L(s,\itPi_v,{\rm As}^-\otimes\chi_v) = \Gamma_\R(s)^2\Gamma_\C(s+\tfrac{\kappa_1-\kappa_2}{2}).
\]
A similar computation show that
\begin{align*}
Z^{c,(2)}\left(\rho(\overline{k}_\infty)W_{\left(\underline{\kappa}^\vee, {\sf w}(\itPi);\,(\kappa_1-\kappa_2)/2\right)}\right) = 2^{1-(\kappa_1-\kappa_2)/2}(\sqrt{-1})^{-(3\kappa_1+\kappa_2)/2} \cdot \frac{\Gamma_\C(\tfrac{1}{2})^2\Gamma_\C(1+\tfrac{\kappa_1-\kappa_2}{2})}{\Gamma_\C(1)}.
\end{align*}
This completes the proof.
\end{proof}

\begin{lemma}\label{L:Asai 2}
Let 
\[
k_\infty' = \displaystyle {\rm diag}(\tfrac{1}{\sqrt{2}},\tfrac{1}{\sqrt{2}},1)\bp 1 & \sqrt{-1} & 0 \\ \sqrt{-1} & 1 & 0 \\ 0&0&1\ep \in {\rm SU}(3).
\]
Let $v \in S_\infty$ and $\underline{\ell}_v\succ\underline{\kappa}_v'$.
\begin{itemize}
\item[(1)] We have $Z_v^{(3)}\left(\rho(k_\infty')W_{\left(\underline{\ell}_v,{\sf w}(\itSigma);\,\underline{\kappa}_v',i\right)}\right) =0$
unless $\kappa_{1,v}' \equiv \kappa_{2,v}' \equiv {\sf w}(\itSigma)\,({\rm mod}\,2)$ and $i = \tfrac{\kappa_{1,v}'-\kappa_{2,v}'}{2}$.
\item[(2)] Suppose $\kappa_{1,v}' \equiv \kappa_{2,v}' \equiv {\sf w}(\itSigma)\,({\rm mod}\,2)$, then we have
\begin{align*}
Z_v^{(3)}\left(\rho(k_\infty')W_{\left(\underline{\ell}_v,{\sf w}(\itSigma);\,\underline{\kappa}_v',(\kappa_{1,v}'-\kappa_{2,v}')/2\right)}\right)
& = (-1)^{(\kappa_{1,v}'-\kappa_{2,v}')/2}\cdot Z_v^{c,(3)}\left(\rho(\overline{k}_\infty')W_{\left(\underline{\ell}_v^\vee,{\sf w}(\itSigma);\,(\underline{\kappa}_v')^\vee,(\kappa_{1,v}'-\kappa_{2,v}')/2\right)}\right)\\
& \in (\sqrt{-1})^{(\kappa_{1,v}'-\kappa_{2,v}')/2}\cdot \pi^{-2-\ell_{1,v}+\ell_{3,v}}\cdot\Q^\times.
\end{align*}
\end{itemize}
\end{lemma}

\begin{proof}
We drop the subscript $v$ for brevity.
Note that
\[
(k_\infty')^{-1}\bp \cos\theta & \sin\theta & 0  \\ -\sin\theta & \cos\theta  & 0 \\ 0 &0& 1 \ep k_\infty' = {\rm diag}( e^{\sqrt{-1}\,\theta}, e^{-\sqrt{-1}\,\theta},1).
\]
We see that $Z^{(3)}(\rho(k_\infty')W_{(\underline{\ell},{\sf w}(\itSigma);\,\underline{\kappa}',i)})=0$ unless $\kappa_1'-\kappa_2'$ is even and $i=\tfrac{\kappa_1'-\kappa_2'}{2}$.
For $\varepsilon = \pm 1$, we have
\[
(k_\infty')^{-1}{\rm diag}(-\varepsilon,\varepsilon)k_\infty' = \bp 0 & -\sqrt{-1}\,\varepsilon & 0  \\ \sqrt{-1}\,\varepsilon & 0  & 0 \\ 0& 0& 1\ep
\]
and 
\begin{align}\label{E:Asai 2 proof 3}
\begin{split}
\rho\left(\bp 0 & -\sqrt{-1}\,\varepsilon & 0  \\ \sqrt{-1}\,\varepsilon & 0  & 0 \\ 0& 0& 1\ep\right)W_{\left(\underline{\ell},{\sf w}(\itSigma);\,\underline{\kappa}',(\kappa_1'-\kappa_2')/2\right)}& = (-1)^{\kappa_2'}\cdot W_{\left(\underline{\ell},{\sf w}(\itSigma);\,\underline{\kappa}',(\kappa_1'-\kappa_2')/2\right)}\\
&= (-1)^{{\sf w}(\itSigma)}\cdot W_{\left(\underline{\ell},{\sf w}(\itSigma);\,\underline{\kappa}',(\kappa_1'-\kappa_2')/2\right)}.
\end{split}
\end{align}
Hence $Z^{(3)}\left(\rho(k_\infty')W_{\left(\underline{\ell},{\sf w}(\itSigma);\,\underline{\kappa}',(\kappa_1'-\kappa_2')/2\right)}\right) =0$ unless $\kappa_1' \equiv \kappa_{2}' \equiv {\sf w}(\itSigma)\,({\rm mod}\,2)$. 
Now we prove the second assertion. Suppose $\kappa_{1,v}' \equiv \kappa_{2,v}' \equiv {\sf w}(\itSigma)\,({\rm mod}\,2)$.
By (\ref{E:Asai 2 proof 3}), we have
\begin{align*}
& W_{\left(\underline{\ell},{\sf w}(\itSigma);\,\underline{\kappa}',(\kappa_1'-\kappa_2')/2\right)}({\rm diag}(-\varepsilon a_1a_2,\sqrt{-1}\,\varepsilon a_2,-1)k_\infty')\\
& = (-1)^{{\sf w}(\itSigma)}\cdot W_{\left(\underline{\ell},{\sf w}(\itSigma);\,\underline{\kappa}',(\kappa_1'-\kappa_2')/2\right)}({\rm diag}(a_1a_2,\sqrt{-1}\, a_2,-1)k_\infty').
\end{align*}
Therefore, we have
\begin{align*}
&Z^{(3)}\left(\rho(k_\infty')W_{\left(\underline{\ell},{\sf w}(\itSigma);\,\underline{\kappa}',(\kappa_1'-\kappa_2')/2\right)}\right) \\
&=\int_{\R^\times} \frac{d^\times a_1}{|a_1|_\R} \int_{\R^\times} d^\times a_2\,W_{\left(\underline{\ell},{\sf w}(\itSigma);\,\underline{\kappa}',(\kappa_1'-\kappa_2')/2\right)}({\rm diag}(a_1a_2,\sqrt{-1}\, a_2,-1)k_\infty'){\rm sgn}^{{\sf w}(\itSigma)}(a_1)|a_1a_2^2|_\R^{-{\sf w}(\itSigma)}\\
&=4\int_0^\infty \frac{d^\times a_1}{a_1} \int_0^\infty d^\times a_2\,W_{\left(\underline{\ell},{\sf w}(\itSigma);\,\underline{\kappa}',(\kappa_1'-\kappa_2')/2\right)}({\rm diag}(a_1a_2,\sqrt{-1}\, a_2,-1)k_\infty')(a_1a_2^2)^{-{\sf w}(\itSigma)}\\
& = 2^{2-(\kappa_1'-\kappa_2')/2}(\sqrt{-1})^{(\kappa_1'+\kappa_2')/2}\\
&\times\sum_{j=0}^{\tfrac{\kappa_1'-\kappa_2'}{2}}(-1)^{\ell_3+j}{\tfrac{\kappa_1'-\kappa_2'}{2} \choose j}\int_0^\infty \frac{d^\times a_1}{a_1} \int_0^\infty d^\times a_2\,W_{\left(\underline{\ell},{\sf w}(\itSigma);\,\underline{\kappa}',2j\right)}({\rm diag}(a_1a_2,a_2,1))(a_1a_2^2)^{-{\sf w}(\itSigma)}. 
\end{align*}
By (\ref{E:U(2) to U(3)}), we have
\begin{align*}
{\bf v}_{(\underline{\ell};\,\underline{\kappa}',i)} & = {\kappa_1'-\kappa_2' \choose i}^{-1}\sum_{m= \max\{0,\kappa_2'-\ell_2+i\}}^{\min\{i,\kappa_1'-\ell_2\}}{\kappa_1'-\ell_2 \choose m}{\ell_2-\kappa_2' \choose i-m}\cdot {\bf v}_{(\underline{\ell};\, \kappa_1'-\ell_2-m, m,\ell_1-\kappa_1',i-m, \ell_2-\kappa_2'-i+m, \kappa_2'-\ell_3)}
\end{align*}
for $0 \leq i \leq \kappa_1'-\kappa_2'$.
For $0 \leq j \leq \tfrac{\kappa_1'-\kappa_2'}{2}$ and $ \max\{0,\kappa_2'-\ell_2+2j\}\leq m \leq \min\{2j,\kappa_1'-\ell_2\}$, we have
\begin{align*}
& \int_0^\infty \frac{d^\times a_1}{a_1} \int_0^\infty d^\times a_2\,W_{(\underline{\ell},{\sf w}(\itSigma);\, \kappa_1'-\ell_2-m, m,\ell_1-\kappa_1',2j-m, \ell_2-\kappa_2'-2j+m, \kappa_2'-\ell_3)}({\rm diag}(a_1a_2,a_2,1))(a_1a_2^2)^{-{\sf w}(\itSigma)}\\
& = (\sqrt{-1})^{2\kappa_1'+\kappa_2'-2j}\\
&\times\int_{L_1}\frac{ds_1}{2\pi\sqrt{-1}}\int_{L_2}\frac{ds_2}{2\pi\sqrt{-1}}\,a_1^{-s_1+1}a_2^{-s_2+2} \Gamma_\C(\tfrac{s_1+\ell_1-\kappa_1'+2j}{2})\Gamma_\C(\tfrac{s_1+\kappa_1'-\ell_2+2j-2m}{2})\Gamma_\C(\tfrac{s_1+\kappa_1'-\ell_3-2j}{2})\\
&\quad\times \Gamma_\C(\tfrac{s_2+\kappa_1'-\ell_2+\kappa_2'-\ell_3}{2}) \Gamma_\C(\tfrac{s_2+\ell_1-\kappa_1'+\kappa_2'-\ell_3}{2}) \Gamma_\C(\tfrac{s_2+\ell_1-\kappa_1'+\ell_2-\kappa_2'}{2})\Gamma_\C(\tfrac{s_1+s_2+\ell_1-\ell_3+\kappa_2'-\ell_2+2j-2m}{2})^{-1}\\
&= (\sqrt{-1})^{2\kappa_1'+\kappa_2'-2j}\cdot\Gamma_\C(\tfrac{1+\ell_1-\kappa_1'+2j}{2})\Gamma_\C(\tfrac{1+\kappa_1'-\ell_2+2j-2m}{2})\Gamma_\C(\tfrac{1+\kappa_1'-\ell_3-2j}{2})\Gamma_\C(\tfrac{3+\ell_1-\ell_3+\kappa_2'-\ell_2+2j-2m}{2})^{-1}\\
&\times \Gamma_\C(\tfrac{2+\kappa_1'-\ell_2+\kappa_2'-\ell_3}{2}) \Gamma_\C(\tfrac{2+\ell_1-\kappa_1'+\kappa_2'-\ell_3}{2}) \Gamma_\C(\tfrac{2+\ell_1-\kappa_1'+\ell_2-\kappa_2'}{2}).
\end{align*}
Here the last equality follows from the Mellin inversion formula. Since $\Gamma(\tfrac{1}{2}+n) \in \pi^{1/2}\cdot\Q_{>0}$ for $n \in \Z_{\geq 0}$, we see that the product of the gamma factors above belongs to $\pi^{-2-\ell_1+\ell_3}\cdot \Q_{>0}$.
Therefore, we have
\begin{align*}
Z^{(3)}\left(\rho(k_\infty')W_{\left(\underline{\ell},{\sf w}(\itSigma);\,\underline{\kappa}',(\kappa_{1}'-\kappa_{2}')/2\right)}\right) &\in (-1)^{\ell_3}\cdot(\sqrt{-1})^{(5\kappa_1'+3\kappa_2')/2}\cdot \pi^{-2-\ell_1+\ell_3}\cdot\Q_{>0} \\
&\subset (\sqrt{-1})^{(\kappa_1'-\kappa_2')/2}\cdot \pi^{-2-\ell_1+\ell_3}\cdot\Q^\times.
\end{align*}
To complete the proof of the second assertion, it remains to show that
\begin{align}\label{E:Asai 2 proof 2}
&Z^{(3)}\left(\rho(k_\infty')W_{\left(\underline{\ell},{\sf w}(\itSigma);\,\underline{\kappa}',(\kappa_1'-\kappa_2')/2\right)}\right) = (-1)^{(\kappa_{1}'-\kappa_{2}')/2}\cdot Z^{c,(3)}\left(\rho(\overline{k}_\infty')W_{\left(\underline{\ell}^\vee,{\sf w}(\itSigma);\,(\underline{\kappa}')^\vee,(\kappa_1'-\kappa_2')/2\right)}\right).
\end{align}
First note that, when $\underline{\kappa}' = (\ell_2,\ell_2)$, we deduce from the above computation that
\begin{align}\label{E:Asai 2 proof 1}
\begin{split}
&Z^{(3)}\left(\rho(k_\infty')W_{\left(\underline{\ell},{\sf w}(\itSigma);\,(\ell_{2},\ell_{2}),0\right)}\right)\\
& = Z^{c,(3)}\left(\rho(\overline{k}_\infty')W_{\left(\underline{\ell}^\vee,{\sf w}(\itSigma);\,(-\ell_{2},-\ell_{2}),0\right)}\right)\\
&=(-1)^{\sf w(\itSigma)}\cdot2^2\cdot\Gamma_\C(\tfrac{1}{2})\Gamma_\C(\tfrac{1+\ell_1-\ell_2}{2})\Gamma_\C(\tfrac{1+\ell_2-\ell_3}{2})\Gamma_\C(\tfrac{3+\ell_1-\ell_3}{2})^{-1}\Gamma_\C(\tfrac{2+\ell_2-\ell_3}{2}) \Gamma_\C(\tfrac{2+\ell_1-\ell_3}{2}) \Gamma_\C(\tfrac{2+\ell_1-\ell_2}{2}).
\end{split}
\end{align}
Define 
\[
f_1 \in {\rm Hom}_{{\rm SO}(3)}(V^{(3)}_{\underline{\ell}},\C),\quad f_2 \in {\rm Hom}_{{\rm SO}(3)}(V^{(3)}_{\underline{\ell}^\vee},\C)
\]
by
\[
f_1({\bf v}_{(\underline{\ell};\,\underline{\kappa}',i)}) = Z^{(3)}(W_{(\underline{\ell},{\sf w}(\itSigma);\,\underline{\kappa}',i)}),\quad f_2({\bf v}_{(\underline{\ell}^\vee;\,(\underline{\kappa}')^\vee,i)}) = Z^{c,(3)}(W_{(\underline{\ell}^\vee,{\sf w}(\itSigma);\,(\underline{\kappa}')^\vee,i)})
\]
for $ \underline{\ell} \succ \underline{\kappa}'$ and $0 \leq i \leq \kappa_1'-\kappa_2'$.
Note that $f_2 \circ c_{\underline{\ell}}^{(3)} \in {\rm Hom}_{{\rm SO}(3)}(V^{(3)}_{\underline{\ell}},\C)$. 
By (\ref{E:Asai 2 proof 1}), we have
\[
f_1(\rho_{\underline{\ell}}(k_\infty')\cdot{\bf v}_{(\underline{\ell};\,(\ell_2,\ell_2),0)}) = f_2 \circ c_{\underline{\ell}}^{(3)}(\rho_{\underline{\ell}}(k_\infty')\cdot{\bf v}_{(\underline{\ell};\,(\ell_2,\ell_2),0)})\neq 0.
\]
Since ${\rm dim}_\C\,{\rm Hom}_{{\rm SO}(3)}(V^{(3)}_{\underline{\ell}},\C) \leq 1$, we conclude that $f_1 = f_2 \circ c_{\underline{\ell}}^{(3)}$ and (\ref{E:Asai 2 proof 2}) then follows from (\ref{E:conjugate equiv. 3}). This completes the proof.
\end{proof}

\begin{thm}\label{T:Asai}
\noindent
\begin{itemize}
\item[(1)] We have
\begin{align*}
&\sigma\left(\frac{{\rm Res}_{s=1}L^{(\infty)}(s,\itPi,{\rm As}^-\otimes\chi)}{|D_\F|^{1/2}\cdot (\sqrt{-1})^{d{\sf w}(\itPi)}\cdot G(\omega_{\E/\F})\cdot G(\chi)\cdot{\rm Reg}_\F\cdot\pi^{d+\sum_{v \in S_\infty}(\kappa_{1,v}-\kappa_{2,v})/2}\cdot  p^t(\itPi)}\right) \\
& = \frac{{\rm Res}_{s=1}L^{(\infty)}(s,{}^\sigma\!\itPi,{\rm As}^-\otimes{}^\sigma\!\chi)}{|D_\F|^{1/2}\cdot (\sqrt{-1})^{d{\sf w}(\itPi)}\cdot G(\omega_{\E/\F})\cdot G({}^\sigma\!\chi)\cdot{\rm Reg}_\F\cdot\pi^{d+\sum_{v \in S_\infty}(\kappa_{1,v}-\kappa_{2,v})/2}\cdot  p^t({}^\sigma\!\itPi)}
\end{align*}
for all $\sigma \in {\rm Aut}(\C)$.
\item[(2)] We have
\begin{align*}
&\sigma\left(\frac{{\rm Res}_{s=1}L^{(\infty)}(s,\itSigma,{\rm As}^+\otimes\eta)}{|D_\E|^{1/2}\cdot(\sqrt{-1})^{d{\sf w}(\itSigma)}\cdot G(\eta)^3\cdot{\rm Reg}_\F\cdot\pi^{3d+\sum_{v \in S_\infty}(\ell_{1,v}-\ell_{3,v})}\cdot p^t(\itSigma)} \right)\\
& = \frac{{\rm Res}_{s=1}L^{(\infty)}(s,{}^\sigma\!\itSigma,{\rm As}^+\otimes{}^\sigma\!\eta)}{|D_\E|^{1/2}\cdot(\sqrt{-1})^{d{\sf w}(\itSigma)}\cdot  G({}^\sigma\!\eta)^3\cdot{\rm Reg}_\F\cdot\pi^{3d+\sum_{v \in S_\infty}(\ell_{1,v}-\ell_{3,v})}\cdot p^t({}^\sigma\!\itSigma)}
\end{align*}
for all $\sigma \in {\rm Aut}(\C)$.
\end{itemize}
\end{thm}
\begin{proof}
By Lemma \ref{L:Gauss sum} and Theorem \ref{T:GHL}, to prove the first assertion, it suffices to show that
\begin{align}
L_v^{(2)}([\itPi_v]_t) &\in \pi^{-(\kappa_{1,v}-\kappa_{2,v})/2}\cdot {\Q},\label{E:Asai proof 1}\\
L_v^{c,(2)}([\itPi_v^c]_t) &= (-1)^{{\sf w}(\itPi)}\cdot L_v^{(2)}([\itPi_v]_t)\label{E:Asai proof 2}
\end{align}
for all $v \in S_\infty$. 
Indeed, by Lemma \ref{L:Gauss sum} and (\ref{E:Asai proof 2}), we have
\begin{align*}
L_{\sigma,\infty}^{(2)}([{}^\sigma\!\itPi_\infty]_t) &= \prod_{v \in S_\infty,\, \sigma^{-1}\circ\iota_v=\iota_{\sigma^{-1}\circ v}}L_{\sigma^{-1}\circ v}^{(2)}([\itPi_{\sigma^{-1}\circ v}]_t) \prod_{v \in S_\infty,\, \sigma^{-1}\circ\iota_v=\overline{\iota}_{\sigma^{-1}\circ v}}L_{\sigma^{-1}\circ v}^{c,(2)}([\itPi_{\sigma^{-1}\circ v}^c]_t)\\
& = \prod_{v \in S_\infty,\, \sigma\circ\iota_v=\iota_{\sigma\circ v}}L_{v}^{(2)}([\itPi_v]_t) \prod_{v \in S_\infty, \,\sigma\circ\iota_v=\overline{\iota}_{\sigma\circ v}}L_v^{c,(2)}([\itPi_v^c]_t)\\
& = (-1)^{{}^\sharp\{v \in S_\infty\,\vert\,\sigma\circ\iota_v=\overline{\iota}_{\sigma\circ v}\}\cdot {\sf w}(\itPi)}\cdot L_{\infty}^{(2)}([\itPi_\infty]_t)\\
& = \frac{\sigma(|D_\E|^{{\sf w}(\itPi)/2}\cdot(\sqrt{-1})^{d{\sf w}(\itPi)})}{|D_\E|^{{\sf w}(\itPi)/2}\cdot(\sqrt{-1})^{d{\sf w}(\itPi)}}\cdot L_{\infty}^{(2)}([\itPi_\infty]_t).
\end{align*}
The first assertion then follows from Theorem \ref{T:GHL} and (\ref{E:Asai proof 1}).
Fix $v \in S_\infty$. We drop the subscript $v$ for brevity. 
Recall that
\[
[\itPi_v]_t  = \sum_{i=0}^{\kappa_1-\kappa_2}(-1)^i{\kappa_1-\kappa_2 \choose i}\cdot \rho(k_\infty)W_{(\underline{\kappa}, {\sf w}(\itPi);\,i)}\otimes \xi_t^{(2)}(\rho_{\underline{\kappa}^\vee}(k_\infty)x^iy^{\kappa_1-\kappa_2-i}).
\]
By Lemma \ref{L:Asai 1}, we have
\begin{align*}
L^{(2)}([\itPi_v]_t) &= (-1)^{(\kappa_1-\kappa_2)/2}{\kappa_1-\kappa_2 \choose \frac{\kappa_1-\kappa_2}{2}}\cdot Z^{(2)}\left(\rho(k_\infty)W_{(\underline{\kappa},{\sf w}(\itPi);\,(\kappa_1-\kappa_2)/2)}\right)\\
&\times (s_\R^{(2)}\otimes f_{\lambda,\C})\left(\xi_t^{(2)}\left(\rho_{\underline{\kappa}^\vee}(k_\infty)(xy)^{(\kappa_1-\kappa_2)/2}\right)\right)\\
& = (-1)^{(\kappa_1-\kappa_2)/2}\cdot2^{1-(\kappa_1-\kappa_2)/2}\cdot{\kappa_1-\kappa_2 \choose \frac{\kappa_1-\kappa_2}{2}}\cdot(\sqrt{-1})^{(3\kappa_1+\kappa_2)/2} \cdot \frac{\Gamma_\C(\tfrac{1}{2})^2\Gamma_\C(1+\tfrac{\kappa_1-\kappa_2}{2})}{\Gamma_\C(1)}\\
&\times (s_\R^{(2)}\otimes f_{\lambda,\C})\left(\xi_t^{(2)}\left(\rho_{\underline{\kappa}^\vee}(k_\infty)(xy)^{(\kappa_1-\kappa_2)/2}\right)\right)\\
& = (-1)^{(\kappa_1+\kappa_2)/2}\cdot2^{1-\kappa_1+\kappa_2}\cdot{\kappa_1-\kappa_2 \choose \frac{\kappa_1-\kappa_2}{2}}\cdot \frac{\Gamma_\C(\tfrac{1}{2})^2\Gamma_\C(1+\tfrac{\kappa_1-\kappa_2}{2})}{\Gamma_\C(1)}\\
&\times (s_\R^{(2)}\otimes f_{\lambda,\C})\left(\xi_t^{(2)}\left((x^2+y^2)^{(\kappa_1-\kappa_2)/2}\right)\right).
\end{align*}
To prove that (\ref{E:Asai proof 1}) holds, it remains to show that
\[
(s_\R^{(2)}\otimes f_{\lambda,\C})\left(\xi_t^{(2)}\left((x^2+y^2)^{(\kappa_1-\kappa_2)/2}\right)\right) \in \Q.
\]
By (\ref{E:Lie algebra action 1}), we have
\begin{align*}
\xi_t^{(2)}\left((x^2+y^2)^{(\kappa_1-\kappa_2)/2}\right)\in  \wedge^2 (\frak{g}_{2,\C} / \frak{k}_{2,\C})^*_\Q \,\otimes M_{\lambda_{\iota}}^\vee\otimes M_{\lambda_{\overline{\iota}}}^\vee.
\end{align*}
It is clear that
\[
f_{\lambda,\C}({\bf u}_{\iota} \otimes {\bf u}_{\overline{\iota}}) \in \Q
\]
for all ${\bf u}_{\iota} \otimes {\bf u}_{\overline{\iota}} \in M_{\lambda_{\iota}}^\vee\otimes M_{\lambda_{\overline{\iota}}}^\vee$.
Thus we are reduced to show that
\[
s_\R^{(2)}\left(\wedge^2(\frak{g}_{2,\C} / \frak{k}_{2,\C})^*_\Q\right) = \Q.
\]
Indeed, the image is equal to a scalar multiple of $\Q$. To determine the scalar, note that the map $\frak{g}_{2,\C}' / \frak{k}_{2,\C}' \rightarrow \frak{g}_{2,\C} / \frak{k}_{2,\C}$ induced by $\delta_{(2),\R}$ is given by
\[
e_{11} \longmapsto \frac{1}{2}\cdot Y_{(0,0)},\quad e_{12} \longmapsto -\frac{1}{2}\cdot Y_{(1,-1)} +\frac{1}{2} \cdot Y_{(-1,1)}.
\]
Hence we have
\begin{align}\label{E:Asai proof 3}
\delta_{(2),\R}^* Y_{(0,0)}^* = \frac{1}{2}\cdot e_{11}^*,\quad \delta_{(2),\R}^* Y_{\pm(1,-1)}^* = \mp\frac{1}{2}\cdot e_{12}^*.
\end{align}
In particular, we have
\[
\delta_{(2),\R}^* Y_{(0,0)}^* \wedge \delta_{(2),\R}^* Y_{(1,-1)}^* = -\frac{1}{4} \cdot e_{11}^* \wedge e_{12}^*.
\]
In other words, $s_\R^{(2)}(Y_{(0,0)}^* \wedge Y_{(1,-1)}^*) = -\tfrac{1}{4}$.
Now we show that (\ref{E:Asai proof 2}) holds. 
Similarly, by Lemma \ref{L:Asai 1}, we have
\begin{align*}
L^{c,(2)}([\itPi_v^c]_t) &= (-1)^{(\kappa_1-\kappa_2)/2}{\kappa_1-\kappa_2 \choose \frac{\kappa_1-\kappa_2}{2}}\cdot Z^{c,(2)}\left(\rho(\overline{k}_\infty)W_{(\underline{\kappa}^\vee,{\sf w}(\itPi);\,(\kappa_1-\kappa_2)/2)}\right)\\
&\times (-1)^{{\sf w}(\itPi)}\cdot (s_\R^{(2)}\otimes f_{\lambda^c,\C})\left(\xi_t^{c,(2)}\left(\rho_{\underline{\kappa}}(\overline{k}_\infty)(xy)^{(\kappa_1-\kappa_2)/2}\right)\right)\\
& = (-1)^{(\kappa_1+\kappa_2)/2}\cdot2^{1-\kappa_1+\kappa_2}\cdot{\kappa_1-\kappa_2 \choose \frac{\kappa_1-\kappa_2}{2}}\cdot \frac{\Gamma_\C(\tfrac{1}{2})^2\Gamma_\C(1+\tfrac{\kappa_1-\kappa_2}{2})}{\Gamma_\C(1)}\\
&\times (-1)^{{\sf w}(\itPi)}\cdot (s_\R^{(2)}\otimes f_{\lambda^c,\C})\left(\xi_t^{c,(2)}\left((x^2+y^2)^{(\kappa_1-\kappa_2)/2}\right)\right).
\end{align*}
Therefore, to prove (\ref{E:Asai proof 2}), it suffices to show that
\[
(s_\R^{(2)}\otimes f_{\lambda,\C})\left(\xi_t^{(2)}(x^{\kappa_1-\kappa_2-2i}y^{2i})\right) = (s_\R^{(2)}\otimes f_{\lambda^c,\C})\left(\xi_t^{c,(2)}(x^{2i}y^{\kappa_1-\kappa_2-2i})\right)
\]
for all $0 \leq i \leq \tfrac{\kappa_1-\kappa_2}{2}$.
By (\ref{E:conjugate equiv. 1}) and the definition of $f_{\lambda^c,\C}$, we have
\begin{align*}
(s_\R^{(2)}\otimes f_{\lambda^c,\C})\left(\xi_t^{c,(2)}(x^{2i}y^{\kappa_1-\kappa_2-2i})\right) & = \left(s_\R^{(2)}\circ \wedge^2(c_{(1,-1)}^{(2)})^* \otimes f_{\lambda^c,\C}\circ {\bf c}_{\lambda^\vee}^{(2)}\right)\left(\xi_t^{(2)}(x^{\kappa_1-\kappa_2-2i}y^{2i})\right)\\
& = \left(s_\R^{(2)}\circ \wedge^2(c_{(1,-1)}^{(2)})^* \otimes f_{\lambda,\C}\right)\left(\xi_t^{(2)}(x^{\kappa_1-\kappa_2-2i}y^{2i})\right).
\end{align*}
Thus we need only to show that
\[
s_\R^{(2)}\circ \wedge^2(c_{(1,-1)}^{(2)})^* = s_\R^{(2)}.
\]
Indeed, the two functionals on $\wedge^2(\frak{g}_{2,\C} / \frak{k}_{2,\C})^*$ differ by a scalar. By (\ref{E:conjugate equiv. Lie 1}) and (\ref{E:Asai proof 3}), we have
\begin{align*}
\delta_{(2),\R}^* \circ(c_{(1,-1)}^{(2)})^*(Y_{(0,0)}^*) \wedge \delta_{(2),\R}^* \circ(c_{(1,-1)}^{(2)})^*(Y_{(1,-1)}^*)
& = -\delta_{(2),\R}^* Y_{(0,0)}^* \wedge \delta_{(2),\R}^* Y_{(-1,1)}^* \\
& = \delta_{(2),\R}^* Y_{(0,0)}^* \wedge \delta_{(2),\R}^* Y_{(1,-1)}^*.
\end{align*}
Thus the scalar between the functionals is equal to $1$.
This completes the proof of the first assertion.

By Lemma \ref{L:Gauss sum} and Theorem \ref{T:GHL}, to prove the second assertion, it suffices to show that
\begin{align}
L_v^{(3)}([\itSigma_v]_t) &\in \pi^{-2-\ell_{1,v}+\ell_{3,v}}\cdot {\Q},\label{E:Asai proof 4}\\
L_v^{c,(3)}([\itSigma_v^c]_t) &= (-1)^{{\sf w}(\itSigma)}\cdot L_v^{(3)}([\itSigma_v]_t)\label{E:Asai proof 5}
\end{align}
for all $v \in S_\infty$. Fix $v \in S_\infty$. We drop the subscript $v$ for brevity. 
Recall that 
\[
[\itSigma_v]_t  = \sum_{\underline{\ell} \succ \underline{\kappa}'}C_{(\underline{\ell};\,\underline{\kappa}')}\sum_{i=0}^{\kappa_1'-\kappa_2'}(-1)^i{\kappa_1'-\kappa_2' \choose i}\cdot \rho(k_\infty')W_{(\underline{\ell},{\sf w}(\itSigma);\,\underline{\kappa}',i)}\otimes \xi_t^{(3)}(\rho_{\underline{\ell}^\vee}(k_\infty'){\bf v}_{(\underline{\ell};\,\underline{\kappa}',\kappa_1'-\kappa_2'-i)}).
\]
By Lemma \ref{L:Asai 2}-(1), we have
\begin{align*}\label{E:Asai proof 6}
\begin{split}
L^{(3)}([\itSigma_v]_t) &=  \sum_{\underline{\ell} \succ \underline{\kappa}',\,\kappa_1' \equiv \kappa_2' \equiv {\sf w}(\itSigma)\,({\rm mod}\,2)}C_{(\underline{\ell};\,\underline{\kappa}')}(-1)^{(\kappa_1'-\kappa_2')/2}{\kappa_1'-\kappa_2' \choose \tfrac{\kappa_1'-\kappa_2'}{2}}\cdot Z^{(3)}\left(\rho(k_\infty')W_{\left(\underline{\ell},{\sf w}(\itSigma);\,\underline{\kappa}',(\kappa_1'-\kappa_2')/2\right)}\right)\\
&\quad\times (s_\R^{(3)}\otimes f_{\mu,\C})\left(\xi_t^{(3)}\left(\rho_{\underline{\ell}^\vee}(k_\infty'){\bf v}_{\left(\underline{\ell}^\vee;\,(\underline{\kappa}')^\vee,(\kappa_1'-\kappa_2')/2\right)}\right)\right)\\
&=  \sum_{\underline{\ell} \succ \underline{\kappa}',\,\kappa_1' \equiv \kappa_2' \equiv {\sf w}(\itSigma)\,({\rm mod}\,2)}C_{(\underline{\ell};\,\underline{\kappa}')}(-1)^{(\kappa_1'-\kappa_2')/2}{\kappa_1'-\kappa_2' \choose \tfrac{\kappa_1'-\kappa_2'}{2}}\cdot Z^{(3)}\left(\rho(k_\infty')W_{\left(\underline{\ell},{\sf w}(\itSigma);\,\underline{\kappa}',(\kappa_1'-\kappa_2')/2\right)}\right)\\
&\quad\times 2^{-(\kappa_1'-\kappa_2')/2}(\sqrt{-1})^{(\kappa_1'-\kappa_2')/2}\sum_{j=0}^{\tfrac{\kappa_1'-\kappa_2'}{2}}{\tfrac{\kappa_1'-\kappa_2'}{2} \choose j}(s_\R^{(3)}\otimes f_{\mu,\C})\left(\xi_t^{(3)}\left({\bf v}_{\left(\underline{\ell}^\vee;\,(\underline{\kappa}')^\vee,2j\right)}\right)\right).
\end{split}
\end{align*}
To prove (\ref{E:Asai proof 4}), by Lemma \ref{L:Asai 2}-(2), it suffices to show that
\[
(s_\R^{(3)}\otimes f_{\mu,\C})\left(\xi_t^{(3)}\left({\bf v}_{\left(\underline{\ell}^\vee;\,(\underline{\kappa}')^\vee,i\right)}\right)\right) \in \Q
\]
for all $\underline{\ell} \succ \underline{\kappa}'$ and $0 \leq i \leq \kappa_1'-\kappa_2'$.
By (\ref{E:U(2) to U(3)}), we have
\[
\xi_t^{(3)}\left({\bf v}_{\left(\underline{\ell}^\vee;\,(\underline{\kappa}')^\vee,i\right)}\right) \in \wedge^5 (\frak{g}_{3,\C} / \frak{k}_{3,\C})^*_\Q \,\otimes M_{\mu_{\iota}}^\vee\otimes M_{\mu_{\overline{\iota}}}^\vee.
\]
It is clear that
\[
f_{\mu,\C}({\bf v}_{\iota} \otimes {\bf v}_{\overline{\iota}}) \in \Q
\]
for all ${\bf v}_{\iota} \otimes {\bf v}_{\overline{\iota}} \in M_{\mu_{\iota}}^\vee\otimes M_{\mu_{\overline{\iota}}}^\vee$.
Thus we are reduced to show that
\[
s_\R^{(3)}\left(\wedge^5(\frak{g}_{3,\C} / \frak{k}_{3,\C})^*_\Q\right) = \Q.
\]
Indeed, the image is equal to a scalar multiple of $\Q$. To determine the scalar, note that the map $\frak{g}_{3,\C}' / \frak{k}_{3,\C}' \rightarrow \frak{g}_{3,\C} / \frak{k}_{3,\C}$ induced by $\delta_{(3),\R}$ is given by
\begin{gather*}
e_{11} \longmapsto \frac{2}{3}\cdot Z_{12}+ \frac{1}{3}\cdot Z_{23}, \quad e_{22} \longmapsto -\frac{1}{3}\cdot Z_{12}+ \frac{1}{3}\cdot Z_{23},\\
e_{12} \longmapsto -\frac{1}{2}\cdot X_{(1,-1,0)}+\frac{1}{2}\cdot X_{(-1,1,0)},\quad e_{13} \longmapsto -\frac{1}{2}\cdot X_{(1,0,-1)}+\frac{1}{2}\cdot X_{(-1,0,1)},\\
e_{23} \longmapsto -\frac{1}{2}\cdot X_{(0,1,-1)}+\frac{1}{2}\cdot X_{(0,-1,1)}.
\end{gather*}
Hence we have
\begin{align}\label{E:Asai proof 7}
\begin{split}
\delta_{(3),\R}^* Z_{12}^* &= \frac{2}{3}\cdot e_{11}^* - \frac{1}{3}\cdot e_{22}^*,\quad \delta_{(3),\R}^* Z_{23}^* = \frac{1}{3}\cdot e_{11}^* + \frac{1}{3}\cdot e_{22}^*,\\
\delta_{(3),\R}^* X_{\pm(1,-1,0)}^* &= \mp \frac{1}{2}\cdot e_{12}^*,\quad \delta_{(3),\R}^* X_{\pm(1,0,-1)}^* = \mp \frac{1}{2}\cdot e_{13}^*,\quad 
\delta_{(3),\R}^* X_{\pm(0,1,-1)}^* = \mp \frac{1}{2}\cdot e_{23}^*.\\
\end{split}
\end{align}
In particular, we have
\begin{align*}
& \delta_{(3),\R}^* Z_{12}^* \wedge \delta_{(3),\R}^* Z_{23}^* \wedge \delta_{(3),\R}^* X_{(1,-1,0)}^* \wedge \delta_{(3),\R}^* X_{(1,0,-1)}^* \wedge \delta_{(3),\R}^* X_{(0,1,-1)}^* \\
&= -\frac{1}{24}\cdot e_{11}^* \wedge e_{22}^* \wedge e_{12}^* \wedge e_{13}^* \wedge e_{23}^*.
\end{align*}
In other words, $s_\R^{(3)}(Z_{12}^* \wedge Z_{23}^* \wedge X_{(1,-1,0)}^* \wedge X_{(1,0,-1)}^* \wedge X_{(0,1,-1)}^*) = -\tfrac{1}{24}$.
Now we show that (\ref{E:Asai proof 5}) holds.
Similarly, by Lemma \ref{L:Asai 2}, we have
\begin{align*}
L^{c,(3)}([\itSigma_v^c]_t) &=  \sum_{\underline{\ell} \succ \underline{\kappa}',\,\kappa_1' \equiv \kappa_2' \equiv {\sf w}(\itSigma)\,({\rm mod}\,2)}C_{(\underline{\ell};\,\underline{\kappa}')}(-1)^{(\kappa_1'-\kappa_2')/2}{\kappa_1'-\kappa_2' \choose \tfrac{\kappa_1'-\kappa_2'}{2}}\cdot Z^{(3)}\left(\rho(\overline{k}_\infty')W_{\left(\underline{\ell}^\vee,{\sf w}(\itSigma);\,(\underline{\kappa}')^\vee,(\kappa_1'-\kappa_2')/2\right)}\right)\\
&\quad\times (-1)^{{\sf w}(\itSigma)}\cdot (s_\R^{(3)}\otimes f_{\mu^c,\C})\left(\xi_t^{c,(3)}\left(\rho_{\underline{\ell}}(\overline{k}_\infty'){\bf v}_{\left(\underline{\ell};\,\underline{\kappa}',(\kappa_1'-\kappa_2')/2\right)}\right)\right)\\
&=  \sum_{\underline{\ell} \succ \underline{\kappa}',\,\kappa_1' \equiv \kappa_2' \equiv {\sf w}(\itSigma)\,({\rm mod}\,2)}C_{(\underline{\ell};\,\underline{\kappa}')}(-1)^{(\kappa_1'-\kappa_2')/2}{\kappa_1'-\kappa_2' \choose \tfrac{\kappa_1'-\kappa_2'}{2}}\cdot Z^{(3)}\left(\rho(k_\infty')W_{\left(\underline{\ell},{\sf w}(\itSigma);\,\underline{\kappa}',(\kappa_1'-\kappa_2')/2\right)}\right)\\
&\quad\times 2^{-(\kappa_1'-\kappa_2')/2}(\sqrt{-1})^{(\kappa_1'-\kappa_2')/2}\sum_{j=0}^{\tfrac{\kappa_1'-\kappa_2'}{2}}{\tfrac{\kappa_1'-\kappa_2'}{2} \choose j}(-1)^{{\sf w}(\itSigma)}\cdot (s_\R^{(3)}\otimes f_{\mu^c,\C})\left(\xi_t^{c,(3)}\left({\bf v}_{\left(\underline{\ell};\,\underline{\kappa}',2j\right)}\right)\right).
\end{align*}
Therefore, to prove (\ref{E:Asai proof 5}), it suffices to show that
\[
(s_\R^{(3)} \otimes f_{\mu,\C})\left(\xi_t^{(3)}\left({\bf v}_{\left(\underline{\ell}^\vee;\,(\underline{\kappa}')^\vee,2j\right)}\right)\right) = (s_\R^{(3)} \otimes f_{\mu^c,\C})\left(\xi_t^{c,(3)}\left({\bf v}_{\left(\underline{\ell};\,\underline{\kappa}',\kappa_1'-\kappa_2'-2j\right)}\right)\right)
\]
for all $\underline{\ell} \succ \underline{\kappa}'$ with $\kappa_1' \equiv \kappa_2' \equiv {\sf w}(\itSigma)\,({\rm mod}\,2)$ and $0 \leq j \leq \tfrac{\kappa_1'-\kappa_2'}{2}$.
By (\ref{E:conjugate equiv. 3}) and the definition of $f_{\mu^c,\C}$, we have
\begin{align*}
&(s_\R^{(3)} \otimes f_{\mu^c,\C})\left(\xi_t^{c,(3)}\left({\bf v}_{\left(\underline{\ell};\,\underline{\kappa}',\kappa_1'-\kappa_2'-2j\right)}\right)\right)\\
& = (-1)^{\ell_2-\kappa_2'}\cdot \left(s_\R^{(3)}\circ \wedge^5(c_{(1,0,-1)}^{(3)})^* \otimes f_{\mu^c,\C} \circ {\bf c}_{\mu^\vee}^{(3)}\right)\left(\xi_t^{(3)}\left({\bf v}_{\left(\underline{\ell}^\vee;\,(\underline{\kappa}')^\vee,2j\right)}\right)\right)
\\
& = \left(s_\R^{(3)}\circ \wedge^5(c_{(1,0,-1)}^{(3)})^* \otimes f_{\mu,\C}\right)\left(\xi_t^{(3)}\left({\bf v}_{\left(\underline{\ell}^\vee;\,(\underline{\kappa}')^\vee,2j\right)}\right)\right)
\end{align*}
Thus we need only to show that
\[
s_\R^{(3)}\circ \wedge^5(c_{(1,0,-1)}^{(3)})^* = s_\R^{(3)}.
\]
Indeed, the two functionals on $\wedge^5(\frak{g}_{2,\C} / \frak{k}_{2,\C})^*$ differ by a scalar. 
By (\ref{E:conjugate equiv. Lie 2}), (\ref{E:conjugate equiv. Lie 3}), and (\ref{E:Asai proof 7}), we have
\begin{align*}
&\wedge^5 \delta_{(3),\R}^*\circ (c_{(1,0,-1)}^{(3)})^* \left( Z_{12}^* \wedge Z_{23}^* \wedge X_{(-1,1,0)}^*\wedge X_{(0,-1,1)}^* \wedge X_{(-1,0,1)}^*\right)\\
& = \wedge^5 \delta_{(3),\R}^*\left( Z_{12}^* \wedge Z_{23}^* \wedge (-X_{(1,-1,0)}^*)\wedge (-X_{(0,1,-1)}^*) \wedge (-X_{(1,0,-1)}^*)\right)\\
& = \wedge^5 \delta_{(3),\R}^*\left(Z_{12}^* \wedge Z_{23}^* \wedge X_{(-1,1,0)}^*\wedge X_{(0,-1,1)}^* \wedge X_{(-1,0,1)}^*\right).
\end{align*}
Thus the scalar between the functionals is equal to $1$.
This completes the proof.
\end{proof}

As a consequence of Theorems \ref{T:adjoint} and \ref{T:Asai}, we obtain the following refinement of Theorem \ref{T:GHL 2} for $n=2,3$.

\begin{thm}\label{T:Asai 2}
\noindent
\begin{itemize}
\item[(1)] We have
\begin{align*}
&\sigma\left(\frac{L^{(\infty)}(1,\itPi,{\rm As}^+ \otimes \chi)}{|D_\F|^{1/2}\cdot (\sqrt{-1})^{d{\sf w}(\itPi)}\cdot  G(\chi)\cdot\pi^{3d+\sum_{v \in S_\infty}(\kappa_{1,v}-\kappa_{2,v})/2}\cdot  p^b(\itPi)}\right) \\
& = \frac{L^{(\infty)}(1,{}^\sigma\!\itPi,{\rm As}^+ \otimes {}^\sigma\!\chi)}{|D_\F|^{1/2}\cdot (\sqrt{-1})^{d{\sf w}(\itPi)}\cdot G({}^\sigma\!\chi)\cdot\pi^{3d+\sum_{v \in S_\infty}(\kappa_{1,v}-\kappa_{2,v})/2}\cdot  p^b({}^\sigma\!\itPi)}
\end{align*}
for all $\sigma \in {\rm Aut}(\C)$.
\item[(2)] We have
\begin{align*}
&\sigma\left(\frac{L^{(\infty)}(1, \itSigma,{\rm As}^- \otimes \eta)}{(\sqrt{-1})^{d{\sf w}(\itSigma)}\cdot G(\omega_{\E/\F})\cdot G(\eta)^3\cdot \pi^{6d+\sum_{v \in S_\infty}(\ell_{1,v}-\ell_{3,v})}\cdot p^b(\itSigma)} \right)\\
& = \frac{L^{(\infty)}(1, \itSigma,{\rm As}^- \otimes \eta)}{(\sqrt{-1})^{d{\sf w}(\itSigma)}\cdot G(\omega_{\E/\F})\cdot G({}^\sigma\!\eta)^3\cdot\pi^{6d+\sum_{v \in S_\infty}(\ell_{1,v}-\ell_{3,v})}\cdot p^b({}^\sigma\!\itSigma)}
\end{align*}
for all $\sigma \in {\rm Aut}(\C)$.
\end{itemize}
\end{thm}

\begin{rmk}
Under the normalization of the generators in \cite[\S\,1.5.3]{GL2020}, similar result was proved by Grobner and Lin \cite[Theorem B]{GL2020} for $\GL_n(\A_\E) \times \GL_{n-1}(\A_\E)$ with additional assumption that $\itPi^\vee = \itPi^c$, $\itSigma^\vee = \itSigma^c$, and the Galois-equivariance was proved for $\sigma \in {\rm Aut}(\C/\E^{\Gal})$.
When $\lambda$ or $\mu$ is not sufficiently regular, they also need to assume certain non-vanishing hypotheses on central critical values.
\end{rmk}

Finally, we have the following immediate corollary of Theorems \ref{T:RS} and \ref{T:Asai 2}.

\begin{corollary}\label{T:GGP 2}
Assume $\ell_{1,v} > -\kappa_{2,v} > \ell_{2,v} > -\kappa_{1,v} > \ell_{3,v}$ for all $v \in S_\infty$. Let $m+\tfrac{1}{2}$ be critical for $L(s,\itSigma \times \itPi)$. We have
\begingroup
\small
\begin{align*}
&\sigma\left(\frac{L^{(\infty)}(m+\tfrac{1}{2},\itSigma \times \itPi)}{|D_\F|^{1/2}\cdot(2\pi\sqrt{-1})^{3d(2m-2+{\sf w}(\itSigma)+{\sf w}(\itPi))}\cdot G(\eta^3\cdot \chi)^{-1}\cdot G(\omega_\itPi)\cdot L^{(\infty)}(1,\itSigma,{\rm As}^-\otimes\eta)\cdot L^{(\infty)}(1,\itPi,{\rm As}^+\otimes\chi)}\right)\\
& = \frac{L^{(\infty)}(m+\tfrac{1}{2},{}^\sigma\!\itSigma \times {}^\sigma\!\itPi)}{|D_\F|^{1/2}\cdot(2\pi\sqrt{-1})^{3d(2m-2+{\sf w}(\itSigma)+{\sf w}(\itPi))}\cdot G({}^\sigma\!\eta^3\cdot {}^\sigma\!\chi)^{-1}\cdot G({}^\sigma\!\omega_\itPi)\cdot L^{(\infty)}(1,{}^\sigma\!\itSigma,{\rm As}^-\otimes{}^\sigma\!\eta)\cdot L^{(\infty)}(1,{}^\sigma\!\itPi,{\rm As}^+\otimes{}^\sigma\!\chi)}
\end{align*}
\endgroup
for all $\sigma \in {\rm Aut}(\C)$.
\end{corollary}


\subsection*{Acknowledgement}

The author would like to thank Ming-Lun Hsieh, Atsushi Ichino, and Kazuki Morimoto for suggestions and helpful conversations. Thanks are also due to the referee for his/her careful reading and valuable comments. The author was partially supported by the MOST grant 108-2628-M-001-009-MY4 and 110-2628-M-001-004 -.

\appendix
\section{Period relation for $\GL_n$ under duality}\label{S:appendix}

The purpose of this appendix is to prove Proposition \ref{P:period relation duality}, which is an automorphic analogue of the motivic period relation in \cite[Proposition 5.1]{Deligne1979}. By the global functional equation, algebraicity of critical $L$-values holds in the right-half critical range for a cohomological irreducible cuspidal automorphic representation if and only if it holds in the left-half critical range for its contragredient representation. 
We establish period relation for the conjectural periods under duality. We use the result to pass from the right-half critical range to the left-half critical range in the proof of Theorem \ref{T:Deligne}.

Fix an integer $n \in \Z_{\geq 2}$.
Let $\itPi$ be a cohomological irreducible cuspidal automorphic representation of $\GL_n(\A_\F)$ with central character $\omega_\itPi$. We have $|\omega_\itPi|=|\mbox{ }|_{\A_\F}^{n{\sf w}/2}$ for some ${\sf w} \in \Z$. Note that $\sf w$ must be even if $n$ is odd.
Let 
\[
L(s,\itPi)
\]
be the standard $L$-function of $\itPi$. We denote by $L^{(\infty)}(s,\itPi)$ the $L$-function obtained by excluding the archimedean $L$-factors. A critical point for $L(s,\itPi)$ is a half-integer $m+\tfrac{n-1}{2} \in \Z+\tfrac{n-1}{2}$ which is not a pole of the archimedean local factors $L(s,\itPi_v)$ and $L(1-s,\itPi_v^\vee)$ for all $v \in S_\infty$. We denote by ${\rm Crit}(\itPi)$ the set of critical points. The explicit description of ${\rm Crit}(\itPi)$ is given as follows: Write $n=2r$ (resp.\,$n=2r+1$) for some $r \in \Z_{\geq 1}$ if $n$ is even (resp.\,odd). 
For each $v \in S_\infty$, we have
\begin{align*}
\itPi_v = \begin{cases}
{\rm Ind}_{P_{2,\cdots,2}(\R)}^{\GL_n(\R)}(D_{\kappa_{1,v}}\boxtimes\cdots\boxtimes D_{\kappa_{r,v}})\otimes |\mbox{ }|_\R^{{\sf w}/2} & \mbox{ if $n=2r$},\\
{\rm Ind}_{P_{2,\cdots,2,1}(\R)}^{\GL_n(\R)}(D_{\ell_{1,v}}\boxtimes\cdots\boxtimes D_{\ell_{r,v}}\boxtimes {\rm sgn}^\delta)\otimes |\mbox{ }|_\R^{{\sf w}/2} & \mbox{ if $n=2r+1$},
\end{cases}
\end{align*}
for some $\kappa_{1,v} > \cdots > \kappa_{r,v} \geq 2$ such that $\kappa_{1,v} \equiv \cdots \kappa_{r,v} \equiv {\sf w}\,({\rm mod}\,2)$ if $n$ is even, and some odd integers $\ell_{1,v}>\cdots>\ell_{r,v}\geq 3$ and some $\delta\in\{0,1\}$ independent of $v$ if $n$ is odd.
Here $P_{2,\cdots,2}$ and $P_{2,\cdots,2,1}$ are the standard parabolic subgroups of $\GL_{2r}$ and $\GL_{2r+1}$ with Levi parts $\GL_2 \times\cdots \times\GL_2$ and $\GL_2 \times\cdots \times\GL_2\times\GL_1$, respectively.
Let ${\rm Crit}^+(\itPi)$ and ${\rm Crit}^-(\itPi)$ be the sets of right-half and left-half critical points, respectively, defined by
\begin{align*}
{\rm Crit}^+(\itPi)& = \begin{cases}
\left\{m+\tfrac{1}{2} \in \Z+\tfrac{1}{2}\,\left\vert\,-\tfrac{\sf w}{2} \leq m \leq \tfrac{\min_{v \in S_\infty}\{\kappa_{r,v}\}-{\sf w}}{2}-1\right.\right\} & \mbox{ if $n=2r$},\\
\left\{m\in\Z\,\vert\, 1-\tfrac{\sf w}{2} \leq m \leq \tfrac{\min_{v \in S_\infty}\{\ell_{r,v}\}-1}{2}-\tfrac{\sf w}{2},\,(-1)^{m+{\sf w}/2} = (-1)^\delta \right\} & \mbox{ if $n=2r+1$},
\end{cases}\\
{\rm Crit}^-(\itPi)& = \begin{cases}
\left\{-m-{\sf w}+\tfrac{1}{2}\,\left\vert\,m+\tfrac{1}{2} \in {\rm Crit}^+(\itPi)\right.\right\} & \mbox{ if $n=2r$},\\
\left\{1-m-{\sf w}\,\left\vert\,m \in {\rm Crit}^+(\itPi)\right.\right\} & \mbox{ if $n=2r+1$}.
\end{cases}
\end{align*}
Then ${\rm Crit}(\itPi)$ is the union of ${\rm Crit}^+(\itPi)$ and ${\rm Crit}^-(\itPi)$.
The union is not disjoint if and only if both $n$ and ${\sf w}$ are even. In this case, the intersection contains a unique critical point $\tfrac{1-{\sf w}}{2}$, which is the central critical point.
We have the following conjecture on the algebraicity of the critical $L$-values.
\begin{conj}[Deligne \cite{Deligne1979}, Clozel \cite{Clozel1990}]\label{C:DC}
There exist non-zero complex numbers $p({}^\sigma\!\itPi,\varepsilon)$ defined for each $\sigma \in {\rm Aut}(\C)$ and $\varepsilon=\pm1$ such that
\begin{align*}
\sigma\left(\frac{L^{(\infty)}(m+\tfrac{1}{2},\itPi)}{(2\pi\sqrt{-1})^{rdm}\cdot p(\itPi,(-1)^m)}\right) = \frac{L^{(\infty)}(m+\tfrac{1}{2},{}^\sigma\!\itPi)}{(2\pi\sqrt{-1})^{rdm}\cdot p({}^\sigma\!\itPi,(-1)^m)}
\end{align*}
for all $m+\tfrac{1}{2} \in {\rm Crit}(\itPi)$ if $n=2r$, and
\begin{align*}
\sigma\left(\frac{L^{(\infty)}(m-\tfrac{\sf w}{2},\itPi)}{(2\pi\sqrt{-1})^{\left(r+\tfrac{1\pm1}{2}\right)dm}\cdot p(\itPi,\pm)}\right) = \frac{L^{(\infty)}(m-\tfrac{\sf w}{2},{}^\sigma\!\itPi)}{(2\pi\sqrt{-1})^{\left(r+\tfrac{1\pm1}{2}\right)dm}\cdot p({}^\sigma\!\itPi,\pm)}
\end{align*}
for all $m -\tfrac{\sf w}{2} \in {\rm Crit}^\pm(\itPi)$ if $n=2r+1$.
\end{conj}

For each place $v$ of $\F$, let $\varepsilon(s,\itPi_v,\psi_{\F_v})$ and $\gamma(s,\itPi_v,\psi_{\F_v})$ be the $\varepsilon$-factor and $\gamma$-factor of $\itPi_v$ with respect to $\psi_{\F_v}$ defined by the local zeta integrals of Godement and Jacquet \cite{GJ1972}.

\begin{lemma}\label{L:critical gamma factor}
We have
\begin{align*}
\begin{cases}
\prod_{v \in S_\infty}\gamma(-m+\tfrac{1}{2},\itPi_v,\psi_{\F_v}) \in (2\pi\sqrt{-1})^{-ndm+rd{\sf w}}\cdot \Q^\times
& \mbox{ if $n=2r$ and $m+\tfrac{1}{2} \in {\rm Crit}^-(\itPi^\vee)$},\\
\prod_{v \in S_\infty}\gamma(1-m-\tfrac{\sf w}{2},\itPi_v,\psi_{\F_v}) \in (2\pi\sqrt{-1})^{-ndm+rd+d}\cdot \Q^\times
& \mbox{ if $n=2r+1$ and $m+\tfrac{\sf w}{2} \in {\rm Crit}^-(\itPi^\vee)$}.
\end{cases}
\end{align*}
\end{lemma}

\begin{proof}
Let $v \in S_\infty$. Suppose $n=2r$ and $m+\tfrac{1}{2} \in {\rm Crit}^-(\itPi^\vee)$, then we have
\[
\gamma(-m+\tfrac{1}{2},D_{\kappa_{i,v}}\otimes |\mbox{ }|_\R^{{\sf w}/2},\psi_{\F_v}) = (\sqrt{-1})^{\kappa_{i,v}}\cdot\frac{\Gamma_\C(m+\tfrac{\kappa_{i,v}-{\sf w}}{2})}{\Gamma_\C(-m+\tfrac{\kappa_{i,v}+{\sf w}}{2})} \in (2\pi\sqrt{-1})^{-2m+{\sf w}}\cdot\Q^\times
\]
for $1 \leq i \leq r$. 
Suppose $n=2r+1$ and $m+\tfrac{\sf w}{2} \in {\rm Crit}^-(\itPi^\vee)$, then we have
\begin{align*}
\gamma(1-m-\tfrac{\sf w}{2},{\rm sgn}^\delta|\mbox{ }|_\R^{{\sf w}/2},\psi_{\F_v})& = (\sqrt{-1})^{\delta}\cdot \frac{\Gamma_\R(m+\delta)}{\Gamma_\R(1-m+\delta)} \in (2\pi\sqrt{-1})^{-m+1}\cdot\Q^\times,\\
\gamma(1-m-\tfrac{\sf w}{2},D_{\ell_{i,v}}\otimes|\mbox{ }|_\R^{{\sf w}/2},\psi_{\F_v})& = (\sqrt{-1})^{\ell_{i,v}}\cdot \frac{\Gamma_\C(m+\tfrac{\ell_{i,v}-1}{2})}{\Gamma_\C(1-m+\tfrac{\ell_{i,v}-1}{2})} \in (2\pi\sqrt{-1})^{-2m+1}\cdot\Q^\times
\end{align*}
for $1 \leq i \leq r$.  
This completes the proof.
\end{proof}

\begin{lemma}\label{L:Galois-equiv. root number}
We have
\begin{align*}
\sigma \left( \frac{\prod_{v \nmid \infty}\varepsilon(m+\tfrac{n-1}{2},\itPi_v,\psi_{\F_v})}{|D_\F|^{n/2}\cdot G(\omega_\itPi)}\right) = \frac{\prod_{v \nmid \infty}\varepsilon(m+\tfrac{n-1}{2},{}^\sigma\!\itPi_v,\psi_{\F_v})}{|D_\F|^{n/2}\cdot G({}^\sigma\!\omega_\itPi)}
\end{align*}
for all $\sigma \in {\rm Aut}(\C)$ and $m\in\Z$.
\end{lemma}

\begin{proof}
For a rational function $Q$ and $\sigma \in {\rm Aut}(\C)$, let ${}^\sigma\!Q$ be the rational function obtained by acting $\sigma$ on its coefficients. Let $v$ be a finite place of $\F$ lying over a rational prime $p$. Let $D_{\F_v}$ be the discriminant of $\F_v/ \Q_p$, $q_v$ the cardinality of the residue field of $\F_v$, and $|\mbox{ }|_{\F_v}$ the absolute value on $\F_v$ normalized so that it takes value $q_v^{-1}$ on a uniformizer of $\frak{o}_{\F_v}$.
For a Schwartz function $\Phi$ on the $n$ by $n$ matrices ring ${\rm M}_n(\F_v)$ and a matrix coefficient $f$ of $\itPi_v$, we have the local zeta integral integral
\[
Z_v(\Phi,f,s) = \int_{\GL_n(\F_v)}\Phi(g)f(g)|\det(g)|_{\F_v}^{s+(n-1)/2}\,dg
\]
defined by Godement and Jacquet \cite{GJ1972}.
The integral converges absolutely for ${\rm Re}(s)$ sufficiently large and defines a rational function in $q_v^{-s}$. Moreover, by \cite[Theorem 3.3]{GJ1972}, the ratio $\frac{Z_v(\Phi,f,s)}{L(s,\itPi_v)}$ belongs to $\C[q_v^{-s},q_v^s]$ and we have the local functional equational 
\[
\frac{Z_v(\widehat{\Phi},f^\vee,1-s)}{L(1-s,\itPi_v^\vee)} = \varepsilon(s,\itPi_v,\psi_{\F_v})\cdot \frac{Z_v(\Phi,f,s)}{L(s,\itPi_v)}.
\]
Here $f^\vee(g) = f(g^{-1})$ is a matrix coefficient of $\itPi_v^\vee$ and $\widehat{\Phi}$ is the Fourier transform of $\Phi$ with respect to $\psi_{\F_v}$ defined by
\[
\widehat{\Phi}(x) = \int_{{\rm M}_n(\F_v)}\Phi(y)\psi_{\F_v}({\rm tr}(xy))\,d_{\psi_{\F_v}}y
\]
with ${\rm vol}({\rm M}_n(\frak{o}_{\F_v}),d_{\psi_{\F_v}}y)= |D_{\F_v}|_{\Q_p}^{-n^2/2}$. Note that the measure $d_{\psi_{\F_v}}y$ is the self-dual Haar measure with respect to $\psi_{\F_v}\circ {\rm tr}$. 
Let $\sigma \in {\rm Aut}(\C)$. Define ${}^\sigma\!\Phi$ and ${}^\sigma\!f$ by ${}^\sigma\!\Phi(x) = \sigma(\Phi(x))$ and ${}^\sigma\!f(g) = \sigma(f(g))$.
It is easy to see that ${}^\sigma\!f$ is a matrix coefficient of ${}^\sigma\!\itPi_v^\vee$. We have the following identity for rational functions in $q_v^{-s}$, which is established implicitly in the proof of \cite[Lemme 4.6]{Clozel1990}:
\begin{align}\label{E:root proof 1}
{}^\sigma\!Z_v(\Phi,f,s+\tfrac{n-1}{2}) = Z_v({}^\sigma\!\Phi,{}^\sigma\!f,s+\tfrac{n-1}{2}).
\end{align}
We have
\begin{align*}
{}^\sigma\!(\widehat{\Phi})(x) & = \frac{\sigma(|D_{\F_v}|_{\Q_p}^{-n/2})}{|D_{\F_v}|_{\Q_p}^{-n/2}}\cdot\int_{{\rm M}_n(\F_v)}{}^\sigma\!\Phi(y)\sigma(\psi_{\F_v}({\rm tr}(xy)))\,d_{\psi_{\F_v}}y\\
& = \frac{\sigma(|D_{\F_v}|_{\Q_p}^{-n/2})}{|D_{\F_v}|_{\Q_p}^{-n/2}}\cdot \widehat{{}^\sigma\!\Phi}(u_{\sigma,p}x).
\end{align*}
Here $u_{\sigma,p} \in \Z_p^\times\subset \F_v^\times$ is the unique element such that $\sigma(\psi_{\Q_p}(x)) = \psi_{\Q_p}(u_{\sigma,p}x)$ for $x \in \Q_p$.
Therefore, a simple change of variable implies that
\begin{align}\label{E:root proof 2}
Z_v({}^\sigma\!(\widehat{\Phi}),{}^\sigma\!f^\vee,s) = \frac{\sigma(|D_{\F_v}|_{\Q_p}^{-n/2})}{|D_{\F_v}|_{\Q_p}^{-n/2}}\cdot {}^\sigma\!\omega_{\itPi_v}(u_{\sigma,p})\cdot Z_v(\widehat{{}^\sigma\!\Phi},{}^\sigma\!f^\vee,s).
\end{align}
Now we apply $\sigma$ to both sides of the local functional equation and using (\ref{E:root proof 1}) and (\ref{E:root proof 2}). We conclude that
\begin{align*}
\frac{Z_v(\widehat{{}^\sigma\!\Phi},{}^\sigma\!f^\vee,-s+\tfrac{3-n}{2})}{L(-s+\tfrac{3-n}{2},{}^\sigma\!\itPi_v^\vee)} = \frac{\sigma(|D_{\F_v}|_{\Q_p}^{n/2})}{|D_{\F_v}|_{\Q_p}^{n/2}}\cdot {}^\sigma\!\omega_{\itPi_v}(u_{\sigma,p})^{-1}\cdot{}^\sigma\!\varepsilon(s+\tfrac{n-1}{2},\itPi_v,\psi_{\F_v})\cdot \frac{Z_v({}^\sigma\!\Phi,{}^\sigma\!f,s+\tfrac{n-1}{2})}{L(s+\tfrac{n-1}{2},{}^\sigma\!\itPi_v)}.
\end{align*}
In other words, we have
\begin{align*}
{}^\sigma\!\varepsilon(s+\tfrac{n-1}{2},\itPi_v,\psi_{\F_v}) = \frac{\sigma(|D_{\F_v}|_{\Q_p}^{-n/2})}{|D_{\F_v}|_{\Q_p}^{-n/2}}\cdot {}^\sigma\!\omega_{\itPi_v}(u_{\sigma,p})\cdot\varepsilon(s+\tfrac{n-1}{2},{}^\sigma\!\itPi_v,\psi_{\F_v}).
\end{align*}
Since the $\varepsilon$-factors are monomials in $q_v^{-s}$, we deduce that 
\begin{align*}
\sigma\left(\varepsilon(m+\tfrac{n-1}{2},\itPi_v,\psi_{\F_v})\right) = \frac{\sigma(|D_{\F_v}|_{\Q_p}^{-n/2})}{|D_{\F_v}|_{\Q_p}^{-n/2}}\cdot {}^\sigma\!\omega_{\itPi_v}(u_{\sigma,p})\cdot\varepsilon(m+\tfrac{n-1}{2},{}^\sigma\!\itPi_v,\psi_{\F_v})
\end{align*}
for all $m \in \Z$.
Finally, note that $\varepsilon(s,\itPi_v,\psi_{\F_v})=1$ for all but finitely many $v$ and 
\[
\prod_{p}\prod_{v \mid p}|D_{\F_v}|_{\Q_p} = |D_\F|^{-1},\quad \frac{\sigma(G(\omega_\itPi))}{G({}^\sigma\!\omega_\itPi)} = \prod_p \prod_{v \mid p}{}^\sigma\!\omega_{\itPi_v}(u_{\sigma,p}).
\]
This completes the proof.
\end{proof}

\begin{prop}\label{P:period relation duality}
Conjecture \ref{C:DC} holds for $\itPi$ in the right-half critical set ${\rm Crit}^+(\itPi)$ if and only if it holds for $\itPi^\vee$ in the left-half critical set ${\rm Crit}^-(\itPi^\vee)$. In this case, the periods can be chosen so that
\begin{align*}
\begin{cases}
p({}^\sigma\!\itPi^\vee,\pm) = (2\pi\sqrt{-1})^{-rd{\sf w}}\cdot G({}^\sigma\!\omega_\itPi)^{-1}\cdot p({}^\sigma\!\itPi,\pm)
& \mbox{ if $n=2r$},\\
p({}^\sigma\!\itPi^\vee,-) = |D_\F|^{1/2}\cdot G({}^\sigma\!\omega_\itPi)^{-1}\cdot p({}^\sigma\!\itPi,+) & \mbox{ if $n=2r+1$}
\end{cases}
\end{align*}
for $\sigma \in {\rm Aut}(\C)$.
\end{prop}

\begin{proof}
The assertion follows from Lemmas \ref{L:critical gamma factor}, \ref{L:Galois-equiv. root number}, and the global functional equation
\[
L^{(\infty)}(s,\itPi) = \prod_{v \nmid \infty}\varepsilon(s,\itPi_v,\psi_{\F_v})\cdot \prod_{v \in S_\infty}\gamma(s,\itPi_v,\psi_{\F_v})\cdot L^{(\infty)}(1-s,\itPi^\vee).
\]
\end{proof}


\begin{thebibliography}{JPSS79}

\bibitem[BP21]{Raphael2018}
R.~Beuzart-Plessis.
\newblock {Archimedean theory and $\epsilon$-factors for the Asai
  Rankin-Selberg integrals}.
\newblock In {\em {Relative trace formulas}}, Simons Symposia, pages 1--50.
  Springer, 2021.

\bibitem[BR17]{BR2017}
B.~Balasubramanyam and A.~Raghuram.
\newblock {Special values of adjoint $L$-functions and congruences for
  automorphic forms on ${\rm GL}(n)$ over a number field}.
\newblock {\em Amer. J. Math.}, 139(3):641--679, 2017.

\bibitem[CCI20]{CCI2020}
S.-Y. Chen, Y.~Cheng, and I.~Ishikawa.
\newblock {Gamma factors for Asai representations of ${\rm GL}_2$}.
\newblock {\em J. Number Theory}, 209:83--146, 2020.

\bibitem[Che20]{Cheng2020}
Y.~Cheng.
\newblock {Special value formula for the twisted triple product $L$-function
  and an application to the restricted $L^2$-norm problem}.
\newblock {\em Forum Math.}, 2020.
\newblock DOI:10.1515/forum-2018-0292.

\bibitem[Che22a]{Chen2021}
S.-Y. Chen.
\newblock {Algebraicity of critical values of triple product $L$-functions in
  the balanced case}.
\newblock {\em Pacific J. Math.}, 2022.
\newblock accepted, arXiv:2108.02111.

\bibitem[Che22b]{Chen2020}
S.-Y. Chen.
\newblock {Algebraicity of the near central non-critical values of symmetric
  fourth $L$-functions for Hilbert modular forms}.
\newblock {\em J. Number Theory}, 231:269--315, 2022.

\bibitem[CHT08]{CHT2008}
L.~Clozel, M.~Harris, and R.~Taylor.
\newblock {Automorphy for some $l$-adic lifts of automorphic mod $l$ Galois
  representations}.
\newblock {\em Publ. Math. Inst. Hautes \'Etudes Sci.}, 108:1--181, 2008.

\bibitem[Clo90]{Clozel1990}
L.~Clozel.
\newblock {Motifs et Formes Automorphes: Applications du Principe de
  Fonctorialit\'e}.
\newblock In {\em Automorphic Forms, Shimura Varieties, and L-functions, Vol.
  I}, Perspectives in Mathematics, pages 77--159, 1990.

\bibitem[Del79]{Deligne1979}
P.~Deligne.
\newblock {Valeurs de fonctions $L$ et périodes d’intégrales}.
\newblock In {\em Automorphic Forms, Representations and $L$-Functions},
  volume~33, pages 313--346. Proceedings of Symposia in Pure Mathematics, 1979.
\newblock Part 2.

\bibitem[Fli88]{Flicker1988}
Y.~Z. Flicker.
\newblock {Twisted ttensor and Euler products}.
\newblock {\em Bull. Soc. Math. France}, 116:295--313, 1988.

\bibitem[FZ95]{FZ1995}
Y.~Z. Flicker and D.~Zinoviev.
\newblock {On poles of twisted tensor $L$-functions}.
\newblock {\em Proc. Japan Acad. Ser. A Math. Sci.}, 71(6):114--116, 1995.

\bibitem[GH93]{GH1993}
P.~Garrett and M.~Harris.
\newblock {Special values of triple product $L$-functions}.
\newblock {\em Amer. J. Math.}, 115(1):161--240, 1993.

\bibitem[GHL16]{GHL2016}
H.~Grobner, M.~Harris, and E.~Lapid.
\newblock {Whittaker rational structures and special values of the Asai
  $L$-function}.
\newblock {\em Contemp. Math.}, 664:119--134, 2016.

\bibitem[GJ72]{GJ1972}
R.~Godement and H.~Jacquet.
\newblock {\em {Zeta functions of simple algebras}}, volume 260 of {\em Lecture
  Notes in Mathematics}.
\newblock Springer, 1972.

\bibitem[GL21]{GL2020}
H.~Grobner and J.~Lin.
\newblock {Special values of $L$-functions and the refined Gan-Gross-Prasad
  conjecture}.
\newblock {\em Amer. J. Math.}, 143(3):859--937, 2021.

\bibitem[Gro18]{Grobner2018}
H.~Grobner.
\newblock {Rationality results for the exterior and the symmetric square
  $L$-function (with an appendix by Nadir Matringe)}.
\newblock {\em Math. Ann.}, 370:1639--1679, 2018.

\bibitem[GT11]{GT2011}
W.~T. Gan and S.~Takeda.
\newblock {The local Langlands conjecture for ${\rm GSp}(4)$}.
\newblock {\em Ann. of Math.}, 173:1841--1882, 2011.

\bibitem[GT19]{GT2019}
T.~Gee and O.~Ta\"{i}bi.
\newblock {Arthur's multiplicity formula for ${\rm GSp}_4$ and restriction to
  ${\rm Sp}_4$}.
\newblock {\em J. \'{E}c. polytech. Math.}, 6:469--535, 2019.

\bibitem[GW09]{GW2009}
R.~Goodman and N.~Wallach.
\newblock {\em {Symmetry, representations, and invariants}}, volume 255 of {\em
  Graduate Texts in Mathematics}.
\newblock Springer, 2009.

\bibitem[HIM12]{HIM2012}
M.~Hirano, T.~Ishii, and T.~Miyazaki.
\newblock {Archimedean Whittaker functions on ${\rm GL}(3)$}.
\newblock In {\em Geometry and analysis of automorphic forms of several
  variables}, volume~7 of {\em Series on Number Theory and its applications}.
  World Scientific, 2012.

\bibitem[HIM16]{HIM2016}
M.~Hirano, T.~Ishii, and T.~Miyazaki.
\newblock {The archimedean zeta integrals for ${\rm GL}(3) \times {\rm
  GL}(2)$}.
\newblock {\em Proc. Japan Acad. Ser. A Math. Sci.}, 92(2), 2016.

\bibitem[HO09]{HO2009}
M.~Hirano and T.~Oda.
\newblock {Calculus of principal series Whittaker functions on ${\rm
  GL}(3,\mathbb{C})$}.
\newblock {\em J. Funct. Anal.}, 256:2222--2267, 2009.

\bibitem[HR20]{HR2020}
G.~Harder and A.~Raghuram.
\newblock {\em {Eisenstein cohomology for ${\rm GL}_N$ and the special values
  of Rankin--Selberg $L$-functions}}.
\newblock Annals of Mathematics Studies. Princeton University Press, 2020.

\bibitem[Ich05]{Ichino2005}
A.~Ichino.
\newblock {Pullbacks of Saito-Kurokawa lifts}.
\newblock {\em Invent. Math.}, 162:551--647, 2005.

\bibitem[IK14]{IK2014}
T.~Ibukiyama and H.~Katsurada.
\newblock {Exact critical values of the symmetric fourth $L$-function and
  vector valued Siegel modular forms}.
\newblock {\em J. Math. Soc. Japan}, 66(1):139--160, 2014.

\bibitem[Im91]{Im1991}
J.~Im.
\newblock {Special values of Dirichlet series attached to Hilbert modular
  forms}.
\newblock {\em Amer. J. Math.}, 113(6):975--1017, 1991.

\bibitem[Jac09]{Jacquet2009}
H.~Jacquet.
\newblock {Archimedean Rankin-Selberg integrals}.
\newblock In {\em {Automorphic forms and $L$-functions II. Local aspects}},
  volume 489 of {\em Contemp. Math.}, pages 57--172, 2009.

\bibitem[Jan19]{Januszewski2019}
F.~Januszewski.
\newblock {On period relations for automorphic $L$-functions I}.
\newblock {\em Trans. Amer. Math. Soc.}, 371(9):6547--6580, 2019.

\bibitem[JPSS79]{JPSS1979}
H.~Jacquet, I.~I. Piatetski-Shapiro, and J.~A. Shalika.
\newblock {Automorphic forms on ${\rm GL}(3)$, I}.
\newblock {\em Ann. of Math.}, 109(1):169--212, 1979.

\bibitem[JS76]{JS1976}
H.~Jacquet and J.~A. Shalika.
\newblock {A non-vanishing theorem for zeta functions of ${\rm GL}_n$}.
\newblock {\em Invent. Math.}, 38:1--16, 1976.

\bibitem[JS81]{JS1981}
H.~Jacquet and J.~A. Shalika.
\newblock {On Euler products and the classification of automorphic forms I}.
\newblock {\em Amer. J. Math.}, 103(3):499--558, 1981.

\bibitem[JS90]{JS1990}
H.~Jacquet and J.~A. Shalika.
\newblock {Rankin--Selberg convolutions: Archimedean theory}.
\newblock In {\em Festschrift in honor of I. I. Piatetski-Shapiro on the
  occasion of his sixtieth birthday, Part I}, pages 125--207. Weizmann Science
  Press of Israel, 1990.

\bibitem[Kab04]{Kable2004}
A.~C. Kable.
\newblock {Asai $L$-functions and Jacquet's conjecture}.
\newblock {\em Amer. J. Math.}, 126(4):789--820, 2004.

\bibitem[Kem15]{Kemarsky2015}
A.~Kemarsky.
\newblock {Distinguished representations of ${\rm GL}_n(\mathbb{C})$}.
\newblock {\em Israel J. Math.}, 207:435--448, 2015.

\bibitem[Kim03]{Kim2003}
H.~H. Kim.
\newblock {Functoriality for the exterior square of ${\rm GL}_4$ and the
  symmetric fourth of ${\rm GL}(2)$}.
\newblock {\em J. Amer. Math. Soc.}, 16(1):139--183, 2003.

\bibitem[KS02a]{KS2002b}
H.~H. Kim and F.~Shahidi.
\newblock Cuspidality of symmetric powers with applications.
\newblock {\em Duke Math. J.}, 112(1):177--197, 2002.

\bibitem[KS02b]{KS2002}
H.~H. Kim and F.~Shahidi.
\newblock {Functorial products for $\GL_2 \times \GL_3$ and the symmetric cube
  for $\GL_2$}.
\newblock {\em Ann. of Math.}, 155(2):837--893, 2002.

\bibitem[Liu21]{Liu2019b}
Z.~Liu.
\newblock {The doubling Archimedean zeta integrals for $p$-adic interpolation}.
\newblock {\em Math. Res. Lett.}, 28(1):145--173, 2021.

\bibitem[Mat11]{Matringe2011}
N.~Matringe.
\newblock {Distinguished generic representations of ${\rm GL}(n)$ over $p$-adic
  fields}.
\newblock {\em Int. Math. Res. Not.}, (1):74--95, 2011.

\bibitem[Mok15]{Mok2015}
C.-P. Mok.
\newblock {Endoscopic classification of representations of quasi-split unitary
  groups}.
\newblock {\em Mem. Amer. Math. Soc.}, 235(1108), 2015.

\bibitem[Mor21]{Morimoto2021}
K.~Morimoto.
\newblock {On algebraicity of special values of symmetric $4$-th and $6$-th
  power $L$-functions for ${\rm GL}(2)$}.
\newblock {\em Math. Z.}, 299:1331--1350, 2021.

\bibitem[Off11]{Offen2011}
O.~Offen.
\newblock {On local root numbers and distinction}.
\newblock {\em J. Reine Angew. Math.}, 652:165--205, 2011.

\bibitem[PSS20]{PSS2020}
A.~Pitale, A.~Saha, and R.~Schmidt.
\newblock {On the standard $L$-function for ${\rm GSp}_{2n} \times {\rm GL}_1$
  and algebraicity of symmetric fourth $L$-values for ${\rm GL}_2$}.
\newblock {\em Ann. Math. Qu\'e.}, 2020.
\newblock DOI:10.1007/s40316-020-00134-6.

\bibitem[Rag16]{Raghuram2016}
A.~Raghuram.
\newblock {Critical values for Rankin-Selberg $L$-functions for ${\rm GL}_n
  \times {\rm GL}_{n-1}$ and the symmetric cube $L$-functions for ${\rm
  GL}_2$.}
\newblock {\em Forum Math.}, 28(3):457--489, 2016.

\bibitem[Rag20]{Raghuram2020}
A.~Raghuram.
\newblock {Eisenstein cohomology for ${\rm GL}_N$ and the special values of
  Rankin--Selberg $L$-functions - II}.
\newblock 2020.
\newblock Preprint.

\bibitem[RS07a]{RS2007c}
A.~Raghuram and F.~Shahidi.
\newblock Raghuram shahidi - functoriality and special values of $l$-vunctions.
\newblock In {\em {Eisenstein Series and Applications}}, volume 258 of {\em
  {Progress in Mathematics}}, pages 271--293. Birkhauser, 2007.

\bibitem[RS07b]{RS2007b}
D.~Ramakrishnan and F.~Shahidi.
\newblock {Siegel modular forms of genus $2$ attached to elliptic curves}.
\newblock {\em Math. Res. Lett.}, 14:315--332, 2007.

\bibitem[RS08]{RS2008}
A.~Raghuram and F.~Shahidi.
\newblock {On certain period relations for cusp forms on ${\rm GL}_n$}.
\newblock {\em Int. Math. Res. Not.}, 2008.
\newblock doi:10.1093/imrn/rnn077.

\bibitem[RS18]{RS2018}
A.~Raghuram and M.~Sarnobat.
\newblock Cohomological representations and functorial transfer from classical
  groups.
\newblock In {\em Cohomology of arithmetic groups}, volume 245 of {\em
  {Springer Proceedings in Mathematics and Statistics}}, pages 157--176.
  Springer, 2018.

\bibitem[Sha81]{Shahidi1981}
F.~Shahidi.
\newblock {On certain $L$-functions}.
\newblock {\em Amer. J. Math.}, 103:297--355, 1981.

\bibitem[Shi78]{Shimura1978}
G.~Shimura.
\newblock {The special values of the zeta functions associated with Hilbert
  modular forms}.
\newblock {\em Duke Math. J.}, 45(3):635--679, 1978.

\bibitem[Stu89]{Sturm1989}
J.~Sturm.
\newblock {Evaluation of the symmetric square at the near center point}.
\newblock {\em Amer. J. Math.}, 111(4):585--598, 1989.

\bibitem[Sun17]{Sun2017}
B.~Sun.
\newblock {The nonvanishing hypothesis at infinity for Rankin-Selberg
  convolutions}.
\newblock {\em J. Amer. Math. Soc.}, 30(1):1--25, 2017.

\bibitem[Sun19]{Sun2019}
B.~Sun.
\newblock {Cohomologically induced distinguished representations and
  cohomological test vectors}.
\newblock {\em Duke Math. J.}, 168(1):85--126, 2019.

\bibitem[Tat79]{Tate1979}
J.~Tate.
\newblock Number theoretic background.
\newblock In {\em {Automorphic forms, representations, and L-functions}},
  volume~33, pages 3--26. Proceedings of Symposia in Pure Mathematics, 1979.
\newblock Part 2.

\bibitem[Wal84]{Wallach1984}
N.~Wallach.
\newblock On the constant term of a square integrable automorphic form.
\newblock In {\em Operator algebras and group representations}, volume~18 of
  {\em {Monographs and Studies in Mathematics}}, pages 227--237, 1984.

\bibitem[Was96]{Washingtonbook}
L.~C. Washington.
\newblock {\em Introduction to cyclotomic fields}, volume~85 of {\em Graduate
  Texts in Mathematics}.
\newblock Springer, 2 edition, 1996.

\bibitem[Zha14]{Zhang2014}
Wei Zhang.
\newblock {Automorphic period and the central value of Rankin--Selberg
  $L$-function}.
\newblock {\em J. Amer. Math. Soc.}, 27(2):541--612, 2014.

\end{thebibliography}

\end{document}